\magnification=\magstep1

\input amstex

\loadeusm 
\loadeusb 
\documentstyle{amsppt}
\voffset=-0.5in
\hfuzz=7pt

\baselineskip=20pt
\parskip=7pt plus 5pt
\parindent=15pt
\hoffset=-0.3in
 \pagewidth{7truein}
\def\average{{\mathchoice {\kern1ex\vcenter{\hrule height.4pt width 6pt
depth0pt} \kerQ-9.7pt} {\kern1ex\vcenter{\hrule height.4pt width 4.3pt
depth0pt}
\kerQ-7pt} {} {} }}

\leftheadtext{Y. Han and G. Lu}
\rightheadtext{Discrete Littlewood-Paley-Stein multi-parameter analysis }

\def\aaverage{{\mathchoice {\kern1ex\vcenter{\hrule height.4pt width 12pt
depth0pt} \kerQ-16pt} {\kern1ex\vcenter{\hrule height.4pt width 9pt depth0pt}
\kerQ-12pt} {} {} }}

\def\vbar{\mathchoice{\vrule height6.3ptdepth-.5ptwidth.8pt\kerQ-.8pt}
   {\vrule height6.3ptdepth-.5ptwidth.8pt\kerQ-.8pt}
   {\vrule height4.1ptdepth-.35ptwidth.6pt\kerQ-.6pt}
   {\vrule height4.1ptdepth-.25ptwidth.5pt\kerQ-.5pt}}
\def\fudge{\mathchoice{}{}{\mkern.5mu}{\mkern.8mu}}
\def\bbc#1#2{{\rm \mkern#2mu\vbar\mkerQ-#2mu#1}}
\def\bbb#1{{\rm I\mkerQ-4.5mu #1}}
\def\bba#1#2{{\rm #1\mkerQ-#2mu\fudge #1}}
\def\bb#1{{\count4=`#1 \advance\count4by-64 \ifcase\count4\or\bba A{11.5}\or
   \bbb B\or\bbc C{5}\or\bbb D\or\bbb E\or\bbb F \or\bbc G{5}\or\bbb H\or
   \bbb I\or\bbc J{3}\or\bbb K\or\bbb L \or\bbb M\or\bbb N\or\bbc O{5} \or
   \bbb P\or\bbc B{5}\or\bbb R\or\bba S{8}\or\bba T{10.5}\or\bbc U{5}\or
   \bba V{12}
\or\bba W{16.5}\or\bba X{11}\or\bba Y{11.7}\or\bba Z{7.5}\fi}}

\def\be{\begin{equation}}
\def\ee{\end{equation}}

\def\C0{C^\infty_0(\Omega)}
\def\K1{K_\epsilon(x,h)}

\def\phpf{H^p_F(R^n\times R^m)}
\def\prods{R^n\times R^m}
\def\bmof{BMO_F(R^n\times R^m)}

\def\ph1f{H^1_F(R^n\times R^m)}

\def\sf{{\Cal S}_F(R^n\times R^m)}

\NoBlackBoxes
\goodbreak
\topmatter
\title
Discrete Littlewood-Paley-Stein theory and multi-parameter
Hardy spaces associated with flag singular integrals
\endtitle
 \author
Yongsheng Han\\{\rm Department of Mathematics\\
Auburn University\\
Auburn, AL 36849, U.S.A\\
E-mail: hanyong\@mail.auburn.edu}\\
\smallskip
Guozhen Lu$^{*}$\\{\rm Department of Mathematics\\Wayne State University\\Detroit, MI 48202\\E-mail: gzlu\@math.wayne.edu }
\endauthor

\thanks {(*) Research was supported partly by
 the U.S. NSF Grant DMS\#0500853. \endgraf}\endthanks

\keywords {Flag singular integrals, Multiparameter Hardy spaces,  Discrete Calder\'on reproducing formulas, Discrete Littlewood-Paley-Stein analysis} \endkeywords

\thanks \noindent \endthanks
\goodbreak

\abstract

The main purpose of this paper is to develop a unified  approach of multi-parameter Hardy space theory using  the discrete
 Littlewood-Paley-Stein analysis in the setting of implicit multi-parameter structure.
It is motivated by the goal to establish and develop the Hardy space
theory for the flag singular integral operators studied by
Muller-Ricci-Stein [MRS] and Nagel-Ricci-Stein [NRS]. This approach
enables us to avoid the use of transference method of Coifman-Weiss
[CW] as often used in the $L^p$ theory for $p>1$ and establish the
Hardy spaces $H^p_F$ and its dual spaces associated with the flag
singular integral operators for all $0<p\leq 1$. We also prove the
boundedness of flag singular integral operators on $BMO_F$ and
$H^p_F$, and from $H^p_F$ to $L^p$  for all $0<p\le 1$ without using
the deep atomic decomposition. As a result, it bypasses the use of
Journe's type covering lemma in this implicit multi-parameter
structure. The method used here provides
 alternate approaches of   those developed by Chang-R. Fefferman [CF1-3], Chang [Ch], R.
Fefferman [F], Journe [J1-2], Pipher [P] in their important work in
pure product setting.     A Calder\'on-Zygmund decomposition and
interpolation theorem are also proved for the implicit
multi-parameter Hardy spaces.

\endabstract

\endtopmatter
\goodbreak

\document

\centerline{\bf Table of Contents}
\vskip .05cm
\noindent {\bf 1.  Introduction and statement of results}
\vskip .05cm
\noindent {\bf 2.  $L^p$ estimates for Littlewood-Paley-Stein square function: Proofs of Theorems 1.1 and 1.4}
\vskip .05cm
\noindent {\bf 3.  Test function spaces, almost orthogonality estimates, discrete Calder\'on reproducing formula: Proofs of Theorems 1.8 and 1.9}
\vskip .05cm
\noindent {\bf 4.  Discrete Littlewood-Paley-Stein square function and Hardy spaces, boundedness of flag singular integrals on Hardy spaces $H^p_F$, from $H^p_F$ to $L^p$: Proofs of Theorems 1.10 and 1.11}
\vskip .05cm
\noindent {\bf 5.  Duality of Hardy spaces $H^p_F$ and the boundedness of flag singular integrals on $BMO_F$ space: Proofs of Theorems 1.14, 1.16 and 1.18}
\vskip .05cm
\noindent{\bf 6. Calder\'on-Zgymund decomposition and interpolation on flag Hardy spaces $H^p_F$: Proofs of Theorems 1.19 and 1.20}
\vskip .05cm

\noindent {\bf References}
\vskip 0.2cm
\centerline {\bf 1.  Introduction and statement of results }
\vskip .2cm
The multi-parameter structures play a significant role
in Fourier analysis. On the one hand, the classical Calder\'on-Zygmund theory can be regarded as centering around
the Hardy-Littlewood maximal operator and certain singular integrals which commute with
the usual dilations on $R^n$, given by $\delta\cdot x=(\delta x_1,..., \delta x_n)$ for
$\delta >0$. On the other hand, if  we consider the multi-parameter dilations on $R^n$, given by
$\delta\cdot x =(\delta_1x_1,..., \delta_nx_n)$,
where $\delta =(\delta_1,...,\delta_n)\in R^n_+=(R_+)^n,$ then
these $n$-parameter dilations are naturally associated with the strong maximal function ([JMZ]),
given by
\vskip -0.3cm
$$M_s(f)(x)=\sup\limits_{x\in R}{{1}\over {\vert R\vert}}\int\limits_{R}\vert f(y)\vert dy,\tag1.1$$
\vskip -0.5cm
\noindent where the supremum is taken over the family of all rectangles with sides parallel to the axes.

This multi-parameter pure product theory has been developed by many authors over the past thirty
years or so.
For Calder\'on-Zygmund theory in this setting, one considers operators of the form
$Tf=K*f,$ where $K$ is homogeneous, that is, ${\delta_1}...{\delta_n}K(\delta\cdot
 x)=K(x),$ or, more generally, $K(x)$ satisfies the certain differential inequalities and cancellation
conditions such that ${\delta_1}...{\delta_n}K(\delta\cdot x)$ also satisfy  the same
bounds. This type of operators has been the subject of extensive investigations in the literature,  see
 for instances the  fundamental works of Gundy-Stein ([GS]), R. Fefferman and Stein [FS1], R. Fefferman ([F]), Chang and R. Fefferman ([CF1], [CF2], [CF3]), Journe ([J1], [J2]),  Pipher [P],   etc.

It is well-known that
there is a basic obstacle to the pure product Hardy space theory
and pure product $BMO$ space. Indeed, the role of cubes in the classical
atomic decomposition of $H^p(R^n)$ was replaced by arbitrary open
sets of finite measures in the product $H^p(R^n\times R^m).$ Suggested by a counterexample constructed
by L. Carleson [Car],
the very deep product $BMO(R^n\times R^m)$ and Hardy space $H^p(\prods)$ theory was developed
by Chang and R. Fefferman ([Ch],[CF3]). Because of the complicated nature of atoms in product space,  a certain
geometric lemma, namely Journe's covering lemma([J1], [J2] and [P]), played
an important role in the study of the boundedness of product
singular integrals on $H^p(R^n\times R^m)$ and $BMO(R^n\times R^m).$

While great progress has been made in the case of pure product structure for both $L^p$ and $H^p$
theory, multi-parameter
analysis has only been developed in recent years for the $L^p$ theory when the underlying multi
-parameter structure is not
 explicit, but implicit, such as the flag multi-parameter structure studied in [MRS] and [NRS].
 The main goal of this paper is to develop a theory of Hardy space in this setting.
 One of the main ideas of our program  is to develop a discrete version of Calder\'on reproducing formula associated
 with the given multiparameter structure, and thus prove a Plancherel-P\^olya type inequality in this setting.
 This discrete scheme of Littlewood-Paley-Stein analysis is
particularly useful in dealing with the Hardy spaces $H^p$  for
$0<p\le 1$.

 We now recall  two instances of implicit multiparameter structures which are of interest to us in this paper.
 We begin with reviewing one of these cases first. In the work of Muller-Ricci-Stein [MRS], by considering an implicit multi-parameter
 structure on Heisenberg(-type) groups,
the Marcinkiewicz multipliers on the Heisenberg
groups yield a new class of flag singular integrals.
To be more precise, let $m({\Cal L}, iT)$ be the Marcinkiewicz multiplier
operator, where ${\Cal L}$ is the sub-Laplacian, $T$ is the central
element of the Lie algebra on the Heisenberg group ${\Bbb H}^n,$
and $m$ satisfies the Marcinkiewicz conditions. It was proved in [MRS]
that the kernel of $m({\Cal L}, iT)$
satisfies the standard one-parameter
Calder\'on-Zygmund type estimates associated with automorphic dilations in the region where
$\vert  t \vert <\vert z\vert ^2,$
and the multi-parameter product kernel  in the region where $\vert t \vert \geq \vert z\vert ^2$
on the space ${\Bbb  C}^n\times R.$
The proof of the $L^p, 1<p<\infty,$ boundedness of $m({\Cal L}, iT)$ given in [MRS] requires lifting
 the operator to a larger
group, ${\Bbb  H}^n\times R.$ This lifts $K$, the kernel of $m({\Cal L}, iT)$ on ${\Bbb  H}^n,$
to a product kernel $\widetilde K$ on ${\Bbb H}^n\times R.$
The lifted kernel $\widetilde K$ is constructed
so that it projects to $K$ by
$$ K(z, t) =\int \limits^\infty \limits_{-\infty } {\widetilde K}(z,
t-u, u) du$$
taken in the sense of distributions.

The operator ${\widetilde T}$ corresponding to product kernel ${\widetilde K}$ can be dealt with
in terms of tensor products of operators, and one can obtain their
$L^p, 1<p<\infty,$ boundedness by the known pure product theory. Finally, the $L^p, 1<p<\infty,$
 boundedness
of operator with kernel $K$ follows from transference method of Coifman and Weiss ([CW]), using
the projection $\pi:{\Bbb H}^n\times R\rightarrow {\Bbb H}^n$ by $\pi ((z, t), u) =(z, t+u).$

Another example of implicit multi-parameter structure is the flag singular integrals on $\prods$
studied by Nagel-Ricci-Stein [NRS]. The simplest form of flag singular integral kernel $K(x,y)$ on $R^n\times R^m$
is defined through a projection of a product kernel ${\widetilde K}(x,y,z)$ defined
on $R^{n+m}\times R^m$ given by
$$ K(x, y) =\int \limits_{R^m} {\widetilde K}(x, y-z, z) dz.\tag1.2$$
A more general definition of flag singular kernel was introduced
in [NRS], see more details of definitions and applications of flag singular integrals there. We will also briefly recall them later in the introduction. Note that convolution with a flag singular kernel is a
special case of product singular kernel. As a consequence, the $L^p,
1<p<\infty ,$ boundedness of flag singular integral follows directly
from the product theory on $R^n\times R^m$. We note the regularity satisfied by flag singular
kernels is better than that of the product singular kernels. More precisely,
the singularity of the standard pure product kernel on $R^n\times R^m,$
is sets $\lbrace (x, 0)\rbrace \cup \lbrace (0, y)\rbrace$
 while the singularity of $K(x,y),$ the flag singular kernel on $R^n
\times R^m$ defined by (1.2), is a flag set given by $\lbrace (0, 0)\rbrace \subseteq \lbrace (0, y)\rbrace.$ For example, $K_1(x,y)={{1}\over{xy}}$ is a product
 kernel on $R^2$ and $K_2(x,y)={{1}\over {x(x+iy)}}$ is a flag kernel on $R^2.$

   The work of [NRS] suggests that a satisfactory Hardy space theory should be developed and boundedness of flag singular integrals on such spaces should be established.
Thus some natural questions arise. From now on, we will use the subscript $"F"$ to express function spaces or functions associated with the multi-parameter flag structure without further explanation.

{\bf Question 1:} What is the analogous estimate when $p=1$? Namely, do we have a satisfactory
flag Hardy space $H^1_F(R^n\times R^m)$ theory associated with the flag singular integral operators? More generally, can we
 develop the flag Hardy space $\phpf$ theory for all $0<p\leq 1$ such that the flag singular integral operators are
 bounded on such spaces?

{\bf Question 2:} Do we have a boundedness result on a certain
 type of $\bmof$ space for flag singular integral operators considered in [NRS]? Namely, does  an
 endpoint estimate of the result by Nagel-Ricci-Stein hold when $p=\infty$?

{\bf Question 3:} What is the duality theory of so defined flag Hardy space? More
precisely, do we have an analogue of  $BMO$ and Lipchitz type function spaces which are dual spaces
of the flag Hardy spaces.

{\bf Question 4:} Is there a Calder\'on-Zygmund decomposition in terms of functions in flag Hardy spaces $\phpf$? Furthermore, is there a satisfactory theory of interpolation on such spaces?

{\bf Question 5:} What is the difference and relationship between the  Hardy space $H^p(R^n\times R^m)$ in the pure product setting and $H^p_F(R^n\times R^m)$ in flag multiparameter setting?

  The original goal of our work is to address these questions. As in the $L^p$ theory for $p>1$
considered in [MRS], one is naturally  tempted to establish the Hardy space theory under the
 implicit multi-parameter structure  associated with the flag singular kernel by lifting method to the
 pure product setting together with the transference method in [CW]. However, this direct lifting method is not adaptable directly to the case of $p \le 1$ because the transference method is not known to be valid when $p \le 1$.
This suggests
 that a different approach in dealing with the Hardy $\phpf$ space associated with this implicit multi-parameter structure  is necessary. This motivated our
 work in this paper.   In fact, we will develop a unified approach to study multi-parameter Hardy space theory. Our approach will be carried out in the order of the following steps.

(1) We  first establish the theory of Littlewood-Paley-Stein square function $g_F$
associated with the  implicit multi-parameter structure
and the $L^p$ estimates of $g_F$ $(1<p<\infty)$. We then develop a discrete Calder\'on reproducing formula
 and a Plancherel-Polya type inequality in a test function space associated to this structure.
 As in the classical case of pure product setting, these $L^p$ estimates can be used to provide a new proof of
  Nagel-Ricci-Stein's $L^p (1<p<\infty)$ boundedenss of flag singular integral operators.

(2) We next develop the theory of Hardy spaces $H_F^p$ associated to
the multi-parameter flag structures and the boundedness of flag
singular integrals on these spaces; We then establish the
boundedness of flag singular integrals from $H^p_F$ to $L^p$.  We
refer to the reader the work of product multi-parameter Hardy space
theory by Chang-R. Fefferman [CF1-3], R. Fefferman [F1-3], Journe
[J1-2] and Pipher [P].

(3) We then establish the duality theory of the flag Hardy space $H_F^p$
and introduce the dual space $CMO^p_F$, in particular, the duality
of $H_F^1$ and the space $BMO_F.$ We then establish the boundedness of flag singular integrals on $BMO_F.$ It is worthwhile to point out that in the classical one-parameter  or pure product case, $BMO(R^n)$ or $BMO(\prods)$ is related to the Carleson measure. The space $CMO^p_F$ for all
$0<p\leq 1,$ as the dual space of $H^p_F$ introduced in this paper,  is then defined by a generalized Carleson measure.

(4) We further establish a Calder\'on-Zygmund decomposition lemma for any $\phpf$ function ($0<p<\infty$) in terms of functions in $H^{p_1}_F(\prods)$ and $H^{p_2}_F(\prods)$ with $0<p_1<p<p_2<\infty$.   Then an interpolation theorem is established between  $H^{p_1}_F(\prods)$ and $H^{p_2}_F(\prods)$ for any $0<p_2<p_1<\infty$ (it is noted that $\phpf=L^p(R^{n+m})$ for $1<p<\infty$).

In the present paper, we will use the above approach to study the Hardy space theory associated with the implicit
multi-parameter structures induced by the flag singular integrals.
We now describe our approach and results in more details.

We first introduce the continuous version of the Littlewood-Paley-Stein square function $g_F$. Inspired by the idea of lifting method of proving the $L^p(\prods)$ boundedness given in [MRS],
we will use a  lifting method  to
construct a test function defined on $\prods$, given by the non-standard convolution $*_2$ on the second variable only:
$$\psi (x,y)=\psi^{(1)}*_2\psi^{(2)}(x,y)=\int\limits_{R^m} \psi^{(1)}(x, y-z)\psi^{(2)}(z) dz,\tag1.3$$
where $\psi^{(1)}\in {\Cal S}(R^{n+m}), \psi^{(2)}\in
{\Cal S}(R^m),$ and satisfy
$$\sum \limits_j\vert  \widehat{\psi^{(1)}}(2^{-j}\xi _1,2^{-j}
\xi _2)\vert ^2=1$$
for all $(\xi _1,\xi _2)\in {R^n\times R^m}\backslash \lbrace
(0,0)\rbrace,$ and
$$ \sum \limits_k\vert {\widehat \psi }^{(2)}(2^{-k}
\eta )\vert ^2= 1$$
for all $\eta \in {R^m}\backslash \lbrace 0\rbrace,$ and the moment conditions
$$\int\limits_{R^{n+m}}x^{\alpha}y^{\beta}\psi^{(1)}(x,y)dxdy=\int\limits_{R^m}
z^{\gamma}\psi^{(2)}(z)dz=0$$ for all multi-indices $\alpha, \beta,$
and $\gamma.$ We remark here that it is this subtle convolution
$*_2$ which provides a  rich theory for the implicit multi-parameter
analysis.

For $f\in L^p, 1<p<\infty, g_F(f),$ the Littlewood-Paley-Stein square function of $f,$ is defined by
$$ g_F(f)(x,y) =\left\lbrace \sum \limits_j
\sum \limits_k\vert \psi _{j,k}*f(x, y)\vert
^2\right\rbrace ^{{{1}\over{2}}}\tag1.4$$
where functions $$\psi _{j,k}(x,y)=\psi_j ^{(1)}*_2\psi_k^{(2)}(x,y),$$
$$\psi_j^{(1)}(x,y)=2^{(n+m)j}\psi^{(1)}(2^jx,
2^jy)\,\,  \text{and}\,\, \psi_k^{(2)}(z)=2^{mk}\psi^{(2)}(2^kz).$$

We remark here that the terminology "implicit multi-parameter
structure" is  clear from the fact that the dilation $\psi
_{j,k}(x,y)$ is not induced from $\psi(x,y)$ explicitly.

By taking the Fourier transform, it is easy to see the following continuous
version of the Calder\'on reproducing formula holds
on $L^2(R^{n+m}),$
$$ f(x,y) =\sum \limits_j\sum
\limits_k \psi_{j,k}*\psi _{j,k}*f(x, y).\tag1.5$$
Note that if one considers the summation on the right hand side of (1.5) as an operator then, by the construction of
function $\psi,$ it is a flag singular integral and has the implicit multi-parameter structure as
mentioned before. Using iteration and the vector-valued Littlewood-Paley-Stein
estimate together with the Calder\'on reproducing formula on $L^2$ allows us to obtain the $L^p, 1<p<\infty,$
estimates of $g_F$.

\proclaim{ Theorem 1.1:} Let $1<p<\infty$. Then there exist constants $C_1$ and $C_2$ depending on $p$ such that for
$$ C_1\Vert f \Vert _p\leq \Vert g_F(f)
\Vert _p\leq C_2\Vert f \Vert _p.$$
\endproclaim

In order to state our results for flag singular integrals, we need to recall some definitions given in [NRS]. Following
closely from [NRS], we begin with the definitions of a class of distributions on an Euclidean space $R^N$.
A $k-normalized$ bump function on a space $R^N$ is a $C^k-$function supported on the unit ball with $C^k-$norm bounded by $1$.
As pointed out in [NRS], the definitions given below are independent of the choices of $k$, and thus we will simply refer to "normalized bump function" without specifying $k$.

For the sake of simplicity of presentations, we will restrict our considerations to the case $R^N=R^{n+m}\times R^m$. We will rephrase Definition 2.1.1 in [NRS] of product kernel in this case as follows:

\proclaim{Definition 1.2:} A product kernel on $R^{n+m}\times R^m$ is a distribution $K$ on $R^{n+m+m}$ which coincides with a $C^\infty$ function away from the coordinate subspaces $(0,0, z)$ and $(x,y, 0)$, where $(0, 0)\in R^{n+m}$ and $(x,y)\in R^{n+m}$, and satisfies

(1) (Differential Inequalities) For any multi-indices $\alpha=(\alpha_1, \cdots, \alpha_n)$, $\beta=(\beta_1,\cdots, \beta_m)$ and $\gamma_m=(\gamma_1, \cdots, \gamma_m)$
$$|\partial_{x}^{\alpha} \partial_{y}^{\beta}\partial_{z}^{\gamma}K(x, y, z)|\le C_{\alpha, \beta, \gamma} (|x|+|y|)^{-n-m-|\alpha|-|\beta|}\cdot |z|^{-m-|\gamma|}$$
for all $(x,y,z)\in R^n\times R^m\times R^m$ with $|x|+|y|\not=0$ and $|z|\not=0$.

(2) (Cancellation Condition)
$$|\int_{R^m} \partial_{x}^{\alpha} \partial_{y}^{\beta}K(x, y, z)\phi_1(\delta z)dz|\le C_{\alpha, \beta}
(|x|+|y|)^{-n-m-|\alpha|-|\beta|}$$
for all multi-indices $\alpha, \beta$ and every normalized  bump function $\phi_1$ on $R^m$ and every $\delta>0$;
$$|\int_{R^m} \partial^{\gamma}_{z} K(x, y, z)\phi_2(\delta x, \delta y)dxdy|\le C_{\gamma}
|z|^{-m-|\gamma|}$$
for every multi-index $\gamma$ and every normalized bump function $\phi_2$ on $R^{n+m}$ and every $\delta>0$;
and
$$|\int_{R^{n+m+m}}K(x,y,z)\phi_3(\delta_1 x, \delta_1 y, \delta_2 z)dxdydz|\le C$$
for every normalized bump function $\phi_3$ on $R^{n+m+m}$ and every $\delta_1>0$ and $\delta_2>0$.

\endproclaim

\proclaim{Definition 1.3:} A flag kernel on $R^n\times R^m$ is a distribution on $R^{n+m}$ which coincides with a $C^\infty$ function away from the coordinate subspace $\{(0, y)\}\subset R^{n+m}$, where $0\in R^n$ and $y\in R^m$ and satisfies

(1) (Differential Inequalities) For any multi-indices $\alpha=(\alpha_1, \cdots, \alpha_n)$, $\beta=(\beta_1,\cdots, \beta_m)$
$$|\partial_{x}^{\alpha} \partial_{y}^{\beta}K(x, y)|\le C_{\alpha, \beta} |x|^{-n-|\alpha|}\cdot (|x|+|y|)^{-m-|\beta|}$$
for all $(x,y)\in R^n\times R^m$ with $|x|\not=0$.

(2) (Cancellation Condition)
$$|\int_{R^m} \partial_{x}^{\alpha} K(x, y)\phi_1(\delta y)dy|\le C_{\alpha}
|x|^{-n-|\alpha|}$$
for every multi-index $\alpha$ and every normalized  bump function $\phi_1$ on $R^m$ and every $\delta>0$;
$$|\int_{R^n} \partial_{y}^{\beta}K(x, y)\phi_2(\delta x)dx|\le C_{\gamma}
|y|^{-m-|\beta|}$$
for every multi-index $\beta$ and every normalized bump function $\phi_2$ on $R^{n}$ and every $\delta>0$;
and
$$|\int_{R^{n+m}}K(x,y)\phi_3(\delta_1 x, \delta_2 y)dxdy|\le C$$
for every normalized bump function $\phi_3$ on $R^{n+m}$ and every $\delta_1>0$ and $\delta_2>0$.
\endproclaim

By a result in [MRS], we may assume first that a flag kernel $K$ lies in $L^1(R^{n+m})$. Thus, there exists a product kernel $K^\sharp$ on $R^{n+m}\times R^m$ such that $$K(x,y)=\int_{R^m}K^\sharp(x, y-z, z)dz.$$
Conversely, if a product kernel $K^\sharp$ lies in $L^1(R^{n+m}\times R^m)$, then
$K(x,y)$ defined as above is a flag kernel on $R^n\times R^m$. As pointed out in [MRS], we may always assume that $K(x,y)$, a flag kernel, is integrable on $\prods$ by using a smooth truncation argument.

As a consequence of
Theorem 1.1, we give a new proof of the
$L^p, 1<p<\infty,$ boundedness of
flag singular integrals due to Nagel, Ricci and Stein in [NRS].
More precisely, let $T(f)(x,y)=K*f(x,y)$ be a flag singular integral on
$R^n\times R^m.$ Then $K$ is a projection of a product kernel
$K^\sharp$ on $R^{n+m}\times R^m.$

\vskip -0.3cm

\proclaim{Theorem 1.4:} Suppose that $T$ is a flag singular integral
defined on $R^n\times R^m$ with the flag
kernel $K(x,y)=\int\limits_{R^m}K^\sharp(x,y-z,z)dz,$ where the product kernel
$K^\sharp$ satisfies the conditions of Definition 1.2 above.
Then $T$ is bounded on $L^p$ for $1<p<\infty.$ Moreover,
there exists a constant $C$ depending on $p$ such that for $f\in L^p, 1<p<\infty,$
$$\Vert  T(f)\Vert _p\leq C\Vert f\Vert_p. $$
\endproclaim

\vskip -0.4cm
In order to use
the Littlewood-Paley-Stein square function $g_F$ to
define the Hardy space, one needs to extend
the Littlewood-Paley-Stein square
function to be defined on a
suitable distribution space. For this purpose, we first introduce the product test function space on
$R^{n+m}\times R^m$.

\proclaim{Definition 1.5:} A Schwartz test function $f(x,y,z)$ defined on $R^n\times R^m
\times R^m$ is said to be a product test function on $R^{n+m}\times R^m$
if
$$ \int f(x,y,z)x^{\alpha}y^{\beta} dxdy =\int f(x,y,z)z^{\gamma} dz = 0 \tag1.6$$
for all multi-indices $\alpha, \beta, \gamma$ of nonnegative integers.

If $f$ is a product test function on $R^{n+m}\times R^m$ we denote $f\in {\Cal
 S}_{\infty}(R^{n+m}\times R^m)$ and the norm of $f$ is defined by the norm of Schwartz
test function.
\endproclaim

We now
define the test function space ${\Cal S}_F$ on $R^n\times R^m$ associated with the flag structure.
\proclaim{Definition 1.6:} A function $f(x,y)$ defined on $R^n\times R^m$ is said to be a test
function in
${\Cal S}_F(R^n\times R^m)$ if there exists
a function $f^\sharp\in {\Cal S}_{\infty}
(R^{n+m}\times R^m)$ such that
$$ f(x, y) =\int \limits_{R^m} f^\sharp(x, y-z, z) dz. \tag1.7$$
If $f\in {\Cal S}_F(R^n\times R^m),$ then
the norm of $f$ is defined by
$$\Vert f \Vert _{{\Cal S}_F(R^n\times R^m)}=
 \inf\lbrace  \Vert
f^\sharp\Vert _{{\Cal S}_{\infty}(R^{n+m}\times R^m)
}:\,\, \hbox{for all representations of $f$ in (1.7)}\rbrace. $$
We denote by $({\Cal S}_F
)^\prime$ the dual space of ${\Cal S}_F$.
\endproclaim

We would like to point out that the implicit multi-parameter structure is involved in ${\Cal
S}_F.$ Since the functions $\psi _{j,k}$ constructed above belong to ${\Cal S}_F(
R^n\times R^m),$
so the Littlewood-Paley-Stein
square function $g_F$ can be defined for all distributions in $({\Cal S}_F
)^\prime.$ Formally, we can define the flag Hardy space as follows.

\proclaim{Definition 1.7:} Let $0<p\le 1$. $\phpf=\left\lbrace  f\in ({\Cal S}_F
)^\prime: g_F(f)\in L^p(\prods)\right\rbrace.$

If $f\in \phpf,$ the norm of $f$ is defined by
$$\Vert  f \Vert _{H_F^p}=\Vert
g_F(f)\Vert _p.\tag1.8$$
\endproclaim
\vskip -0.3cm
A natural question arises
whether this
definition is independent of the choice of functions $\psi_{j,k}$. Moreover, to study the
$H_F^p$-boundedness of flag singular integrals
and establish the duality result of $H_F^p$, this formal definition is not sufficiently good. We need
to discretize
the norm of $H_F^p$. In order to obtain such a discrete $H_F^p$ norm we will prove the
Plancherel-P\^olya-type inequalities. The main tool to provide such
inequalities is the Calder\'on reproducing formula  (1.5).
To be more specific,
we will prove that the formula  (1.5) still holds on test function space ${\Cal S}_F(R^n\times R^m)
$ and its dual space $({\Cal S}_F
)^\prime$ (see Theorem 3.6 below). Furthermore, using an approximation procedure and the almost
orthogonality argument, we prove the following discrete Calder\'on reproducing formula.
\proclaim {Theorem 1.8:} Suppose that $\psi_{j,k}$ are the same as in (1.4).
Then
$$ f(x,y) =\sum \limits_j\sum \limits_k\sum \limits_J\sum
\limits_I \vert I\vert \vert J\vert {\widetilde \psi }_{j,k}(x,y, x_I,y_J)
\psi _{j,k}*f(x_I,y_J) \tag1.9$$
where ${\widetilde\psi}_{j,k}(x,y,x_I,y_J)
\in {\Cal S}_F(R^n\times R^m),
I\subset R^n, J\subset R^m,$ are dyadic cubes
with side-length $\ell(I)=2^{-j-N}$ and
$\ell(J)=2^{-k-N}+2^{-j-N}$ for a fixed large integer $N, x_I, y_J$ are any fixed points in $I, J,$
respectively, and the series in (1.9) converges in the norm of ${\Cal S}_F(R^n\times R^m)
$ and in the dual space $({\Cal S}_F)^\prime.$
\endproclaim
\vskip -0.3cm
The discrete Calder\'on reproducing formula  (1.9) provides the following Plancherel-P\^olya-type
inequalities. We  use the notation $A\approx B$ to denote that two quantities $A$ and $B$ are comparable independent of other substantial quantities involved in the context.
\proclaim{Theorem 1.9:} Suppose $\psi^{(1)}, \phi^{(1)}\in {\Cal S}(R^{n+m}),\psi^{(2)},
\phi^{(2)}\in
 {\Cal S}(R^m)$ and $$\psi(x,y)=\int \limits_{R^m}\psi^{(1)}(x, y-z)\psi^{(2)}(z)dz,$$ $$\phi(x,y)=\int
 \limits_{R^m}\phi^{(1)}(x, y-z)\psi^{(2)}(z)dz,$$ and $\psi_{jk}$, $\phi_{jk}$ satisfy the conditions in (1.4).
Then for $f\in ({\Cal S}_F)^\prime$ and
$0<p<\infty,$
$$ \Vert \lbrace \sum \limits_j\sum \limits_k\sum
\limits_J\sum \limits_I \sup_{u\in I,v\in J}\vert \psi _{j,k}*f(u,v)
\vert ^2\chi _I(x)\chi _J(y)\rbrace ^{{{1}\over{2}}}\Vert _p$$
$$\approx
\Vert \lbrace \sum \limits_j\sum \limits_k\sum \limits_J
\sum \limits_I \inf_{u\in I,v\in J}\vert \phi _{j,k}*f(u,v)
\vert ^2\chi _I(x)\chi _J(y)\rbrace^{{{1}\over{2}}}\vert \vert _p \tag1.10$$
where $\psi _{j,k}(x,y)$ and $\phi _{j,k}(x,y)$ are defined as in (1.4), $I\subset R^n, J\subset R^m,$ are dyadic cubes
with side-length $\ell(I)=2^{-j-N}$ and
$\ell(J)=2^{-k-N}+2^{-j-N}$ for a fixed large integer $N,  \chi _I$ and $\chi _J$ are indicator
functions of $I$ and $J$, respectively.
\endproclaim
 The Plancherel-P\^olya-type inequalities in Theorem 1.9 give the discrete Littlewood-Paley-Stein
square function
$$ g_F^d(f)(x,y)=\left\lbrace \sum \limits_j\sum \limits_k\sum
\limits_J\sum \limits_I \vert \psi _{j,k}*f(x_I,y_J)
\vert ^2\chi _I(x)\chi _J(y)\right\rbrace ^{{{1}\over{2}}}\tag1.11$$
where $I, J, x_I,$ and $ y_J$ are the same as in Theorem 1.8 and Theorem 1.9.

From this it is easy to see that the
Hardy space $H_F^p$ in (1.8) is well defined and the $H_F^p$ norm of $f$ is equivalent to the
$L^p$ norm of $g_F^d$.
By use of the Plancherel-P\^olya-type inequalities, we will prove the boundedness of
flag singular integrals on $H_F^p.$
\vskip -0.3cm
\proclaim{Theorem 1.10:} Suppose that $T$ is a flag singular integral with the
kernel $K(x,y)$ satisfying the same conditions as in Theorem 1.4.
Then $T$ is bounded on $H_F^p,$ for $0<p\leq
1.$ Namely, for all $0<p\le 1$ there exists a constant $C_p$ such that
$$ \Vert  T(f)\Vert  _{H_F^p}\leq C_p\Vert
 f\Vert _{H_F^p}.$$
\endproclaim
\vskip -0.5cm
To obtain the $H_F^p\to L^p$ boundedness of flag singular integrals, we prove the following general result:
\proclaim{Theorem 1.11}Let $0<p\le 1$.
If $T$ is a linear operator which is bounded on $L^2(R^{n+m})$ and $H_F^p(R^n\times R^m),$ then $T$ can be extended to
a bounded operator from $H_F^p(R^n\times R^m)$ to $L^p(R^{n+m})$.
\endproclaim
\vskip -0.3cm From the proof, we can see that this general result
holds in a very broad setting, which includes the classical
one-parameter and product Hardy spaces and the Hardy spaces on
spaces of homogeneous type.

In particular, for flag singular integral we can deduce from this general result the following
\vskip -0.5cm
\proclaim{Corollary 1.12:} Let $T$ be a flag singular integral as in Theorem 1.4. Then $T$ is
bounded from $H_F^p(R^n\times R^m)$ to $L^p(R^{n+m})$ for $0<p\le 1.$
 \endproclaim
\vskip -0.3cm
To study the duality of $H_F^p,$ we introduce the space $CMO^p_F.$
\vskip -0.7cm
\proclaim{Definition 1.13:} Let $\psi_{j,k}$ be the same as in (1.4). We say that $f\in
CMO^p_F$ if $f\in ({\Cal S}_F
)^\prime$ and it has the finite norm $\Vert f\Vert_{CMO^p_F}$ defined by
\vskip -0.5cm
$$\sup_{\Omega}\left\lbrace {{1}\over {{\vert \Omega\vert}}^{{{2}\over
{p}}-1}}\sum_{j, k}\int
\limits_{\Omega} \sum
\limits_{I, J: I\times J\subseteq \Omega}
\vert \psi _{j,k}*f(x,y)
\vert ^2 \chi_I(x)\chi_J(y)dx dy \right\rbrace^{{1}\over {2}} \tag1.12$$
\vskip -0.3cm
for all open sets $\Omega$ in $R^n\times R^m$ with finite measures, and
$I\subset R^n, J\subset R^m,$ are dyadic cubes
with side-length $\ell(I)=2^{-j}$ and
$\ell(J)=2^{-k}+2^{-j}$ respectively.
\endproclaim

Note that the Carleson measure condition is used and the implicit multi-parameter structure is involved in $CMO^p_F$ space. When $p=1,$ as
usual, we denote by $BMO_F$ the space $CMO^1_F.$
To see the space $CMO^p_F$ is well defined, one needs to show the definition of $CMO^p_F$
is independent of the choice of the functions $\psi _{j,k}.$ This  can be proved, again as in the
Hardy space
$H_F^p,$ by the following Plancherel-P\^olya-type inequality.
\proclaim{Theorem 1.14:} Suppose $\psi, \phi$ satisfy the same conditions as in
Theorem 1.9.
Then for $f\in ({\Cal S}_F)^\prime,$
$$ \sup_{\Omega}\left\lbrace{{1}\over{\vert \Omega \vert^{{{2}\over
{p}}-1}}}\sum \limits_j\sum \limits_k\sum
\limits_{I\times J\subseteq \Omega} \sup_{u\in I,v\in J}\vert \psi _{j,k}*f(u,v)
\vert ^2\vert I\vert \vert J\vert \right\rbrace^{{1}\over {2}}
\approx $$
$$\sup_{\Omega}\left\lbrace{{1}\over{\vert \Omega \vert^{{{2}\over
{p}}-1}}}\sum \limits_j\sum \limits_k\sum
\limits_{I\times J\subseteq \Omega}
\inf_{u\in I,v\in J}\vert \phi _{j,k}*f(u,v)
\vert ^2\vert I\vert \vert J\vert \right\rbrace^{{1}\over {2}} \tag1.13$$
where $I\subset R^n, J\subset R^m,$ are dyadic cubes
with side-length $\ell(I)=2^{-j-N}$ and
$\ell(J)=2^{-k-N}+2^{-j-N}$ for a fixed large integer $N$ respectively, and $\Omega$ are all open
sets in $R^n\times R^m$ with finite measures.
\endproclaim

To show that space $CMO^p_F$ is the dual space of $H_F^p$, we also need to introduce the sequence
spaces.
\proclaim{Definition 1.15:} Let $s^p$ be the collection of all sequences $s =\lbrace
s_{I\times J}\rbrace$ such that
$$\Vert s\Vert_{ s^p}=\left\Vert \left\lbrace \sum\limits_{j,k}\sum \limits_
{I, J}\vert  s_{I\times J}\vert^2\vert I\vert^{-1}\vert J\vert^{-1}
\chi _I(x)\chi _J(y)\right\rbrace^{{1}\over{2}}\right\Vert_{L^p}<\infty,$$
where the sum runs over all dyadic cubes $I\subset R^n, J\subset R^m$
with side-length $\ell(I)=2^{-j-N}$ and
$\ell(J)=2^{-k-N}+2^{-j-N}$ for a fixed large integer $N$, and $\chi _I,$ and $\chi _J$
are indicator functions of $I$ and $J$ respectively.

Let $c^p$ be the collection of all sequences $s =\lbrace s_{I\times J}\rbrace$ such that
$$\Vert s\Vert_{ c^p}=\sup_{\Omega}\left\lbrace {{1}\over {\vert \Omega\vert^{{{2}\over
{p}}-1}}}\sum\limits_{j,k}
\sum_{I,J: I\times J\subseteq \Omega}
\vert s_{I\times J}\vert ^2\right\rbrace^{{1}\over {2}}< \infty,\tag1.14$$
where $\Omega$ are all open sets in $R^n\times R^m$
with finite measures and the sum runs over all dyadic cubes
$I\subset R^n, J\subset R^m,$
with side-length $l(I)=2^{-j-N}$ and
$l(J)=2^{-k-N}+2^{-j-N}$ for a fixed large integer $N$.
\endproclaim

We would like to point out again that certain dyadic rectangles used in $s^p$ and $c^p$ reflect the
implicit multi-parameter structure. Moreover, the Carleson measure condition is used in the definition of $c^p$. Next, we obtain the following duality theorem.
\proclaim{Theorem 1.16:} Let $0<p\le 1$. Then we have $(s^p)^*=c^p$. More precisely, the map which maps $s=
\lbrace s_{I\times J}\rbrace$ to $<s, t>\equiv\sum\limits_{I\times J}s_{I\times J}\overline{t}_{I
\times J}$
defines a continuous linear functional on $s^p$ with operator norm $\Vert t \Vert_{(s^p)^*}
\approx \Vert t \Vert_{c^p}$, and moreover, every $\ell \in (s^p)^*$ is of this form
for some $t\in c^p.$
\endproclaim

When $p=1,$ this theorem in the one-parameter setting on  $R^n$ was proved in [FJ]. The proof given in [FJ] depends on
estimates of certain distribution functions,
which seems to be difficult to apply to the multi-parameter case. For all $0<p\leq 1$ we give a
simple and more constructive proof
of Theorem 1.16, which uses the stopping time argument  for sequence
 spaces. Theorem 1.16 together with the discrete Calder\'on reproducing formula and the
Plancherel-P\^olya-type
inequalities yields the duality of $H_F^p$.
\vskip -0.7cm
\proclaim{Theorem 1.17:} Let $0<p\le 1$. Then $(H_F^p)^*=CMO^p_F$. More precisely, if $g\in CMO^p_F$, the
map $\ell_{g}$ given by $\ell_g(f)=<f,g>,$ defined initially for $f\in {\Cal S}_F,$ extends to a
continuous linear functional on $H_F^p$ with $\Vert \ell_g\Vert \approx \Vert g
\Vert_{CMO^p_F}.$ Conversely, for every $\ell\in (H_F^p)^*$ there exists some $g\in CMO^p_F$
so that $\ell =\ell_g.$ In particular, $(H_F^1)^*=BMO_F.$
\endproclaim

\vskip -0.3cm
As a consequence of the duality of $H_F^1$ and the $H_F^1$-boundedness of flag singular integrals,
we obtain the $BMO_F$-boundedness of flag singular integrals. Furthermore, we will see that
$L^{\infty}\subseteq BMO_F$ and, hence, the $L^{\infty}\to BMO_F$ boundedness of flag singular
integrals is also obtained. These provide the endpoint results of those in [MRS] and [NRS]. These can be summarized as follows:

\vskip -0.4cm

\proclaim{ Theorem 1.18:} Suppose that $T$ is a flag singular integral as in Theorem 1.4. Then
$T$ is bounded on $BMO_F.$ Moreover, there exists a constant $C$ such that
$$\Vert T(f)\Vert_{BMO_F}\leq C\Vert f\Vert_{BMO_F}.$$
\endproclaim

\vskip -0.5cm
Next we prove the Calder\'on-Zygmund decomposition and interpolation theorems on the flag Hardy spaces.
 We note that $\phpf=L^p(R^{n+m})$ for  $1<p<\infty$.

\proclaim{Theorem 1.19}(Calder\'on-Zygmund  decomposition for flag Hardy spaces)
Let $0<p_2\le 1, p_2<p<p_1<\infty$ and let $\alpha>0$ be given and $f\in \phpf$. Then we may write $f=g+b$ where $g\in H^{p_1}_F(\prods)$ with $p<p_1<\infty$ and $b\in H^{p_2}_F(\prods)$ with $0<p_2<p$ such
 that $||g||^{p_1}_{H^{p_1}_F}\le C\alpha^{p_1-p}||f||^p_{H^p_F}$ and
 $||b||^{p_2}_{H^{p_2}_F}\le C\alpha^{p_2-p}||f||^p_{H^p_F}$, where $C$ is an absolute constant.
\endproclaim
\vskip -1cm
\proclaim{Theorem 1.20} (Interpolation theorem on flag Hardy spaces) Let $0<p_2<p_1<\infty$ and $T$ be a linear operator which is bounded from $H^{p_2}_F$ to $L^{p_2}$ and bounded from $H^{p_1}_F$ to  $L^{p_1}$, then $T$ is bounded from $H^p_F$ to $L^p$ for all $p_2<p<p_1$. Similarly, if $T$ is bounded on $H^{p_2}_F$ and $H^{p_1}_F$, then $T$ is bounded on $H^p_F$ for all $p_2<p< p_1$.
\endproclaim

We point out that the Calder\'on-Zygmund decomposition in pure product domains was established for
 all $L^p$ functions ($1<p<2$) into $H^1$ and $L^2$ functions by using atomic decomposition on $H^1$ space
 (see for more precise statement in Section 6).

We end the introduction of this paper with the following remarks.  First of all,
 our approach in this paper    will enable us to revisit the pure product  multi-parameter theory using
 the corresponding discrete Littlewood-Paley-Stein theory. This will provide alternative proofs of some of the known
 results of Chang-R. Fefferman, R. Fefferman, Journe, Pipher and establish some new results in the pure product setting. We will clarify all these in the future.
  Second, as we can see from the definition of flag kernels,
the regularity satisfied by flag singular
kernels is better than that of the product singular kernels.  It is thus
natural to conjecture  that the Hardy space associated
with flag singular integrals should be larger than the classical pure
product Hardy space. This is indeed the case. In fact, if we define the flag kernel on $R^n\times
 R^m$ by
 $$K(x,y)=\int\limits_{R^n}{\widetilde{{\widetilde K}}}(x-z, z, y)dz,$$
 where ${\widetilde{\widetilde K}(x, z, y)}$
is a pure product kernel on $R^n\times R^{n+m},$ and let $\widetilde{H^p_F}$ be the flag Hardy space
associated with this structure, thus we have shown in a forthcoming paper  that $H^p(R^n\times R^m)=
H^p_F(R^n\times R^m)\cap {\widetilde{H^p_F}}(R^n\times R^m).$
Results in [MRS] and [NRS] together with those in this paper  demonstrate  that the
implicit multi-parameter structure, the geometric property of sets of singularities
and regularities of singular kernels and multipliers are closely related.
 Third,  the authors have  carried out in [HL] the discrete Littlewood-Paley-Stein analysis and Hardy space
 theory in the multi-parameter structure  induced by the Zygmund dilation and proved the endpoint estimates such as
  boundedness of singular integral operators considered by Ricci-Stein [RS] on $H^p_Z$ ($0<p\le 1$) and $BMO_Z$, the Hardy and BMO spaces associated with the Zygmund dilation, e.g.,
  on $R^3,$ given by
$\delta \dot (x,y,z)=(\delta_1 x, \delta_2 y, \delta_1\delta_2 z),
\delta_1, \delta_2>0$, where  the $L^p$ ($1<p<\infty$) boundedness
has been established  (see [RS] and [FP]).

This paper is organized as follows. In Section 2, we establish the $L^p$ estimates for the multi-parameter Littlewood-Paley-Stein $g-$ function for $1<p<\infty$ and prove Theorems 1.1 and 1.4. In Section 3 we first introduce the test function spaces associated with the multi-parameter flag structure and
show the Calder\'on reproducing formula in (1.5) still holds on test function space ${\Cal S}_F(R^n
\times R^m)$ and its dual space $({\Cal S}_F)^{\prime},$ and then prove the almost orthogonality estimates and  establish the discrete Calder\'on reproducing formula on the test function spaces, i.e., Theorem 1.8.  Some crucial strong maximal function estimates are given (e.g. Lemma 3.7) and together with the discrete Calder\'on reproducing formula we derive the Plancherel-P\^olya-type inequalities,i.e.,  Theorem 1.9.
Section 4 deals with numerous  properties of Hardy space $H^p_F$ and a  general result of bounding the $L^p$ norm of the function by its $H^p_F$ norm (see Theorem 4.3), and then prove the $H^p_F$ boundedness of flag singular integrals for all $0<p<1$, i.e., Theorem 1.10. The boundedness from $H^p_F$ to $L^p$ for all $0<p\le 1$ for the flag singular integral operators,  i.e., Theorem 1.11, is thus a consequence of Theorem 1.10 and Theorem 4.3. The duality of the Hardy space $H^p_F$ is then established in Section 5. The boundedness of flag singular integral operators on $BMO_F$ space is also proved in Section 5. Thus the proofs of Theorem 1.14, 1.16, 1.17 and 1.18 will be all given in
the Section 5. In Section 6, we prove a Calder\'on-Zygmund decomposition in flag multi-parameter setting and then derive an interpolation theorem.

{\bf Acknowledgement.} The authors wish to express their sincere thanks to Professor E. M. Stein for
 his encouragement over the past ten years to carry out the  program of developing the Hardy space
 theory in the implicit multi-parameter structure and his suggestions
during the course of this work. We also like to thank Professor J. Pipher for her interest in this work and her encouragement to us.

\vskip 0.1cm
\noindent {\bf 2. $L^p$ estimates for Littlewood-Paley-Stein square function: Proofs of Theorems 1.1 and 1.4 }
\vskip 0.1cm

The main purpose of this section is to show that the $L^p$ ($p>1$) norm of $f$ is equivalent to the $L^p$ norm of $g_F(f)$, and thus use this to provide a new proof  of the $L^p$ boundedness of flag singular integral operators given in [MRS]. Our proof here is quite different from those in [MRS] in the sense that we do not need to apply the lifting procedure used in [MRS] directly.
We first prove the $L^p$ estimate of the Littlewood-Paley-Stein
square function $g_F$.

\noindent {\bf Proof of Theorem 1.1:} The proof is similar to that in the pure product case given in [FS] and
follows from iteration and standard vector-valued Littlewood-Paley-Stein inequalities.  To see this, define
$F: R^{n+m}\rightarrow H =\ell^2$ by $F(x, y) =\lbrace\psi_j^{(1)}*f(x,y)\rbrace$ with the norm $$\Vert F\Vert_H =
\lbrace \sum\limits_j \vert \psi _j^{(1)}*f(x,y)\vert^2\rbrace^{{{1}\over{2}}}.$$
When $x$ is fixed, set
$$ {\widetilde g}(F)(x,y) = \lbrace \sum \limits_k\Vert
 \psi^{(2)}_k*_2F(x,\cdot)(y)\Vert _H^2\rbrace
^{{{1}\over{2}}}.$$
It is then easy to see that ${\widetilde g}(F)(x,y) =g_F(f)(x,y).$ If $x$
is fixed, by the vector-valued Littlewood-Paley-Stein inequality,
$$\int\limits_{R^m} {\widetilde g}(F)^p(x,y) dy\leq C\int\limits_{R^m} \Vert  F \Vert
_H^p dy.$$

However, $\Vert F \Vert _H^p= \lbrace \sum
\limits_j\vert \psi _j^{(1)}*f(x, y)\vert ^2
\rbrace^{{{p}\over{2}}},$
so integrating with respect to $x$ together with
the standard Littlewood-Paley-Stein inequality yields
$$\int \limits_{R^n}\int\limits_{R^m} g_F(f)^p(x,y) dydx\leq C\int\limits_{R^n}
 \int\limits_{R^m} \lbrace \sum
\limits_j\vert \psi _j^{(1)}*f(x, y)\vert ^2
\rbrace^{{{p}\over{2}}} dydx\leq C\Vert f \Vert
 ^p_p,$$
which shows that $||g_F(f)||_p\le C||f||_p$.

The proof of the estimate $||f||_p\le C||g_F(f)||_p$ is routine and it follows from the Calderon reproducing formula (1.5) on $L^2(R^{n+m})$, for all $f\in L^2\cap L^p$, $g\in L^2\cap L^{p^\prime}$ and $\frac{1}{p}+\frac{1}{p^\prime}=1$, and
the inequality $||g_F(f)||_p\le C||f||_p$, which was just proved. This completes the proof of Theorem 1.1.
\hfill {\bf Q.E.D.}

\noindent {\bf Remark 2.1:} Let $\psi^{(1)}\in {\Cal S}(R^{n+m})$ be supported in the unit ball in $R^{n+m}$ and $\psi^{(2)}\in {\Cal S}(R^m)$ be  supported in the unit ball of $R^m$
and
satisfy
$$\int_0^\infty |\widehat{\psi^{(1)}}(t\xi_1, t\xi_2)|^4\frac{dt}{t}=1$$
for all $(\xi _1,\xi _2)\in {R^n\times R^m}\backslash \lbrace
(0,0)\rbrace,$ and
$$\int_0^\infty |\widehat{\psi^{(2)}}(s\eta)|^4\frac{ds}{s}=1$$
for all $\eta\in R^m\backslash \{0\}.$
We define $\psi^\sharp(x,y,z)=\psi^{(1)}(x,y-z)\psi^{(2)}(z).$ Set $\psi^{(1)}_t(x,y)=t^{-n-m}\psi^{(1)}(\frac{x}{t}, \frac{y}{t})$ and $\psi^{(2)}_s(z)=s^{-m}\psi(\frac{z}{s})$ and
$$\psi_{t,s}(x,y)=\int_{R^{m}}\psi^{(1)}_t(x,y-z)\psi^{(2)}_s(z)dz.$$
Repeating the same proof as that of Theorem 1.1, we can get for $1<p<\infty$
$$\Vert  \lbrace
\int \limits^\infty \limits_0\int \limits^\infty \limits_0
\vert \psi _{t,s}*f(x, y)\vert ^2{{dt}\over{t}}{{ds}\over{s}}
\rbrace^{{{1}\over{2}}}\Vert _p\leq C\Vert f\Vert_p, \tag 2.1$$
and
$$ \Vert  f \Vert _p\approx \Vert  \lbrace
\int \limits^\infty \limits_0\int \limits^\infty \limits_0
\vert \psi_{t,s}*\psi _{t,s}*f(x, y)\vert ^2{{dt}\over{t}}{{ds}\over{s}}
\rbrace^{{{1}\over{2}}}\Vert _p.\tag2.2$$

The $L^p$ boundedness of flag singular integrals is then a consequence of Theorem 1.1 and
Remark 2.1. We give a detailed proof of this below.

\noindent {\bf Proof of Theorem 1.4:} We may first assume that $K$ is integrable function
and shall prove the $L^p, 1<p<\infty,$ boundedness of $T$ is independent of the $L^1$ norm of $K.$
The conclusion for
general $K$ then follows by the argument used in [MRS]. For all
$f\in L^p, 1<p<\infty,$ by (2.2)
$$\Vert T(f)\Vert _p\leq C
\Vert  \lbrace
\int \limits^\infty \limits_0\int \limits^\infty \limits_0
\vert \psi_{t,s}*\psi _{t,s}*K*f(x, y)\vert ^2{{dt}\over{t}}{{ds}\over{s}}
\rbrace^{{{1}\over{2}}}\Vert _p.\tag2.3$$

Now we claim the following estimate: for $f\in L^p, 1<p<\infty,$
$$ \vert \psi _{t,s}*K*f(x,y)
\vert \leq C M_s(f)(x,y),\tag2.4$$
where $C$ is a constant which is independent of the $L^1$ norm of $K$ and $M_s(f)$ is the
maximal function of $f$ defined in the first section.

Assuming (2.4) for the moment, we obtain from (2.3)
$$\Vert Tf\Vert _p\leq C\Vert \lbrace
\int \limits^\infty \limits_0\int \limits^\infty \limits_0
(M_s(\psi _{t,s}*f)(x,
y) )^2{{dt}\over{t}}{{ds}\over{s}}
\rbrace^{{{1}\over{2}}}\vert
\vert _p
\leq C\Vert
 f \Vert _p,\tag2.5$$
where the last inequality follows from the Fefferman-Stein vector-valued
maximal function and Remark 2.1.

We now prove the claim (2.4). Note that $\psi _{t,s}*K(x,y)=
\int \psi ^\sharp_{t,s}*K^\sharp(x,y-z,z)dz,$ where $
\psi ^\sharp_{t,s}(x, y, z)$ is given in Remark 2.1 and $K(x,y)=\int K^
\sharp (x, y-z, z) dz,$
where $K^\sharp(x,y,z)$ is a product kernel satisfying the conditions of definition 2.1.1 in [NRS] (or Definition 1.2 in our paper).
The estimate in (2.4) will follow by integrating with respect to $z$ variable from the following estimate:
$$ \vert \psi ^\sharp_{t,s}*K^\sharp
(x,y,z)\vert \leq C{{t}\over{(t  + \vert x\vert  +\vert  y\vert
)^{n+m+1}}}{{s}\over{(s + \vert z\vert )^{m+1}}}, \tag2.6$$
where the constant $C$ is
independent of the $L^1$ norm of $K.$ The estimate (2.6) follows from that in the pure product setting
$R^{n+m}\times R^m$ given by R. Fefferman and Stein [FS]. \hfill {\bf Q.E.D.}

\vskip .1cm
\noindent {\bf 3. Test function spaces, almost orthogonality estimates  and discrete Calder\'on reproducing formula: Proofs of Theorems 1.8 and 1.9 }
\vskip .1cm
In this section, we develop the discrete Calder\'on reproducing formula and
the Plancherel-P\^olya-type inequalities on test function spaces. These are crucial tools in establishing the theory of Hardy spaces associated with the flag type multi-parameter dilation
structure. The key ideas to provide the discrete Calder\'on
reproducing formula and the
Plancherel-P\^olya-type
inequalities are the continuous version of the Calder\'on reproducing formula on test function spaces and the almost orthogonality estimates.

To be more precise, we say that a function $a(x,y)$
for $(x,y)\in R^n\times R^n$ belongs to the class ${\Cal S}^{\infty}(R^n\times R^n)$
if $a(x,y)$ is smooth and satisfies
the differential inequalities
$$\vert \partial_x^\alpha \partial_y^\beta a(x,y)\vert\leq A_{N,\alpha,\beta}(1+\vert x-y\vert)
^{-N}\tag3.1$$
and the cancellation conditions
$$ \int a(x,y)x^{\alpha}dx =\int a(x,y)y^{\beta} dy =0\tag3.2$$
for all positive integers $N$ and multi-indices $\alpha, \beta$ of nonnegative integers.

The following  almost orthogonality estimate is the simplest one and its proof can be adapted to the more complicated orthogonal estimates in subsequent steps.

\proclaim{Lemma 3.1}
If  $\psi$ and $\phi$ are in the class
${\Cal S}^{\infty}(R^n\times R^n)$, then for any given
positive integers $L$ and $M,$ there exists a constant $C=C(L,M)$ depending only on $L, M$ and
the constants $A_{N,\alpha,\beta}$ in (3.1) such that for all $t, s>0$
$$\vert \int\limits_{R^n}\psi_t(x,z)\phi_s(z,y)dz\vert\leq C({{t}\over{s}}\wedge {{s}\over{t}})
^{M}
{({{t}\vee{s}})^{L}\over{({t}\vee{s}+\vert
x-y\vert )^{(n+L)}}},
\tag3.3$$
where $\psi_t(x,z)=t^{-n}\psi({{x}\over{t}},{{z}\over{t}})$ and $\phi_s(z,y)=s^{-n}\phi(\frac{z}{s})$, and $t\wedge s = min(t, s), t\vee s =max(t,s).$
\endproclaim

\noindent{\bf Proof:} We only consider the case $M=L=1$ and $t\ge s$. Then
$$
\align
& \int_{R^m}\psi_t(x,z)\phi_s(z,y)dz\\ = & \int_{R^m} \left[\psi_t(x, z)-\psi_t(x,y)\right]\phi_s(z,y)dz \\ = &
\int_{|z-y|\le \frac{1}{2}(t+|x-y|)}+\int_{|z-y|\ge \frac{1}{2}(t+|x-y|)}=I+II
\endalign
$$
We use the smoothness condition for $\psi_t$ and size condition for $\phi_s$ to estimate term $I$.
To estimate term $II$, we use the size condition for both $\psi_t$ and $\phi_s$.

For the case $M>1$ and $L>1$, we only need to use the Taylor expansion of $\psi_t(x, \cdot)$ at $y$ and use the moment condition of $\psi_t$.
We shall not give the details. \hfill {\bf Q.E.D.}

Similarly, if $\psi^\sharp(x,y,z,u,v,w)$
for $(x,y,z), (u,v,w)\in R^n\times R^m\times R^m$
is a smooth function and satisfies
the differential inequalities
$$\vert \partial_x^{\alpha_1} \partial_y^{\beta_1} \partial_z^{\gamma_1} \partial_u^{\alpha_2}
\partial_v^{\beta_2}\partial_w^{\gamma_2}\psi^\sharp(x,y,z,u,v,w)
\vert$$
$$\leq A_{N,M,\alpha_1,\alpha_2,\beta_1,\beta_2,\gamma_1,\gamma_2}
(1+\vert x-u\vert+\vert y-v\vert)
^{-N} (1+\vert z-w\vert)
^{-M}\tag3.4$$
and the cancellation conditions
$$ \int \psi^\sharp(x,y,z,u,v,w)x^{\alpha_1}y^{\beta_1}dxdy
 =\int \psi^\sharp(x,y,z,u,v,w)z^{\gamma_1} dz $$
$$=\int \psi^\sharp(x,y,z,u,v,w)u^{\alpha_2}v^{\beta_2} dudv=
\int \psi^\sharp(x,y,z,u,v,w)w^{\gamma_2} dw=0,\tag3.5$$
 and for fixed $x_0\in R^n,y_0\in R^m, \phi^\sharp(x,y,z,x_0,y_0)\in{\Cal S}_{\infty}(R^{n+m}
\times R^m)$
and satisfies
$$\vert \partial_x^{\alpha_1} \partial_y^{\beta_1} \partial_z^{\gamma_1}
\phi^\sharp(x,y,z,x_0,y_0)
\vert$$
$$\leq B_{N,M,\alpha_1,\beta_1,\gamma_1,}
(1+\vert x-x_0\vert+\vert y-y_0\vert)
^{-N} (1+\vert z\vert)^{-M},\tag3.6$$
for all positive integers $N, M$ and multi-indices $\alpha_1,\alpha_2, \beta_1,\beta_2,
\gamma_1, \gamma_2$ of
nonnegative integers. Then we have the following almost orthogonality estimate:

\proclaim{Lemma 3.2}
For any given positive integers $L_1, L_2$ and $K_1, K_2,$ there exists a constant
 $C=C(L_1,L_2,K_1,K_2)$
depending only on $L_1, L_2, K_1, K_2$ and
the constants in (3.4) and (3.6) such that for all positive $t, s, t^\prime, s^\prime$ we have
$$
\align
& \vert \int\limits_{R^{n+m+m}}\psi^\sharp_{t,s}(x,y,z,u,v,w)\phi^\sharp_{t',s'}(u,v,w,x_0,y_0)dudvdw
\vert\\
\leq & C({{t}\over{t'}}\wedge {{t'}\over{t}})
^{L_1}
({{s}\over{s'}}\wedge {{s'}\over{s}})
^{L_2}
{({{t}\vee{t'}})^{K_1}\over{({t}\vee{t'}+\vert
x-x_0\vert +\vert y-y_0\vert)^{(n+m+K_1)}}} {({{s}\vee{s'}})^{K_2}\over{({s}\vee{s'}+\vert
z\vert )^{(m+K_2)}}},
\tag3.7
\endalign
$$
where $\psi^\sharp_{t,s}(x,y,z,u,v,w)=t^{-n-m}s^{-m}
\psi^\sharp({{x}\over{t}},{{y}\over{t}},{{z}\over{s}},{{u}\over{t}},{{v}\over{t}},
{{w}\over{s}})$
and
$$\phi^\sharp_{t,s}(x,y,z,x_0,y_0)=
t^{-n-m}s^{-m}
\phi^\sharp({{x}\over{t}},{{y}\over{t}},{{z}\over{s}},{{x_0}\over{t}},{{y_0}\over{t}},
).$$
\endproclaim
\vskip -0.5cm
The proofs of the almost orthogonality estimate in (3.7) is similar to that in (3.3).  We will only provide a brief proof.

\noindent{\bf Proof of Lemma 3.2:} We only consider the case $L_1=L_2=K_1=K_2=1$, $t\ge t^\prime$ and $s\le s^\prime$.   Thus
$$
\align
& \vert \int\limits_{R^n\times R^m\times
 R^m}\psi^\sharp_{t,s}(x,y,z,u,v,w)\phi^\sharp_{t',s'}(u,v,w,x_0,y_0)dudvdw
\vert\\ = & \vert \int\limits_{R^n\times R^m\times
 R^m} A\cdot B dudvdw
\vert
\endalign
$$
where $$A=\psi^\sharp_{t,s}(x,y,z,u,v,w)-\psi^\sharp_{t, s}(x,y,z, x_0, y_0, w)$$
and $$B=\phi^\sharp_{t',s'}(u,v,w,x_0,y_0)-\phi^\sharp_{t',s'}(u,v,z,x_0,y_0).$$
In the above, we have used the cancelation properties
$$\int\limits_{R^n\times R^m}\phi^\sharp_{t',s'}(u,v,w,x_0,y_0)u^\alpha v^\beta dudv=0,\, \, \int\limits_{R^m}\psi^\sharp_{t,s}(x,y,z,u,v,w)w^\gamma=0$$
for all multi-indices $\alpha$, $\beta$ and $\gamma$.

 Next,
 $$\align
 &
 \vert \int\limits_{R^n\times R^m\times
 R^m} A\cdot B \,dudvdw\\ = & \int_{|u-x_0|+|v-y_0|\le \frac{1}{2}\left(t+|x-x_0|+|y-y_0|\right), \,\,|w-z|\le\frac{1}{2}(s^\prime+|z|)}A\cdot B \,dudvdw\\ + & \int_{|u-x_0|+|v-y_0|\le \frac{1}{2}\left(t+|x-x_0|+|y-y_0|\right), \,\,|w-z|\ge \frac{1}{2}(s^\prime+|z|)}A\cdot B \,dudvdw\\ + &  \int_{|u-x_0|+|v-y_0|\ge \frac{1}{2}\left(t+|x-x_0|+|y-y_0|\right), \,\,|w-z|\le \frac{1}{2}(s^\prime+|z|)}A\cdot B \,dudvdw\\+ & \int_{|u-x_0|+|v-y_0|\ge \frac{1}{2}\left(t+|x-x_0|+|y-y_0|\right), \,\,|w-z|\ge \frac{1}{2}(s^\prime+|z|)}A\cdot B \,dudvdw\\ = & I+II+III+IV
 \endalign
 $$
 For term $I$, we use the smoothness conditions on $\psi^\sharp_{t,s}$ and $\phi^\sharp_{t^\prime, s^\prime}$; For term $II$, we use the smoothness conditions on $\psi^\sharp_{t,s}$ and size conditions on  $\phi^\sharp_{t^\prime, s^\prime}$; For term $III$, we use the size conditions on $\psi^\sharp_{t,s}$ and smoothness conditions on $\phi^\sharp_{t^\prime, s^\prime}$; For term $IV$, use the size conditions on both $\psi^\sharp_{t,s}$ and $\phi^\sharp_{t^\prime, s^\prime}$. We shall not provide the details here. \hfill {\bf Q.E.D.}

The crucial feature, however, is
that the almost orthogonality estimate still holds for functions in ${\Cal S}_F(R^n\times R^m).$
To see this, we first derive the relationship of convolutions on $R^n\times R^m$ and $R^{n+m}\times R^m$. We will use this relationship frequently in this paper.

\proclaim{Lemma 3.3} Let $\psi, \phi \in {\Cal S}_F(R^n\times R^m)$, and
   $\psi^\sharp, \phi^\sharp \in {\Cal S}_{\infty}(R^{n+m}\times R^m)$
such that
$$\psi(x,y)=\int_{R^m}\psi^\sharp(x, y-z, z)dz, \,\, \phi(x,y)=\int_{R^m}\phi^\sharp(x, y-z, z)dz.$$
Then
$$(\psi*\phi)(x,y)=\int_{R^m}\left(\psi^\sharp*\phi^\sharp\right)(x, y-z, z)dz.$$
\endproclaim

Lemma 3.3 can be proved very easily. Using this lemma and the almost orthogonality estimates on $R^{n+m}\times R^m$, we can get the following
\proclaim{Lemma 3.4} For any given
positive integers $L_1, L_2$ and $K_1, K_2,$ there exists a constant $C=C(L_1,L_2,K_1,K_2)$
depending only on $L_1, L_2, K_1, K_2$
such that
if ${t}\vee{t'}\leq {s}\vee{s'},$ then
$$
\align &
\vert \psi_{t,s}*\phi_{t',s'}(x,y)\vert\\
\leq & C({{t}\over{t'}}\wedge {{t'}\over{t}})
^{L_1} ({{s}\over{s'}}\wedge {{s'}\over{s}})
^{L_2}\cdot {({{t}\vee{t'}})^{K_1}\over{({t}\vee{t'}+\vert
x\vert )^{(n+K_1)}}} {({{s}\vee{s'}})^{K_2}\over{({s}\vee{s'}+\vert
y\vert )^{(m+K_2)}}},
\endalign
$$
and if ${t}\vee{t'}\geq {s}\vee{s'},$ then
$$
\align &
\vert \psi_{t,s}*\phi_{t',s'}(x,y)\vert\\
\leq  & C({{t}\over{t'}}\wedge {{t'}\over{t}})
^{L_1} ({{s}\over{s'}}\wedge {{s'}\over{s}})
^{L_2}\cdot {({{t}\vee{t'}})^{K_1}\over{({t}\vee{t'}+\vert
x\vert )^{(n+K_1)}}} {({{t}\vee{t'}})^{K_2}\over{({t}\vee{t'}+\vert
y\vert )^{(m+K_2)}}}.
\endalign
$$
\endproclaim
\noindent{\bf Proof of Lemma 3.4:} We first remark that we will prove this lemma with $K_1, K_2, L_1, L_2$ replaced by $K_1^\prime, K_2^\prime, L_1^\prime, L_2^\prime$. Thus, we are given any fixed $K_1^\prime, K_2^\prime, L_1^\prime, L_2^\prime$.

Note that
$$\psi_{t,s}*\phi_{t',s'}(x,y) =\int\limits
_{R^m}\psi^\sharp_{t,s}*\phi^\sharp_{t',s'}(x,y-z,z)dz,$$
 where $\psi^\sharp, \phi^\sharp \in
{\Cal S}_{\infty}(R^{n+m}\times R^m),$ and
$$\psi^\sharp_{t,s}*\phi^\sharp_{t',s'}(x,y,z)=\int\limits_{R^m\times R^m\times
 R^n}\psi^\sharp_{t,s}(x-u,y-v,z-w)\phi^\sharp_{t',s'}(u,v,w)dudvdw,$$
Then by the estimate in (3.7), for any given
positive integers $L_1, L_2$ and $K_1, K_2,$ there exists a constant $C=C(L_1,L_2,K_1,K_2)$
depending only on $L_1, L_2, K_1, K_2$
such that
$$
\align
&
\vert\psi_{t,s}*\phi_{t',s'}(x,y)\vert\\
\leq  & C({{t}\over{t'}}\wedge {{t'}\over{t}})
^{L_1}
({{s}\over{s'}}\wedge {{s'}\over{s}})
^{L_2}\int\limits_{R^m}{({{t}\vee{t'}})^{K_1}\over{({t}\vee{t'}+\vert
x\vert +\vert y-z\vert)^{(n+m+K_1)}}} {({{s}\vee{s'}})^{K_2}\over{({s}\vee{s'}+\vert
z\vert )^{(m+K_2)}}}dz.\tag3.8
\endalign
$$
\vskip -0.3cm
{\bf Case 1:} If ${t}\vee{t'}\le {s}\vee{s'}$ and $|y|\ge {s}\vee{s'}$, write
 $$\align &
 \int\limits_{R^m}{({{t}\vee{t'}})^{K_1}\over{({t}\vee{t'}+\vert
x\vert +\vert y-z\vert)^{(n+m+K_1)}}} {({{s}\vee{s'}})^{K_2}\over{({s}\vee{s'}+\vert
z\vert )^{(m+K_2)}}}dz\\ =& \int_{|z|\le \frac{1}{2}|y|,\, \text{or}\,\, |z|\ge 2|y|}+\int_{\frac{1}{2}|y|\le |z|\le 2|y|}=I+II
\endalign
$$
It is easy to see that
$$
\align
|I|\le &  C {({{t}\vee{t'}})^{K_1}\over{({t}\vee{t'}+\vert
x\vert +\vert y\vert)^{(n+m+K_1)}}}\\ \le & C {({{t}\vee{t'}})^{K_1^\prime}\over{({t}\vee{t'}+\vert
x\vert)^{(n+K_1^\prime)}}} {({{t}\vee{t'}})^{K_2\prime}\over{\vert y\vert^{K_2^\prime}}} \\ \le &
C {({{t}\vee{t'}})^{K_1^\prime}\over{({t}\vee{t'}+\vert x\vert )^{(n+K_1^\prime)}}}\cdot {({{s}\vee{s'}})^{m+K_2^\prime}\over{({s}\vee{s'}+\vert
y\vert)^{m+K^\prime_2}}}
\endalign
$$
where we have taken $K_1=K_1^\prime+K_2^\prime$ and used the fact that ${t}\vee{t'}\le {s}\vee{s'}$ and $|y|\ge {s}\vee{s'}$.

Next, we estimate
$$
\align
|II| \le &  {({{s}\vee{s'}})^{K_2}\over{({s}\vee{s'}+\vert
y\vert )^{(m+K_2)}}} \int\limits_{R^m}{({{t}\vee{t'}})^{K_1}\over{({t}\vee{t'}+\vert
x\vert +\vert y-z\vert)^{(n+m+K_1)}}}dz\\ \le  & C {({{s}\vee{s'}})^{K_2}\over{({s}\vee{s'}+\vert
y\vert )^{(m+K_2)}}}\cdot {({{t}\vee{t'}})^{K_1}\over{({t}\vee{t'}+\vert
x\vert )^{(n+K_1)}}}\\ \le & C {({{s}\vee{s'}})^{K_2^\prime}\over{({s}\vee{s'}+\vert
y\vert )^{(m+K_2^\prime)}}}\cdot {({{t}\vee{t'}})^{K_1^\prime}\over{({t}\vee{t'}+\vert
x\vert )^{(n+K_1^\prime)}}}
\endalign
$$
where we have used $K_2\ge K_2^\prime$ and $K_1=K_1^\prime+K_2^\prime>K_1^\prime$.

{\bf Case 2:} If ${t}\vee{t'}\le {s}\vee{s'}$ and $|y|\le {{s}\vee{s'}}$, then
 $$\align &
 \int\limits_{R^m}{({{t}\vee{t'}})^{K_1}\over{({t}\vee{t'}+\vert
x\vert +\vert y-z\vert)^{(n+m+K_1)}}} {({{s}\vee{s'}})^{K_2}\over{({s}\vee{s'}+\vert
z\vert )^{(m+K_2)}}}dz\\ \le & \frac{1}{({s}\vee{s^\prime})^m}
\int\limits_{R^m}{({{t}\vee{t'}})^{K_1}\over{({t}\vee{t'}+\vert
x\vert +\vert y-z\vert)^{(n+m+K_1)}}}dz\\ \le &  C
\frac{({s}\vee{s^\prime})^{K_2}}{({s}\vee{s^\prime}+|y|)^{m+K_2}}{({{t}\vee{t'}})^{K_1}\over{({t}\vee{t'}+\vert
x\vert )^{(n+K_1)}}}.
\endalign
$$

{\bf Case 3:} We now consider the case ${t}\vee{t^\prime}\ge {s}\vee{s^\prime}$ and $|y|\le {t}\vee{t^\prime}$. Then
$$\align &
 \int\limits_{R^m}{({{t}\vee{t'}})^{K_1}\over{({t}\vee{t'}+\vert
x\vert +\vert y-z\vert)^{(n+m+K_1)}}} {({{s}\vee{s'}})^{K_2}\over{({s}\vee{s'}+\vert
z\vert )^{(m+K_2)}}}dz\\ \le & C {({{t}\vee{t'}})^{K_1}\over{({t}\vee{t'}+\vert
x\vert )^{(n+m+K_1)}}} \\ \le & {({{t}\vee{t'}})^{K_1^\prime}\over{({t}\vee{t'}+\vert
x\vert )^{(n+K_1^\prime)}}} \cdot {({{t}\vee{t'}})^{K_2^\prime}\over{({t}\vee{t'}+\vert
x\vert  +|y|)^{(m+K_2^\prime)}}}.
\endalign
$$
by noticing that $K_1=K_1^\prime+K_2^\prime$.

{\bf Case 4:} If we assume ${t}\vee{t^\prime}\ge {s}\vee{s^\prime}$ and $|y|\ge {t}\vee{t^\prime}$, then we divide the integral into two parts
$I$ and $II$ as in the case ${t}\vee{t^\prime}\le {s}\vee{s^\prime}$.
Thus we have
$$
\align
|I| \le & {({{t}\vee{t'}})^{K_1}\over{({t}\vee{t'}+\vert
x\vert +\vert y\vert)^{(n+m+K_1)}}} \\ \le & C {({{t}\vee{t'}})^{K_1^\prime}\over{({t}\vee{t'}+\vert
x\vert )^{(n+m+K_1^\prime)}}}\cdot {({{t}\vee{t'}})^{K_2^\prime}\over{({t}\vee{t'} +\vert y\vert)^{(m+K_2^\prime)}}}
\endalign
$$
where we have taken $K_1=K_1^\prime+K_2^\prime$.

To estimate $II$, we have
$$
\align
|II|\le & C
{({{t}\vee{t'}})^{K_1}\over{({t}\vee{t'}+\vert
x\vert )^{(n+K_1)}}} \cdot {({{s}\vee{s'}})^{K_2}\over{({s}\vee{s'}+\vert
 y\vert )^{(m+K_2)}}}\\ \le & C
{({{t}\vee{t'}})^{K_1}\over{({t}\vee{t'}+\vert
x\vert )^{(n+K_1)}}} \cdot {({{t}\vee{t'}})^{K_2}\over{({t}\vee{t'}+\vert
 y\vert )^{(m+K_2)}}}
 \endalign
 $$
\hfill{\bf Q.E.D.}

Roughly speaking, $\psi_{t,s}*\phi_{t',s'}(x,y)$ satisfies the one-parameter almost orthogonality
when ${t}\vee{t'}\geq {s}\vee{s'}$ and
the product multi-parameter almost orthogonality when ${t}\vee{t'}\leq {s}\vee{s'}.$
More precisely, we have the following

\proclaim{Corollary 3.5}
Given any positive integers $L_1, L_2, K_1$ and $K_2$, there exists a constant $C=C(L_1, L_2, K_1, K_2)>0$ such that

(i) If $t\geq s$ we obtain the
one-parameter almost orthogonality
$$\vert \psi_{t,s}*\phi_{t',s'}(x,y)\vert
\leq C({{t}\over{t'}}\wedge {{t'}\over{t}})
^{L_1} ({{s}\over{s'}}\wedge {{s'}\over{s}})
^{L_2}
{t^{K_1}\over{(t+\vert
x\vert )^{(n+K_1)}}} {t^{K_2}\over{(t+\vert
y\vert )^{(m+K_2)}}}\tag3.8$$
and, if $t\leq s$ the product multi-parameter almost orthogonality is given by
$$\vert \psi_{t,s}*\phi_{t',s'}(x,y)\vert
\leq C({{t}\over {t'}}\wedge {{t'}\over{t}})
^{L_1} ({{s}\over{s'}}\wedge {{s'}\over{s}})
^{L_2}{t^{K_1}\over{(t+\vert
x\vert )^{(n+K_1)}}} {s^{K_2}\over{(s+\vert
y\vert )^{(m+K_2)}}}.\tag3.9$$

(ii)
Similarly, if $t'\geq s',$
$$\vert \psi_{t,s}*\phi_{t',s'}(x,y)\vert
\leq C({{t}\over{t'}}\wedge {{t'}\over{t}})
^{L_1} ({{s}\over{s'}}\wedge {{s'}\over{s}})
^{L_2}
{{t'}^{K_1}\over{(t'+\vert
x\vert )^{(n+K_1)}}} {{t'}^{K_2}\over{(t'+\vert
y\vert )^{(m+K_2)}}},\tag3.10$$
and if $t'\leq s',$
$$\vert \psi_{t,s}*\phi_{t',s'}(x,y)\vert
\leq C({{t}\over {t'}}\wedge {{t'}\over{t}})
^{L_1} ({{s}\over{s'}}\wedge {{s'}\over{s}})
^{L_2}{{t'}^{K_1}\over{(t'+\vert
x\vert )^{(n+K_1)}}} {{s'}^{K_2}\over{(s'+\vert
y\vert )^{(m+K_2)}}}.\tag3.11$$
\endproclaim

 Corollary 3.5 is actually what we will use frequently in the subsequent parts of the paper. The proof of Corollary 3.5 is a case by case study and can be checked with patience. We shall omit the details of the proof here.
All these estimates will be used to prove the following continuous version of the Calder\'on reproducing formula on test function space
${\Cal S}_F(R^n\times R^m)$ and its dual space $({\Cal S}_F)^\prime.$

\proclaim{Theorem 3.6} Suppose that $\psi_{j,k}$ are the same as in (1.4).
Then
$$ f(x,y) =\sum \limits_j\sum
\limits_k \psi_{j,k}*\psi _{j,k}*f(x, y),
\tag3.12$$
where the series converges in the norm of
${\Cal S}_F$ and in
dual space $({\Cal S}_F)^
\prime.$
\endproclaim
\noindent {\bf Proof:} Suppose $f\in {\Cal S}_F$ and $f(x, y) =\int
\limits_{R^m} f^\sharp(x, y-z, z) dz,$ where $f^\sharp\in {\Cal S}_{\infty}(R^{n+m}\times R^m).$
Then, by the classical Calder\'on reproducing formula
as mentioned in the first section, for all $f^\sharp\in {L^2},$
$$f^\sharp(x, y, z) =\sum \limits_j\sum
\limits_k\psi ^\sharp_{j,k}*\psi ^\sharp_{j,k}*f^
\sharp(x, y, z),\tag3.13$$
where $\psi ^\sharp_{j,k}(x,y,z)=\psi^{(1)}_j(x,y)\psi^{(2)}_k(z).$

We {\bf claim} that the above series in (3.13) converges in ${\Cal S}_{\infty}(R^{n+m}\times R^m).$
This claim yields
$$
\align
&
\Vert f(x,y) - \sum\limits_{-N\leq j\leq N}\sum
\limits_{-M\leq k \leq M}\psi_{j,k}*\psi _{j,k}*f(x,
y) \Vert _{{\Cal S}_F}\\
= &
\Vert \int \limits_{R^m} [f^\sharp(x, y-z, z)- \sum
\limits_{-N \leq j \leq N}\sum\limits_{-M\leq k \leq M}\psi^
\sharp_{j,k}*\psi ^\sharp_{j,k}*f^\sharp(x, y-z, z) ]
dz\Vert _{{\Cal S}_F}\\ \leq  & \Vert f^\sharp(x, y, z)- \sum\limits_{-N\leq j \leq N}
\sum\limits_{-M\leq k \leq M}\psi^\sharp_{j,k}*\psi
^\sharp_{j,k}*f^\sharp(x, y, z) \Vert
_{{\Cal S}_{\infty}}
\endalign
$$
where the last term above goes to zero as $N$
and $M$ tend to infinity by the above claim.

To show the claim, it suffices to prove that all the following three summations $$\sum\limits_{\vert j\vert > N}
\sum\limits_{\vert k \vert\leq M}\psi^\sharp_{j,k}*\psi
^\sharp_{j,k}*f^\sharp, \sum\limits_{\vert j\vert \leq N}
\sum\limits_{\vert k \vert > M}\psi^\sharp_{j,k}*\psi
^\sharp_{j,k}*f^\sharp, \sum\limits_{\vert j\vert > N}
\sum\limits_{\vert k \vert > M}\psi^\sharp_{j,k}*\psi
^\sharp_{j,k}*f^\sharp$$ tend to zero in ${\Cal S}_{\infty}(R^{n+m}\times R^m)$ as $N$ and $M$
tend to infinity. Since all proofs are similar for each of the summations,
we only prove the assertion for the first summation which we denote by $ f^\sharp_{N,M}$. Note that
$$f^\sharp_{N,M}(x,y,z)=\sum\limits_{\vert j\vert > N}
\sum\limits_{\vert k \vert\leq M}\int\limits_{R^n\times R^m\times R^m}\psi^\sharp_{j,k}*\psi
^\sharp_{j,k}(x-u,y-v,z-w)f^\sharp(u,v,w)dudvdw$$
where $\psi^\sharp*\psi
^\sharp(x-u,y-v,z-w)$ satisfies the conditions (3.4) and (3.5), and $f^\sharp(u,v,w)$ satisfies
the conditions (3.6) with $x_0=y_0=0.$ The almost orthogonality estimate in (3.7) with $t=2^{-j},
s=2^{-k}, t'=s'=1$ and $x_0=y_0=0,$ implies
$$
\align &
\vert\int\limits_{R^n\times R^m\times R^m}\psi^\sharp_{j,k}*\psi
^\sharp_{j,k}(x-u,y-v,z-w)f^\sharp(u,v,w)dudvdw\vert\\ \leq  & C 2
^{-\vert j\vert L_1}
2^{-\vert k\vert L_2}
{({{2^{-j}}\vee{1}})^{K_1}\over{({2^{-j}}\vee{1}+\vert
x\vert +\vert y\vert)^{(n+m+K_1)}}} {({{2^{-k}}\vee{1}})^{K_2}\over{({2^{-k}}\vee{1}
+\vert z\vert )^{(m+K_2)}}}.
\endalign
$$

This, by taking $L_1>K_1$ and $L_2>K_2,$ gives us
$$\lim_{N,M\to\infty}\sup_{x\in R^n, y\in R^m, z\in R^m}(1+\vert
x\vert +\vert y\vert)^{n+m+K_1}(1+\vert z\vert )^{m+K_2}
\vert f^\sharp_{N,M}(x,y,z)\vert=0.$$

Since $\partial_x^\alpha\partial_y^\beta\partial_z^\gamma (f^\sharp_{N,M})(x,y,z)=
(\partial_x^\alpha\partial_y^\beta\partial_z^\gamma f^\sharp)_{N,M}(x,y,z)$ and applying the
above estimate to $\partial_x^\alpha\partial_y^\beta\partial_z^\gamma f^\sharp$ which
also satisfies the conditions in (3.6) with $x_0=y_0=0,$ we obtain
$$\lim_{N,M\to\infty}\sup_{x\in R^n, y\in R^m, z\in R^m}(1+\vert
x\vert +\vert y\vert)^{n+m+K_1}(1+\vert z\vert )^{m+K_2}
\vert(\partial_x^\alpha\partial_y^\beta\partial_z^\gamma f^\sharp)_{N,M}(x,y,z)\vert=0,$$
which shows the claim.

The convergence in dual space follows from the duality argument. The proof of Theorem 3.6 is
complete. \hfill {\bf Q.E.D.}

Using Theorem 3.6, we prove the discrete Calder\'on reproducing formula.

\noindent {\bf Proof of Theorem 1.8:}
We first discretize (3.12) as follows.
For $f\in {\Cal S}_F,$ by (3.12) and using an idea similar to that of decomposition of the identity operator due to Coifman, we can rewrite
$$ f(x,y) =\sum \limits_{j, k}\sum \limits_{I, J}
\int \limits_J\int \limits_I
\psi_{j,k}(x-u,y-w)
\left(\psi _{j,k}*f\right)(u, w)dudw$$
$$=\sum \limits_{j,k}\sum \limits_{I, J}
\left[\int \limits_J\int \limits_I
\psi_{j,k}(x-u,y-w)
dudw\right]\left(\psi _{j,k}*f\right)(x_I, y_J) + {\Cal R}(f)(x, y).\tag3.14$$

We shall show that ${\Cal R}$ is bounded on ${\Cal S}_F$
with the small norm when $I$ and $J$ are dyadic cubes in $R^n$ and $R^m$ with
side length $2^{-j-N}$ and
$2^{-k-N}+2^{-j-N}$ for a large given integer $N$, and $x_I, y_J$ are any fixed points in $I, J$,
respectively.

To do this, assuming $f(x, y) =\int
\limits_{R^m} f^\sharp(x, y-z, z) dz,$ where $f^\sharp\in {\Cal S}_{\infty}(R^{n+m}\times
 R^m)$. When $k\le j$, we write
 $$
 \align
& \left(\psi _{j,k}*f\right)(u, w)-\left(\psi _{j,k}*f\right)(x_I, y_J)\\= &
\int_{R^m}\left(\psi^\sharp_{jk}*f^\sharp\right)(u, v, w-v)dv-\int_{R^m}\left(\psi^\sharp_{jk}*f^\sharp\right)(x_I, v, y_J-v)dv\\
= & \int \left[\psi^\sharp_{jk}(u-u^\prime, v-v^\prime, w-v-w^\prime)-\psi^\sharp_{jk}(x_I-u^\prime, v-v^\prime, y_J-v-w^\prime)\right]f^\sharp(u^\prime, v^\prime, w^\prime)du^\prime dv^\prime dw^\prime dv
\endalign
$$
 where the last integral above is over $R^n\times R^m\times R^m\times R^m$.

 When $k>j$, we write
  $$
 \align
& \left(\psi _{j,k}*f\right)(u, w)-\left(\psi _{j,k}*f\right)(x_I, y_J)\\= &
\int_{R^m}\left(\psi^\sharp_{jk}*f^\sharp\right)(u, w-v, v)dv-\int_{R^m}\left(\psi^\sharp_{jk}*f^\sharp\right)(x_I, y_J-v, v)dv\\
= & \int \left[\psi^\sharp_{jk}(u-u^\prime, w-v-v^\prime, v-w^\prime)-\psi^\sharp_{jk}(x_I-u^\prime, y_J-v-v^\prime, v-w^\prime)\right]f^\sharp(u^\prime, v^\prime, w^\prime)du^\prime dv^\prime dw^\prime dv
\endalign
$$
We now set
$$
\align & {\Cal R}(f)(x,y)\\
 = & \sum \limits_{j,k}\sum \limits_{I, J}
\int \limits_J\int \limits_I
\psi_{j,k}(x-u,y-w)
\left[\left(\psi _{j,k}*f\right)(u, w)-\left(\psi _{j,k}*f\right)(x_I, y_J)\right]dudw\\
= & \int \int \int \int {\Cal R}^\sharp(x,y-z,z,u^\prime,v^\prime,w^\prime) f^\sharp(u^\prime,v^\prime,w^\prime)
du^\prime dv^\prime dw^\prime dz\\
= & \int\limits_{R^m}{\Cal R}^\sharp(f^\sharp)(x, y-z,z)dz,
\endalign
$$
where ${\Cal R}^\sharp(x,y,z, u^\prime, v^\prime, w^\prime)$ is the kernel of
${\Cal R}^\sharp$ and
$$
\align &{\Cal  R}^\sharp(x,y-z,z,u^\prime,v^\prime,w^\prime)\\ = & \sum \limits_{k\leq j}\sum \limits_{I,J}\int \int \limits_J\int\limits_I {\psi^{(1)}_j}(x-u, y-z-w){\psi^{(2)}_k}(z)\\ \times &
\left[\psi^{(1)}_j (u-u^\prime, v-v^\prime )
\psi^{(2)}_k(w-v-w^\prime)-
\psi^{(1)}_j (x_I-u^\prime ,v-v^\prime )\psi^{(2)}_k(y_J-w^\prime-v)\right]dudwdv\\
+ & \sum \limits_j\sum \limits_{k>j}\sum \limits_J\sum
\limits_I\int \int \limits_J\int\limits_I {\psi^{(1)}_j}(x-u, y-z-w){
\psi^{(2)}_k}(z)\\ \times &
\left[\psi^{(1)}_j (u-u^\prime ,w-v-v^\prime )
-\psi^{(1)}_j (x_I-u^\prime ,y_J-v-v^\prime )\right]\psi^{(2)}_k(v-w^\prime)dudwdv.
\endalign
$$
Using the change of variables from $z$ to $z+v-w$ in the term of $k\le j$, and $z$ to $z-v$ in the term of $k>j$, we can rewrite
$$
\align
&{\Cal R}^\sharp(x,y,z,u^\prime,v^\prime,w^\prime)\\ = & \sum \limits_{k\leq j}\sum \limits_{I, J}\int \int \limits_J\int\limits_I {\psi^{(1)}_j}(x-u, y-v){
\psi^{(2)}_k}(z+v-w)\\ \times &
\left[\psi^{(1)}_j (u-u^\prime ,v-v^\prime )
\psi^{(2)}_k(w-w^\prime-v)-
\psi^{(1)}_j (x_I-u^\prime ,v-v^\prime )\psi^{(2)}_k(y_J-w^\prime-v)\right]dudwdv\\
+ & \sum \limits_j\sum \limits_{k>j}\sum \limits_J\sum
\limits_I\int \int \limits_J\int\limits_I {\psi^{(1)}_j}(x-u, y+v-w){
\psi^{(2)}_k}(z-v)\\ \times &
\left[\psi^{(1)}_j (u-u^\prime ,w-v-v^\prime )
-\psi^{(1)}_j (x_I-u^\prime ,y_J-v-v^\prime )\right]\psi^{(2)}_k(v-w^\prime)dudwdv\\ = & \sum\limits_{k\le j}A_{jk}+\sum\limits_{k>j}B_{jk}.
\endalign
$$
\vskip -0.5cm
We {\bf claim} that ${\Cal R}^\sharp$ is bounded in ${\Cal S}_{\infty}(R^{n+m}\times R^m).$

To see this,
write
$$
\align
A_{jk} = &
\sum \limits_J\sum
\limits_I\int \int \limits_J\int\limits_I {\psi^{(1)}_j}(x-u, y-v){
\psi^{(2)}_k}(z+v-w)\\ \times &
[\psi^{(1)}_j (u-u^\prime ,v-v^\prime )
-\psi^{(1)}_j (x_I-u^\prime ,v-v^\prime )]\psi^{(2)}_k(w-w^\prime-v)dudwdv\\ + &
\sum \limits_J\sum
\limits_I\int \int \limits_J\int\limits_I {\psi^{(1)}_j}(x-u, y-v){
\psi^{(2)}_k}(z+v-w)\\ \times &
\psi^{(1)}_j (x_I-u^\prime ,v-v^\prime )[\psi^{(2)}_k(w-w^\prime-v)-
\psi^{(2)}_k(y_J-w^\prime-v)]dudwdv\\ = & A_{jk}^{(1)}(x,y,z, u^\prime, v^\prime, w^\prime)+A_{jk}^{(2)}(x,y,z, u^\prime, v^\prime, w^\prime)
\endalign
$$
It is not difficult to check that $\int_{R^m}\psi_k^{(2)}(z+v-w)\psi^{(2)}_k(w-w^\prime-v)dw$
 satisfies all the conditions as $\psi_k^{(2)}(z-w^\prime)$ does  with the comparable constants  of ${\Cal S}_F(\prods)$ norm  and that
 $$
 \sum \limits_J\sum
\limits_I\int \int \limits_J\int\limits_I {\psi^{(1)}_j}(x-u, y-v)
[\psi^{(1)}_j (u-u^\prime ,v-v^\prime )
-\psi^{(1)}_j (x_I-u^\prime ,v-v^\prime )]dudv
 $$
 satisfies all conditions as $ {\psi^{(1)}_j}(x-u^\prime, y-v^\prime)$ but with the constants of $\sf$ norm replaced by
 $C2^{-N}$. This follows from the smoothness condition on $\psi_j^{(1)}$ (say the mean-value theorem) and the fact that $u, x_I\in I$ and $l(I)=2^{-N-j}$.

 If we write
 $$
 \align
 &
A_{jk}^{(1)}(x,y,z, u^\prime, v^\prime, w^\prime)\\= & \sum \limits_J\sum
\limits_I\int \int \limits_J\int\limits_I {\psi^{(1)}_j}(x-u, y-v)
[\psi^{(1)}_j (u-u^\prime ,v-v^\prime )
-\psi^{(1)}_j (x_I-u^\prime ,v-v^\prime )]\\ \times & \left[\int_{R^m}{
\psi^{(2)}_k}(z+v-w)\psi^{(2)}_k(w-w^\prime-v)dw\right]dudv
 \endalign
 $$
 then the function $A_{jk}^{(1)}(x,y,z, u^\prime, v^\prime, w^\prime)$ satisfies all conditions as $$\psi_j^{(1)}(x-u^\prime, y-v^\prime)\cdot \psi^{(2)}_k(z-w^\prime)$$ does
  but with the ${\Cal S}_F(\prods)$ norm constant replaced by $C2^{-N}$.

   By the proof of Theorem 3.6, we conclude that
 $$\int A_{jk}^{(1)}(x,y,z, u^\prime, v^\prime, w^\prime)f^\sharp(u^\prime, v^\prime, w^\prime)du^\prime dv^\prime dw^\prime$$
 is a test function in ${\Cal S}_{\infty}(R^{n+m}\times R^m)$ and its test function norm is bounded by $C2^{-N}$.

 Similarly,
 $$\sum \limits_J  \int \limits_J {
\psi^{(2)}_k}(z+v-w)[\psi^{(2)}_k(w-w^\prime-v)-
\psi^{(2)}_k(y_J-w^\prime-v)]dw
 $$
 satisfies all conditions as $\psi_k^{(2)}(z-w^\prime)$ but with the constant replaced by $C2^{-N}$, which follows from the smoothness condition on $\psi_k^{(2)}$ and the fact $w, y_J\in J$ and $l(J)=2^{-N-j}+2^{-N-k}$ and $k\le j$.

 We also note that
 $$
 \sum \limits_J\sum
\limits_I\int \int \limits_J\int\limits_I {\psi^{(1)}_j}(x-u, y-v)
\psi^{(1)}_j (x_I-u^\prime ,v-v^\prime )dudv
$$
satisfies  the same conditions as $\psi_j^{(1)}(x-u^\prime, y-v^\prime)$ does with comparable ${\Cal S}_F(\prods)$ norm constant.
Thus,
we conclude that
 $$\int A_{jk}^{(2)}(x,y,z, u^\prime, v^\prime, w^\prime)f^\sharp(u^\prime, v^\prime, w^\prime)du^\prime dv^\prime dw^\prime$$
 is a test function in ${\Cal S}_{\infty}(R^{n+m}\times R^m)$ and its test function norm is bounded by $C2^{-N}$.
Therefore,
$$\int A_{jk}(x,y,z, u^\prime, v^\prime, w^\prime)f^\sharp(u^\prime, v^\prime, w^\prime)du^\prime dv^\prime dw^\prime$$
 is a test function in ${\Cal S}_{\infty}(R^{n+m}\times R^m)$ and its test function norm is bounded by $C2^{-N}$.

Similarly, we can conclude that
$$\int B_{jk}(x,y,z, u^\prime, v^\prime, w^\prime)f^\sharp(u^\prime, v^\prime, w^\prime)du^\prime dv^\prime dw^\prime$$
 is also a test function in ${\Cal S}_{\infty}(R^{n+m}\times R^m)$ and its test function norm is bounded by $C2^{-N}$.

This shows that
${\Cal R}^\sharp(f^\sharp)(x,y,z)\in {\Cal S}_{\infty}(R^{n+m}\times R^m)$ and
$$||{\Cal R}^\sharp(f^\sharp)||_{{\Cal S}_{\infty}(R^{n+m}\times R^m)}\le C2^{-N}||f^\sharp||_{{\Cal S}_{\infty}(R^{n+m}\times R^m)},$$
which implies that
${\Cal R}(f)\in \sf$ and
$$||{\Cal R}(f)||_{\sf}\le C2^{-N}||f||_{\sf}.\tag3.15$$

By (3.14) together with the boundedness of ${\Cal R}$ on ${\Cal S}_F$
with the norm at most $C2^{-N},$ if $N$ is chosen large enough, then we obtain
$$ f(x,y) =\sum \limits_j\sum \limits_k\sum \limits_J\sum \limits_I
\left[\sum\limits_{i=0}\limits^{\infty}{\Cal R}^i\int \limits_J\int \limits_I
\psi_{j,k}(\cdot-u,\cdot-v)
dudv\right](x,y)\left(\psi _{j,k}*f\right)(x_I, y_J).$$

Set $$\left[\sum\limits_{i=0}\limits^{\infty}{\Cal R}^i\int \limits_J\int \limits_I
\psi_{j,k}(\cdot-u,\cdot-v)
dudv\right](x,y) = \vert I\vert \vert J\vert{\widetilde\psi}_{j,k}(x,y, x_I,y_J).$$
 It remains to show
${\widetilde\psi}_{j,k}(x,y,x_I,y_J)\in {\Cal S}_F.$ This, however, follows easily from (3.15).

\hfill {\bf Q.E.D}

We next establish a relationship between test function in $f^\sharp\in {\Cal S}_{\infty}(R^{n+m}\times R^m)$ and test function
$f(x,y) =\int\limits
_{R^m}f^\sharp(x,y-z,z)dz$ under the actions of ${\Cal R}^\sharp$, ${\Cal R}$ and the implicit multi-parameter dilations.

To be more precise,
we  define $$f_{j,k}(x,y)=\int\limits_{R^m}
(f^\sharp)_{j,k}(x,y-z,z)dz,$$
 where
 $$(f^\sharp)_{j,k}(x,y,z)=
2^{(n+m)j}2^{mk}f^\sharp(2^jx,2^jy,2^kz).$$
 Note that
 $$
 \align
 & {\Cal R}^\sharp
((f^\sharp)_{j,k})(x,y,z)\\ = &
\int\int\int R^\sharp(x,y,z, u^\prime, v^\prime, w^\prime)(f^\sharp)_{jk}(u^\prime, v^\prime, w^\prime)du^\prime dv^\prime dw^\prime\\ =&
2^{(n+m)j}2^{mk}\int\int\int
{\Cal R}^\sharp(x,y,z, u^\prime, v^\prime, w^\prime)f^\sharp(2^ju^\prime, 2^jv^\prime, 2^k w^\prime)du^\prime dv^\prime dw^\prime\\ = &
\int\int\int
{\Cal R}^\sharp(x,y,z, 2^{-j}u^\prime, 2^{-j}v^\prime, 2^{-k}w^\prime)f^\sharp(u^\prime, v^\prime, w^\prime)du^\prime dv^\prime dw^\prime
\endalign
$$
and
$${\Cal R}^\sharp(2^{-j}x,2^{-j}y,2^{-k}z, 2^{-j}u^\prime, 2^{-j}v^\prime, 2^{-k}w^\prime)=2^{(n+m)j}2^{mk}{\Cal R}^\sharp(x,y,z, u^\prime, v^\prime, w^\prime).$$
Thus we have
$$
{\Cal R}^\sharp(
(f^\sharp)_{j,k})(x,y,z)=({\Cal R}^\sharp(f^\sharp))_{jk}(x,y,z).$$
 This implies
$$
\align
{\Cal R}(f_{j,k})(x,y)= & \int_{R^m}{\Cal R}^\sharp((f^\sharp)_{j,k})(x,y-z,z)dz\\ = &
\int_{R^m}({\Cal R}^\sharp(f^\sharp))_{jk}(x,y-z,z)dz=({\Cal R}(f))_{jk}(x,y)
\endalign
$$

 It is worthwhile to point out that
$${\widetilde\psi}_{j,k}(x,y,x_I,y_J)=\int\limits_{R^m}{\widetilde\psi}^
\sharp_{j,k}(x,y-z,z,x_I,y_J)dz,$$
 where
$${\widetilde\psi}^\sharp_{j,k}(x,y,z,x_I,y_J)=({\widetilde\psi}^\sharp)_{j,k}(x,y,z,x_I,y_J),\,\,
{\widetilde\psi}^\sharp(x,y,z,x_I,y_J)\in
{\Cal S}_{\infty}(R^{n+m}\times R^m)$$
 and satisfies
the condition
in (3.6) with $x_0=x_I, y_0=y_J.$

{\bf Remark 3.1:}
If we begin with discretizing (3.12) by
$$ f(x,y) =\sum \limits_j\sum \limits_k\sum \limits_J\sum \limits_I
\psi_{j,k}(x-x_I,y-y_J)\int \limits_J\int \limits_I
\left(\psi _{j,k}*f\right)(u,v)dudv + {\widetilde {\widetilde {\Cal R}}}(f)(x, y),$$
and repeating the similar proof, then the discrete Calder\'on reproducing formula can also be given by
the following form
$$ f(x,y) =\sum \limits_j\sum \limits_k\sum \limits_J\sum
\limits_I \vert I\vert \vert J\vert \psi_{j,k}(x-x_I,y-y_J){\widetilde {\widetilde \psi }}_{j,k}
(f)(x_I,y_J),$$
where $\vert I\vert \vert J\vert{\widetilde {\widetilde \psi }}_{j,k}(f)(x_I,y_J)
=\sum\limits_{i=0}\limits^{\infty}\int \limits_J\int \limits_I
\psi_{j,k}*({\widetilde {\widetilde {\Cal R}}})^i(f)(u,v)dudv$.
We leave the details of these proofs to the reader.

Before we  prove the Plancherel-P\^olya-type inequality, we first prove the following lemma.

\vskip -0.3cm
\proclaim{Lemma 3.7}Let $I, I^\prime$, $J, J^\prime$ be dyadic cubes in $R^n$ and $R^m$ respectively such that
$\ell(I)=2^{-j-N}$, $\ell(J)=2^{-j-N}+2^{-k-N}$, $\ell(I^\prime)=2^{-j^\prime-N}$ and $\ell(J^\prime)=2^{-j^\prime-N}+2^{-k^\prime-N}$. Thus
for any $u\in I$ and $v\in J$ we have
$$
\align
&
\sum\limits_{k^\prime\leq j^\prime}\sum \limits_{I^\prime, J^\prime}
{{2^{-(\vert j-{j^\prime }\vert +\vert k-{k^\prime }\vert) K}
2^{-(j^\prime+k^\prime) K}\vert I^\prime\vert\vert
 J^\prime\vert }\over{\left(2^{-j^\prime }+
\vert u-x_{I^\prime }\vert \right)^{n+K}\left(2^{-k^\prime }+\vert v-y_{J^\prime}\vert
 \right)^{m+K}}} \vert \phi _{{j^\prime },{k^\prime }}*f(x_{I^\prime },y_{J^\prime
})\vert \\
\leq & C\sum \limits_{ k^\prime\leq j^\prime
}2^{-\vert j-{j^\prime }\vert K }\cdot 2^{-\vert k-{k^\prime }\vert K}
\left\lbrace M_s\left(\sum \limits_{J^\prime }\sum \limits_{I^
\prime }\vert \phi _{{j^\prime },{k^\prime }}*f(x_{I^\prime },y_{J^
\prime })\vert \chi _{J^\prime} \chi _{I^\prime }\right)^r\right\rbrace
^{{{1}\over{r}}}(u,v)
\endalign
$$
and
$$
\align
&   \sum\limits_{k^\prime> j^\prime}\sum \limits_{I^\prime, J^\prime}
{{2^{-(\vert j-{j^\prime }\vert +\vert k-{k^\prime }\vert) K}
2^{-2j^\prime K}\vert I^\prime\vert\vert
 J^\prime\vert }\over{\left(2^{-j^\prime }+
\vert u-x_{I^\prime }\vert \right)^{n+K}\left(2^{-j^\prime }+\vert v-y_{J^\prime}\vert
 \right)^{m+K}}} \vert \phi _{{j^\prime },{k^\prime }}*f(x_{I^\prime },y_{J^\prime
})\vert \\ \le & C\sum \limits_{k^\prime > j^\prime
}2^{-\vert j-{j^\prime }\vert K }2^{-\vert k-{k^\prime }\vert K}
\left\lbrace M\left(\sum \limits_{J^\prime }\sum \limits_{I^
\prime }\vert \phi _{{j^\prime },{k^\prime }}*f(x_{I^\prime },y_{J^
\prime })\vert \chi _{I^\prime }\chi _{J^\prime} \right)^r\right\rbrace
^{{{1}\over{r}}}(u,v)
\endalign
$$
where $M$ is the Hardy-Littlewood maximal function on $R^{n+m}$, $M_s$ is the strong maximal
function on $R^n\times R^m$ as defined in (1.1), and
$\max\left\{\frac{n}{n+K}, \frac{m}{m+K} \right\} <r.$
\endproclaim
\vskip -0.3cm
{\bf Proof:}
We set
$$A_0=\{I^\prime: \ell(I^\prime)=2^{-j^\prime-N},\,\, \frac{|u-x_{I^\prime}|}{2^{-j^\prime}}\le 1\}$$
$$B_0= \{J^\prime: \ell(J^\prime)=2^{-j^\prime-N}+2^{-k^\prime-N},\,\, \frac{|v-y_{J^\prime}|}{2^{-k^\prime}}\le 1\}$$
and for $\ell\ge 1$, $i\ge 1$
$$A_\ell=\{I^\prime: \ell(I^\prime)=2^{-j^\prime-N},\,\, 2^{\ell-1}<\frac{|u-x_{I^\prime}|}{2^{-j^\prime}}\le 2^\ell\}.$$
$$
B_i=\{J^\prime: \ell(J^\prime)=2^{-j^\prime-N}+2^{-k^\prime-N},\,\, 2^{i-1}<\frac{|v-y_{J^\prime}|}{2^{-k^\prime}}\le 2^i\}.
$$
Then
$$
\align
&
{{2^{-(j^\prime+k^\prime)K}\vert I^\prime\vert \vert J^\prime\vert}\over{\left(2^{-j^\prime}+\vert u-x_{I^\prime}\vert\right)^{n+K}\left(2^{-k^\prime }+\vert v-y_{J^\prime}\vert\right)^{m+K}}}  \cdot \vert
 \phi_{j^\prime, k^\prime}*f(x_{I^\prime},y_{J^\prime})\vert \\ \le & \sum\limits_{\ell, i\ge 0} 2^{-\ell (n+K)}2^{-i(m+K)}2^{-N(n+m)}
\sum\limits_{I^\prime\in A_\ell, J^\prime\in B_i}\vert \phi _{{j^\prime },{k^\prime }}*f(x_{I^\prime },y_{J^\prime
})\vert \\ \le &  \sum\limits_{\ell, i\ge 0} 2^{-\ell (n+K)}2^{-i(m+K)}2^{-N(n+m)}
\left(\sum\limits_{I^\prime\in A_\ell, J^\prime\in B_i} \left(\vert \phi _{{j^\prime },{k^\prime }}*f(x_{I^\prime },y_{J^\prime})\vert\right)^r\right)^{\frac{1}{r}}\\ = & \sum\limits_{\ell, i\ge 0} 2^{-\ell (n+K) -i(m+K)-N(n+m)}
\left(\int_{\prods}|I^\prime|^{-1}|J^\prime|^{-1}\sum\limits_{I^\prime\in A_\ell, J^\prime\in B_i} \vert \phi _{{j^\prime },{k^\prime }}*f(x_{I^\prime },y_{J^\prime})\vert^r\chi_{I^\prime}\chi_{J^\prime}\right)^{\frac{1}{r}}\\ \le &
\sum\limits_{\ell, i\ge 0} 2^{-\ell (n+K-\frac{n}{r})-i(m+K-\frac{m}{r})+(\frac{1}{r}-1)N(n+m)}
\left(M_s\left(\sum\limits_{I^\prime\in A_\ell, J^\prime\in B_i} \vert \phi _{{j^\prime },{k^\prime }}*f(x_{I^\prime },y_{J^\prime})\vert^r\chi_{I^\prime}\chi_{J^\prime}\right)(u,v)\right)^{\frac{1}{r}}\\ \le & C(N)
\left(M_s\left(\sum\limits_{I^\prime, J^\prime} \vert \phi _{{j^\prime },{k^\prime }}*f(x_{I^\prime },y_{J^\prime})\vert^r\chi_{I^\prime}\chi_{J^\prime}\right)(u,v)\right)^{\frac{1}{r}}
\endalign
$$
The last inequality follows from the assumption that
$r>\frac{n}{n+K}$ and $r>\frac{m}{m+K}$ which can be done by choosing $K$ large enough.

Similarly, we can prove the second inequality of the lemma. \hfill{\bf Q.E.D.}

We now are ready to give the

\noindent {\bf Proof of Theorem 1.9:}
By  Theorem 1.8, $f\in {\Cal S}_F$ can be represented by
$$ f(x,y) =\sum \limits_{j^\prime }\sum \limits_{k^\prime }\sum
\limits_{J^\prime }\sum \limits_{I^\prime }\vert J^\prime\vert\vert I^\prime\vert
{\widetilde \phi }_{{j^\prime },{k^\prime }}(x,y,x_{I^\prime },y_{J^\prime }
)\left(\phi _{{j^\prime },{k^\prime }}*f\right)(x_{I^\prime
},y_{J^\prime }).$$
We write
$$
\align
& \left(\psi _{j,k}*f\right)(u,v)\\= &
\sum \limits_{j^\prime }\sum \limits_{k^\prime
}\sum \limits_{J^\prime }\sum \limits_{I^\prime }\vert I^\prime\vert\vert J^\prime\vert
\left(\psi _{j,k}*{\widetilde \phi }_{{j^\prime },{k^\prime
}}(\cdot,\cdot, x_{I^\prime}, y_{J^\prime})\right)(u,v)\left(\phi _{{j^\prime
},{k^\prime }}*f\right)(x_{I^\prime },y_{J^\prime }).
\endalign
$$

By the almost orthogonality estimates in (3.10) and (3.11), for any given positive integer $K,$ we have
if $j^{\prime}\geq k^{\prime},$
$$
\align
& \vert \psi _{j,k}*{\widetilde \phi }_{{j^\prime },{k^\prime }}(\cdot, \cdot)(u,v)\vert\\
\leq & C2^{-\vert j-{j^\prime }\vert K }\cdot 2^{-\vert k-{k^\prime }\vert K}\cdot {{2^{-j^\prime K}}\over{(2^{-j^\prime }+
\vert u-x_{I^\prime }\vert )^{n+K}}}\cdot {{2^{-k^\prime K}}\over{(2^{-k^\prime }+\vert v-y_{J^\prime}\vert)^{m+K}}}\tag3.15
\endalign
$$
and if $j^{\prime}\leq k^{\prime},$ we have
$$
\align
& \vert \psi _{j,k}*{\widetilde \phi }_{{j^\prime },{k^\prime }}(\cdot, \cdot)(u,v)\vert\\
\leq & C2^{-\vert j-{j^\prime }\vert K }\cdot 2^{-\vert k-{k^\prime }\vert K}\cdot
{{2^{-j^\prime K}}\over{(2^{-j^\prime }+\vert u-x_{I^\prime }\vert )^{n+K}}}\cdot {{2^{-j^\prime K}}\over{(2^{-j^\prime }+\vert v-y_{J^\prime}\vert )^{m+K}}}.\tag3.16
 \endalign
 $$
Using  Lemma 3.7,  for any $u, x_{I^\prime}\in I,
v, y_{J^\prime}\in J,$
$$
\align
&
\vert \psi _{j,k}*f(u,v)\vert\\
\leq & C \sum \limits_{{k^\prime}\leq {j^\prime}
}\sum \limits_{I^\prime }\sum \limits_{J^\prime }
2^{-\vert j-{j^\prime }\vert K }2^{-\vert k-{k^\prime }\vert K}\vert I^\prime\vert\vert
 J^\prime\vert \times\\
& {{2^{-j^\prime K}}\over{(2^{-j^\prime }+
\vert u-x_{I^\prime }\vert )^{n+K}}}\cdot {{2^{-
k^\prime K}}\over{(2^{-k^\prime }+\vert v-y_{J^\prime}\vert
 )^{m+K}}}
\vert \phi _{{j^\prime },{k^\prime }}*f(x_{I^\prime },y_{J^\prime
})\vert \\
+ &  C \sum \limits_{{k^\prime}> {j^\prime}
}\sum \limits_{J^\prime }\sum \limits_{I^\prime }
2^{-\vert j-{j^\prime }\vert K }2^{-\vert k-{k^\prime }\vert K}\vert I^\prime\vert\vert
 J^\prime\vert \times\\ &
{{2^{-j^\prime K}}\over{(2^{-j^\prime }+
\vert u-x_{I^\prime }\vert )^{n+K}}}\cdot {{2^{-
j^\prime K}}\over{(2^{-j^\prime }+\vert v-y_{J^\prime}\vert
 )^{m+K}}}\vert \phi _{{j^\prime },{k^\prime }}*f(x_{I^\prime },y_{J^\prime
})\vert\\
\leq & C\sum \limits_{ k^\prime\leq j^\prime
}2^{-\vert j-{j^\prime }\vert K }\cdot 2^{-\vert k-{k^\prime }\vert K}
\left\lbrace M_s(\sum \limits_{J^\prime }\sum \limits_{I^
\prime }\vert \phi _{{j^\prime },{k^\prime }}*f(x_{I^\prime },y_{J^
\prime })\vert \chi _{J^\prime} \chi _{I^\prime })^r\right\rbrace
^{{{1}\over{r}}}(u,v)\\
+ & C\sum \limits_{k^\prime > j^\prime
}2^{-\vert j-{j^\prime }\vert K }2^{-\vert k-{k^\prime }\vert K}
\left\lbrace M\left(\sum \limits_{J^\prime }\sum \limits_{I^
\prime }\vert \phi _{{j^\prime },{k^\prime }}*f(x_{I^\prime },y_{J^
\prime })\vert \chi _{I^\prime }\chi _{J^\prime} \right)^r\right\rbrace
^{{{1}\over{r}}}(u,v)
\endalign
$$
where $M$ is the Hardy-Littlewood maximal function on $R^{n+m}$, $M_s$ is the strong maximal
function on $R^n\times R^m$, and
$\max\{\frac{n}{n+K}, \frac{m}{m+K}\} <r<p.$

Applying the Holder's inequality and summing over $j, k, I, J$  yields
$$
\align
&
 \left\lbrace \sum \limits_j\sum \limits_k\sum \limits_J\sum \limits_I
\sup_{u\in I,v\in J}\vert \psi _{j,k}*f(u,v)\vert ^2\chi
_I\chi _J\right\rbrace^{{{1}\over{2}}}\\
\leq & C\left\lbrace \sum \limits_{j^\prime }\sum \limits_{k^\prime
}\left\lbrace M_s(\sum \limits_{J^\prime }\sum \limits_{I^\prime
}\vert \phi _{{j^\prime },{k^\prime }}*f(x_{I^\prime },y_{J^\prime
})\vert \chi _{I^\prime }\chi _{J^\prime} )^r\right\rbrace^{{{2}
\over{r}}}\right\rbrace^{{{1}\over{2}}}.
\endalign
$$
Since $x_{I^\prime }$ and $y_{J^\prime }$ are arbitrary points in $I^\prime
$ and $J^\prime ,$ respectively, we have
$$
\align
& \left\lbrace \sum \limits_j\sum \limits_k\sum \limits_J\sum \limits_I
\sup_{u\in I,v\in J}\vert \psi _{j,k}*f(u,v)\vert ^2\chi
_I\chi _J\right\rbrace^{{{1}\over{2}}}\\ \leq & C\left\lbrace \sum \limits_{j^\prime }\sum \limits_{k^\prime }\left\lbrace
M_s(\sum \limits_{J^\prime }\sum \limits_{I^\prime }\inf_{u
\in I',v\in J'}\vert \phi _{{j^\prime },{k^\prime }}*f(u,v)\vert
\chi _{I^\prime }\chi _{J^\prime} )^r\right\rbrace^{{{2}\over{r}}}
\right\rbrace^{{{1}\over{2}}},\,\,
\endalign
$$
and hence, by the Fefferman-Stein vector-valued maximal function
inequality [FS] with $r<p,$ we get
$$
\align
&
\Vert \left\lbrace \sum \limits_j\sum \limits_k\sum \limits_J
\sum \limits_I \sup_{u\in I,v\in J}\vert \psi_{j,k}*f(u,v)
\vert^2\chi_I\chi_J\right\rbrace^{{{1}\over{2}}}\Vert_p\\
\leq &
 C\Vert \left\lbrace \sum \limits_{j^\prime }\sum \limits_{k^
\prime }\sum \limits_{J^\prime }\sum \limits_{I^\prime } \inf_{u
\in I^\prime ,v\in J^\prime }\vert \phi _{{j^\prime },{k^\prime
}}*f(u,v)\vert ^2\chi _{I^\prime }\chi _{J^\prime }\right\rbrace^{{{1}
\over{2}}}\Vert _p.
\endalign
$$
This ends the proof of Theorem 1.9. \hfill{\bf Q.E.D.}

\vskip 0.1cm
\noindent {\bf 4.  Discrete Littlewood-Paley-Stein square function, boundedness of flag singular integrals on  Hardy spaces  $H^p_F$, from $H^p_F$ to $L^p$: Proofs of  Theorems 1.10 and 1.11}
\vskip .1cm

The main purpose of this section is to establish the Hardy space theory associated with the flag multi-parameter structure using the results we have proved in Section 3.
As a consequence of Theorem 1.9, it is easy to see that the Hardy space $H_F^p$ is independent of
the choice of the functions $\psi.$ Moreover, we have
the following characterization of $H^p_F$ using the discrete norm.
\vskip -0.3cm
\proclaim{Proposition 4.1} Let $0<p\le 1$. Then we have
$$\Vert f \Vert_{H_F^p}\approx \Vert \left\lbrace \sum \limits_j\sum
 \limits_k\sum
\limits_J\sum \limits_I \vert \psi _{j,k}*f(x_I,y_J)
\vert ^2\chi _I(x)\chi _J(y)\right\rbrace ^{{{1}\over{2}}}\Vert _p\tag4.1$$
where $j, k, \psi, \chi _I,\chi _J, x_I, y_J$ are same as in Theorem 1.9.
\endproclaim
\vskip -0.3cm
Before we give the proof of the boundedness of flag singular integrals on $H_F^p,$ we show several
properties of $H_F^p.$
\vskip -0.3cm
\proclaim{Proposition 4.2} ${\Cal S}_F(R^n\times R^m)$ is dense in $H_F^p.$
\endproclaim
\vskip -0.3cm
\noindent{\bf Proof:}
Suppose $f\in H_F^p,$
and set $W=\lbrace (j,k,I,J): \vert j\vert\leq L, \vert k\vert \leq M, I\times J\subseteq B(0, r)
\rbrace,$ where $I, J$ are dyadic cubes in $R^n, R^m$ with side length $2^{-j-N}, 2^{-k-N}+2^{-j-N},$
respectively, and $B(0,r)$ are balls in $R^{n+m}$ centered at the origin with radius $r$. It is
easy to see that $$ \sum \limits_{(j,k,I,J)\in W} \vert I\vert \vert J\vert {\widetilde \psi
 }_{j,k}(x,y, x_I,y_J)\psi _{j,k}*f(x_I,y_J)$$
  is a test function in ${\Cal S}_F(R^n\times R^m)$ for any fixed
$L, M, r.$ To show the proposition, it suffices to
prove $$ \sum \limits_{(j,k,I,J)\in W^c} \vert
 I\vert \vert J\vert {\widetilde \psi
 }_{j,k}(x,y,x_I,y_J)\psi _{j,k}*f(x_I,y_J)$$ tends to zero in the $H_F^p$ norm as $L, M, r$ tend to infinity.
This follows from (4.1) and a similar proof
as in the proof of Theorem 1.9. In fact, repeating the same proof as in theorem 1.9 yields
$$\Vert \sum \limits_{(j,k,I,J)\in W^c} \vert I\vert \vert J\vert {\widetilde \psi
 }_{j,k}(x,y,x_I,y_J)\psi _{j,k}*f(x_I,y_J)\Vert_{H_F^p}$$
$$\leq C\Vert \lbrace \sum \limits_{(j,k,I,J)\in W^c}
\vert \psi _{j,k}*f(x_I,y_J)
\vert ^2\chi _I\chi _J\rbrace ^{{{1}\over{2}}}\Vert _p,$$
where the last term tends to zero as $L, M, r$ tend to infinity whenever $f\in H_F^p.$ \hfill{\bf Q.E.D.}

As a consequence of Proposition 4.2, $L^2(R^{n+m})$ is dense in $H_F^p(R^n\times R^m).$ Furthermore, we have

\proclaim{Theorem 4.3} If $f\in L^2(R^{n+m})\cap \phpf, 0<p\leq 1,$ then $f\in L^p(R^{n+m})$ and
there is a constant $C_p>0$ which is independent of the $L^2$ norm of $f$ such that $$\Vert f\Vert_p\leq C\Vert f\Vert_{H_F^p}.\tag4.3$$
\endproclaim
\vskip -0.7cm

To show theorem 4.3, we need a discrete Calder\'on reproducing formula on $L^2(R^{n+m})$. To be more precise, take $\phi^{(1)}\in C_0^\infty(R^{n+m})$ with
$$\int_{R^{n+m}}\phi^{(1)}(x,y)x^\alpha y^\beta dxdy=0, \,\, \text{for all}\, \alpha, \beta\,\, \text{satisfying}\, 0\le |\alpha|\le M_0, \, 0\le|\beta|\le M_0,$$
where $M_0$ is a large positive integer which will be determined later,
and
$$\sum_{j}|\widehat{\phi^{(1)}}(2^{-j}\xi_1, 2^{-j}\xi_2)|^2=1,\,\, \text{for all}\, (\xi_1,\xi_2)\in R^{n+m}\backslash\{(0,0)\}, $$
and take $\phi^{(2)}\in C_0^\infty(R^m)$ with
$$\int_{R^m}\phi^{(2)}(z)z^\gamma dz=0\,\, \text{for all}\,\, 0\le |\gamma|\le M_0,$$
and $\sum_k|\widehat{\phi^{(2)}}(2^{-k}\xi_2)^2=1$ for all $\xi_2\in R^m\backslash\{0\}$.

Furthermore, we may assume that $\phi^{(1)}$ and $\phi^{(2)}$ are radial functions and supported in the unit balls of $R^{n+m}$ and $R^m$ respectively. Set
again
$$\phi_{jk}(x,y)=\int_{R^m}\phi^{(1)}_j(x, y-z)\phi_k^{(2)}(z)dz.$$

By taking the Fourier transform, it is easy to see the following continuous version of Calder\'on reproducing formula on $L^2$: for $f\in L^2(R^{n+m})$,
$$f(x,y)=\sum \limits_j \sum\limits_k \phi_{jk}*\phi_{jk}*f(x,y).$$
\vskip -0.5cm
For our purpose, we need the discrete version of the above reproducing formula.
\proclaim{Theorem 4.4}
There exist functions $\widetilde {\phi}_{jk}$ and an operator $T_N^{-1}$ such that
 $$f(x,y)=\sum \limits_j\sum \limits_k\sum \limits_J\sum \limits_I
\vert I\vert \vert J\vert\widetilde {\phi}_{j,k}(x-x_I, y- y_J)
\phi _{j,k}*\left (T_N^{-1}(f)\right )(x_I, y_J)$$
where functions $\widetilde{\phi}_{jk}(x-x_I,y-y_J)$ satisfy the conditions in (3.6) with $\alpha_1, \beta_1, \gamma_1, N, M$ depending on $M_0$, $x_0=x_I$ and $y_0=y_J$. Moreover,  $T_N^{-1}$ is bounded on $L^2(R^{n+m})$ and $\phpf,$ and the series converges in $L^2(R^{n+m})$.
\endproclaim
\vskip -0.4cm
{\bf Remark 4.1:} The difference between Theorem 4.4 and Theorem 1.8 are that our ${\widetilde \phi}_{jk}$ in Theorem 4.4 has
compact support. The price we pay here is that ${\widetilde \phi}_{jk}$ only satisfies the  moment condition of finite order, unlike
that in Theorem 1.8 where the moment condition of infinite order is satisfied. Moreover, the formula in Theorem 4.4 only
 holds on $L^2(R^{n+m})$ while the formula in Theorem 1.18 holds in test function space ${\Cal S}_F$ and its dual space $({\Cal S}_F)^{\prime}$.

\noindent
{\bf Proof of Theorem 4.4:}
Following the proof of Theorem 1.8, we have
$$f(x,y)=\sum \limits_j\sum \limits_k\sum \limits_J\sum \limits_I
[\int \limits_J\int \limits_I
\phi_{j,k}(x-u,y-v)
dudv]\left(\phi _{j,k}*f\right)(x_I, y_J) + {\Cal R}(f)(x, y).\tag4.4$$
where $I, J, j, k$ and ${\Cal R}$ are the same as in Theorem 1.8.

We need the following

\proclaim{Lemma 4.5} Let $0<p\le 1$. Then
 the operator ${\Cal R}$ is bounded on $L^2(R^{n+m})\cap \phpf$ whenever $M_0$ is chosen to be a large positive integer. Moreover, there exists a constant $C>0$
such that
$$||{\Cal R}(f)||_2\le C2^{-N} ||f||_2$$
and
$$||{\Cal R}(f)||_{\phpf}\le C2^{-N} ||f||_{\phpf}.$$
\vskip -0.4cm
\endproclaim
{\bf Proof of Lemma 4.5:} Following the proofs of Theorems 1.8 and 1.9 and using the discrete Calder\'on reproducing formula
for $f \in L^2(R^{n+m}),$
we have
$$
\align
&
||g_F({\Cal R}(f))||_p \\ \le & \Vert  \left\lbrace \sum \limits_j\sum\limits_k\sum\limits_J\sum \limits_I \vert \left(\psi _{j,k}*{\Cal R}(f)\right)
\vert ^2\chi _I\chi _J\right\rbrace ^{{{1}\over{2}}}\Vert _p\\
= & \left\Vert\left\lbrace \sum \limits_{j,k, J,I}
\sum \limits_{{j^\prime},k^\prime, J^\prime, I^\prime}|J^\prime||I^\prime|
\vert\left(\psi _{j,k}*{\Cal R}\left(\widetilde{\psi_{j^\prime,k^\prime}}(\cdot, x_{I^\prime}, \cdot, y_{J^\prime})\cdot \psi_{j^\prime k^\prime}*f(x_{I^\prime}, y_{J^\prime})\right)\right)\vert ^2\chi _I\chi _J\right\rbrace ^{{{1}\over{2}}}\right\Vert _p
\endalign
$$
where $j, k, \psi, \chi _I,\chi _J, x_I, y_J$ are the same as in Theorem 1.9.

We {\bf claim:}
$$
\align &
\vert \left(\psi _{j,k}*{\Cal R}\left({\widetilde\psi} _{{j^\prime },{k^\prime }}(\cdot, x_{I^
\prime },\cdot, y_{J^\prime })\right)\right)(x,y)\vert\\
\leq  & C2^{-N}2^{-\vert j-{j^
\prime }\vert K }2^{-\vert k-{k^\prime  }\vert K }\cdot
\int\limits_{R^m}{{2^{-(j\wedge {j^
\prime })K}}\over{(2^{-(j\wedge {j^\prime })}+\vert
x-x_{I^\prime }\vert +\vert y-z-y_{J^\prime} )^{n+m+K}}}{{2^{-(k\wedge {k^\prime
})K}}\over{(2^{-(k\wedge {k^\prime })}+\vert z
\vert )^{m+K}}}dz
\endalign
$$
 where $K<M_0, max({{n}\over{n+K}}, {{m}\over{m+K}})<p,$ and $M_0$ is chosen to be a lager integer later.

Assuming the claim for the moment,  repeating a similar proof in Lemma 3.7 and then Theorem 1.9, we obtain
$$
\align &
\Vert \vert g_F({\Cal R}f)\Vert _p \leq  C2^{-N}\Vert
\lbrace \sum \limits_{j^\prime }\sum \limits_{k^\prime }\lbrace
M_s(\sum \limits_{J^\prime }\sum \limits_{I^\prime }\vert \psi _{{j^\prime },{k^\prime }}*f(x_{I^\prime },y_{J^\prime })\vert \chi _{J^\prime} \chi _{I^\prime })^r\rbrace^{{{2}\over{r}}}
\rbrace^{{{1}\over{2}}}\Vert_p\\ \leq  & C2^{-N}\Vert \lbrace \sum \limits_{j^\prime }\sum \limits_{k^\prime }\sum \limits_{J^\prime }\sum \limits_{I^\prime }\vert \psi_{{j^\prime },{k^\prime }}*f(x_{I^\prime },y_{J^\prime })\vert ^2
\chi _{I^\prime }\chi_{J^\prime } \rbrace^{{{1}\over{2}}}\Vert_p\leq C2^{-N}\Vert  f \Vert_{\phpf}.
\endalign
$$
It is clear that the above estimates still hold when $p$ is replaced by 2. These imply the assertion of Lemma 4.5.

We now prove the Claim. Again, by the proof of Theorem 1.8,
$$
{\Cal R}\left({\widetilde\psi} _{{j^\prime },{k^\prime }}(\cdot, x_{I^
\prime },\cdot, y_{J^\prime })\right)(x,y)=\int_{R^m}{\Cal R}^\sharp (x,y-z,z, u^\prime, v^\prime, w^\prime){\widetilde\psi} _{{j^\prime },{k^\prime }}(\cdot, x_{I^\prime },\cdot, y_{J^\prime })du^\prime dv^\prime dw^\prime dz$$
where ${\Cal R}^\sharp(x,y,z, u^\prime, v^\prime, w^\prime)$ is similar to ${\Cal R}^\sharp$ as given in the proof of
Theorem 1.8 but, as we pointed out in Remark 4.1,
that the difference between ${\Cal R}^\sharp$ here and ${\Cal R}^\sharp$ given in the proof of Theorem 1.8 is the moment
conditions. However, the almost orthogonality estimate still holds if we only require sufficiently high order of moment
conditions. More precisely, if we replace the moment conditions in (3.5) "for all $\alpha_1,\beta_1,\gamma_1, \alpha_2,
\beta_2, \gamma_2$" by "for all $|\alpha_1|, |\beta_1|, |\gamma_1|, |\alpha_2|, |\beta_2|, |\gamma_2|\le M_0$
where $M_0$ is a large integer, then the estimate in (3.7) still holds with $L_1, L_2, K_1, K_2$ depending on $M_0$.
Thus, the claim follows by applying the same proof as that of Theorem 1.8, and the proof of Lemma 4.5 is complete.\hfill{\bf Q. E. D.}

We now return to the proof of Theorem 4.4.

Denote $(T_N)^{-1}=\sum_{i=1}^{\infty}{\Cal R}^i,$ where
$$T_N(f)=\sum\limits_j\sum\limits_k\sum\limits_J\sum\limits_I[\int\limits_J\int\limits_I\phi_{j,k}(x-u, y-v)dudvd]\left(
\phi_{j,k}*f\right)(x_I, y_J).$$
\vskip -0.4cm
Lemma 4.5 shows that if $N$ is large enough, then both of $T_N$ and $(T_N)^{-1}$ are bounded on $L^2(R^{n+m})\cap \phpf$. Hence, we can get the following reproducing formula
 $$f(x,y)=\sum \limits_j\sum \limits_k\sum \limits_J\sum \limits_I
\vert I\vert \vert J\vert {\widetilde \phi}_{j,k}(x-x_I,y-y_J)
\phi _{j,k}*\left (T_N^{-1}(f)\right )(x_I, y_J)$$
where ${\widetilde \phi}_{jk}
(x-x_I,y-y_J)={\frac {1}{\vert I\vert}}{\frac {1}{\vert J\vert}}\int\limits_J\int\limits_I\phi_{jk}
(x-x_I-(u-x_I), y-y_J-(v-y_J))dudv$ satisfies the estimate in (3.6) and the series converges in $L^2(R^{n+m})$.

This completes the proof of Theorem 4.4. \hfill {\bf Q.E.D.}

As a consequence of Theorem 4.4, we obtain the following

\proclaim{Corollary 4.6} If $f\in L^2(R^{n+m})\cap \phpf$ and $0<p\le 1$, then
$$
\align
\Vert f\Vert_{H_F^p}\approx  & \Vert \lbrace (\sum\limits_j\sum\limits_k\sum\limits_J\sum\limits_I \vert\phi_{jk}*\left (T_N^{-1}(f)\right )(x_I,y_J)\vert^2\chi_I(x)
\chi_J(y)\rbrace )^{\frac {1}{2}}\Vert_p
\endalign
$$
where the constants are independent of the $L^2$ norm of $f.$
\endproclaim
\vskip -0.4cm
To see the proof of Corollary 4.6, note that if $f\in L^2(R^{n+m}),$ one can apply the Caldero\'n reproducing formulas in Theorem 1.8 and 4.4 and then repeat the same proof as in Theorem 1.9.
We leave the details to the reader. We now start the

\noindent{\bf Proof of Theorem 4.3:} We define a square function by
$$\widetilde{g}(f)(x,y)=\lbrace \sum \limits_j\sum \limits_k\sum
\limits_J\sum \limits_I \vert \phi_{j,k}*\left (T_N^{-1}(f)\right )(x_I,y_J)
\vert ^2\chi _I(x)\chi _J(y)\rbrace ^{{{1}\over{2}}}$$
where $\phi_{jk}$ are the same as in Theorem 4.4.
By Corollary 4.6, for $f\in L^2(R^{n+m})\cap \phpf$ we have,
 $$||\widetilde{g}(f)||_{L^p(R^{n+m})}\le C||f||_{\phpf}$$.

\vskip -0.3cm
To complete the proof of Theorem 4.3, let $f\in L^2(R^{n+m})\cap \phpf.$
Set $$\Omega_i=\lbrace (x,y)\in R^n\times R^m: \widetilde{g}(f)(x,y)
>2^i\rbrace.$$
 Denote $${\Cal B}_i=\lbrace (j,k,I,J): \vert (I\times J)\cap\Omega_i\vert>{{1}\over
 {2}}\vert I\times J\vert, \vert (I\times J)\cap\Omega_{i+1}\vert\leq{{1}\over
 {2}}\vert I\times J\vert\rbrace,$$ where $I, J$ are dyadic cubes in $R^n, R^m$ with side length
$2^{-j-N}
, 2^{-k-N}+2^{-j-N},$ respectively. Since $f\in L^2(R^{n+m})$, by the discrete Calder\'on reproducing formula in Theorem 4.4,
$$
\align
& f(x,y)\\ = & \sum \limits_j\sum \limits_k\sum \limits_J\sum
\limits_I {\widetilde \phi} _{j,k}(x-x_I,y-y_J)\vert I\vert \vert J\vert
\phi_{j,k}*\left (T_N^{-1}(f)\right )(x_I,y_J)\\
= & \sum\limits_i\sum\limits_{(j,k,I,J)\in {\Cal B}_i} \vert I\vert \vert J\vert {\widetilde \phi} _{j,k}(x-x_I,y-y_J)
\phi_{j,k}*\left (T_N^{-1}(f)\right )(x_I,y_J),
\endalign
$$
where the series converges in $L^2$ norm, and hence it also converges almost everywhere.

We {\bf claim}
$$\Vert\sum\limits_{(j,k,I,J)\in {\Cal B}_i} \vert I\vert \vert J\vert {\widetilde \phi} _{j,k}(x-x_I,y-y_J)\phi_{j,k}*\left (T_N^{-1}(f)\right )
(x_I,y_J)\Vert_p^p \leq C 2^{ip}\vert \Omega_i\vert,$$
which together with the fact $0<p\le 1$ yields
$$
\align
||f||_p^p\le & \sum_i \Vert\sum\limits_{(j,k,I,J)\in {\Cal B}_i} \vert I\vert \vert J\vert {\widetilde \phi} _{j,k}(x-x_I,y-y_J)
\phi_{j,k}*\left (T_N^{-1}(f)\right )(x_I,y_J)\Vert_p^p\\ \le & C\sum_i 2^{ip}|\Omega_i|\\ \le & C||\widetilde{g}(f)||_p^p
\le C ||f||_{H^p_F}^p.
\endalign
$$

To show the claim,
note that $\phi^{(1)}$ and $\psi^{(2)}$ are radial functions supported in unit balls.
Hence, if
$(j,k,I,J)\in {\Cal B}_i$ then $\phi _{j,k}(x-x_I,y-y_J)$ are supported in
$${\widetilde{\Omega_i}}=
\lbrace (x,y): M_s(\chi_{\Omega_i})(x,y)>{{1}\over {100}}\rbrace.$$
 Thus, by H\"older's inequality,
$$
\align &
\Vert \sum\limits_{(j,k,I,J)\in {\Cal B}_i}\vert J\vert \vert I\vert {\widetilde \phi} _{j,k}(x-x_I,y-y_J)
\phi_{j,k}*\left (T_N^{-1}(f)\right )(x_I,y_J)\Vert_p^p\\
\leq & \vert {\widetilde{\Omega_i}}\vert^{1-{{p}\over{2}}}
\Vert \sum\limits_{(j,k,I,J)\in {\Cal B}_i}\vert J\vert \vert I\vert {\widetilde \phi}_{j,k}(x-x_I,y-y_J)
\phi_{j,k}*\left (T_N^{-1}(f)\right )(x_I,y_J)\Vert_2^p.
\endalign
$$
By the duality argument, for all $g\in L^2$ with $\Vert g\Vert_2\leq 1,$
$$
\align
& \vert <\sum\limits_{(j,k,I,J)\in {\Cal B}_i}\vert J\vert \vert I\vert {\widetilde \phi}_{j,k}(x-x_I,y-y_J)
\phi_{j,k}*\left (T_N^{-1}(f)\right )(x_I,y_J), g>\vert\\
= & \vert\sum\limits_{(j,k,I,J)\in {\Cal B}_i}\vert J\vert \vert I\vert {\widetilde \phi}_{j,k}*g(x_I,y_J)
\phi_{j,k}*\left (T_N^{-1}(f)\right )(x_I,y_J)\vert\\
\le  & C  \left(\sum\limits_{(j,k,I,J)\in {\Cal B}_i}|I||J|\vert
\phi_{j,k}*\left (T_N^{-1}(f)\right )(x_I,y_J)\vert^2\right)^{\frac{1}{2}}
\cdot  \left(\sum\limits_{(j,k,I,J)\in {\Cal B}_i}|I||J| \vert{\widetilde \phi}_{j,k}*g(x_I,y_J)\vert^2\right)^{\frac{1}{2}}.
\endalign
$$
Since
$$
\align
& \left(\sum\limits_{(j,k,I,J)\in {\Cal B}_i}|I||J| \vert {\widetilde \phi} _{j,k}*g(x_I,y_J)\vert^2\right)^{\frac{1}{2}}\\
\le & \left(\sum\limits_{(j,k,I,J)\in {\Cal B}_i}|I||J|\left( M_s \left({\widetilde \phi}_{j,k}*g\right )(x,y)\chi_I(x)\chi_J(y)\right)^2\right)^{\frac{1}{2}}\\
\le & C \left(\sum_{j,k}\int_{R^n}\int_{R^m}\left(M_s\left({\widetilde \phi} _{j,k}*g\right)^2(x,y)dxdy\right)\right)^{\frac{1}{2}}\le C||g||_2
\endalign
$$
thus
the claim now follows from the fact  that $\vert {\widetilde{\Omega_i}}\vert\leq C\vert{\Omega_i}\vert$ and the following
estimate:
$$
\align
C2^{2i}\vert \Omega_i\vert\geq & \int\limits_{{\widetilde{\Omega_i}}\backslash
 {\Omega_{i+1}}}\widetilde{g}^2(f)(x,y)dxdy\\
\geq & \sum\limits_{(j,k,I,J)\in {\Cal B}_i}
\vert \phi_{j,k}*\left (T_N^{-1}(f)\right )(x_I,y_J)\vert^2\vert (I\times J)\cap
 {{\widetilde{\Omega_i}}\backslash
 {\Omega_{i+1}}}\vert\\
\geq &  {{1}\over
 {2}}\sum\limits_{(j,k,I,J)\in {\Cal B}_i}\vert I\vert
 \vert J\vert \vert \phi_{j,k}*\left (T_N^{-1}(f)\right )(x_I,y_J) \vert^2,
 \endalign
 $$
where the fact that $\vert (I\times J)\cap {{\widetilde{\Omega_i}}\backslash
 {\Omega_{i+1}}}\vert>{{1}\over
 {2}}\vert I\times J\vert$ when $(j,k,I,J)\in {{\Cal B}_i}$ is used in the last inequality.
This finishes the proof of Theorem 4.3. \hfill{\bf Q.E.D.}

As a consequence of Theorem 4.3, we have the following

\noindent{\bf Corollary 4.7}
$H_F^1(R^n\times R^m)$ is a subspace
of $L^1(R^n\times R^m).$

\noindent{\bf Proof:} Given $f\in H_F^1(R^{n+m}),$ by Proposition 4.2, there is a
sequence $\lbrace f_n\rbrace$
such that $f_n\in {L^2(R^{n+m})}\cap {H_F^1(R^{n+m})}$ and $f_n$ converges to $f$ in the norm of $H_F^1(R^{n+m}).$
By Theorem 4.3, $f_n$ converges to $g$ in $L^1(R^{n+m})$ for some $g\in L^1(R^{n+m}).$ Therefore, $f = g$ in
$({\Cal S}_F)^\prime.$\hfill {\bf Q.E.D.}

We now turn to the

\noindent {\bf Proof of Theorem 1.10:}
We assume that $K$ is the kernel of $T.$ Applying
the discrete Calder\'on
reproducing formula in Theorem 4.4 implies that for $f\in L^2(R^{n+m})\cap \phpf$,
$$\Vert \lbrace \sum \limits_j\sum \limits_k\sum \limits_J
\sum \limits_I\vert \phi _{j,k}*K*f(x,y)\vert ^2\chi _I(x)\chi _J(y)
\rbrace^{{{1}\over{2}}}\Vert _p = $$
$$\Vert \lbrace \sum \limits_j\sum \limits_k\sum \limits_J
\sum \limits_I\vert \sum \limits_{j^\prime }\sum \limits_{k^\prime
}\sum \limits_{J^\prime }\sum \limits_{I^\prime }\vert J^\prime\vert \vert I^\prime \vert
\phi _{j,k}*K* {\widetilde \phi}_{{j^\prime },{k^\prime
}}(\cdot- x_{I^\prime },\cdot- y_{J^\prime })(x,y)\times$$
$$\phi_{{j^\prime },{k^\prime }}*\left (T_N^{-1}(f)\right )(x_{I^\prime },y_{J^\prime
})\vert ^2\chi _I(x)\chi _J(y)\rbrace^{{{1}\over{2}}}\Vert
_p,$$
where the discrete Calder\'on reproducing formula in $L^2(R^{n+m})$ is used.

Note that $\phi_{jk}$ are dilations of bump functions, by estimates similar to the those in (2.5), one can easily  check that
$$\vert \phi_{j,k}*K* {\widetilde \phi}_{j^\prime, k^\prime }(\cdot- x_{I^\prime },\cdot- y_{J^\prime })(x,y)\vert \leq C 2^{-\vert j-{j^\prime }\vert K}2^{-\vert k-{k^\prime  }\vert K }$$
$$\int\limits_{R^m}{{2^{-(j\wedge {j^\prime })K}}\over{(2^{-(j\wedge {j^\prime })}+\vert
x-x_{I^\prime }\vert +\vert y-z-y_{J^\prime}\vert )^{n+m+K}}}\cdot {{2^{-(k\wedge {k^\prime
})K}}\over{(2^{-(k\wedge {k^\prime })}+\vert z
\vert )^{m+K}}}dz,$$
where $K$ depends on $M_0$ given in Theorem 4.4 and $M_0$ is chosen to be large enough.
Repeating a similar proof in Theorem 1.9 together with Corollary 4.6, we obtain
$$\Vert Tf\Vert _{H_F^p}\leq C\Vert
\lbrace \sum \limits_{j^\prime }\sum \limits_{k^\prime }\lbrace M_s(\sum \limits_{J^\prime }\sum \limits_{I^\prime }\vert \phi_{{j^
\prime },{k^\prime }}*\left (T_N^{-1}(f)\right )(x_{I^\prime },y_{J^\prime })\vert \chi_{J^\prime} \chi _{I^\prime })^r\rbrace^{{{2}\over{r}}}(x,y)
\rbrace^{{{1}\over{2}}}\Vert_p$$
$$\leq C\Vert \lbrace \sum \limits_{j^\prime }\sum \limits_{k^\prime }\sum \limits_{J^\prime }\sum \limits_{I^\prime }\vert \phi_{{j^\prime },{k^\prime }}*\left (T_N^{-1}(f)\right )(x_{I^\prime },y_{J^\prime })\vert ^2
\chi_{J^\prime }(y)\chi _{I^\prime }(x)\rbrace^{{{1}\over{2}}}
\Vert_p\leq C\Vert  f \Vert_{H_F^p},$$
where the last inequality follows from Corollary 4.6.

Since $L^2(R^{n+m})$ is dense in $\phpf,$ $T$ can extend to a bounded operator on $\phpf$. This ends the proof of Theorem 1.10.

{\bf Proof of Theorem 1.11} We note that $H^p_F\cap L^2$ is dense in $H^p_F$, so we only have to show this for $f\in H^p_F\cap L^2$. Thus Theorem 1.11 follows from Theorems 4.3 and 1.10 immediately.
\hfill {\bf Q.E.D.}

\vskip 0.2cm

\noindent {\bf 5. Duality of Hardy spaces $H_F^p$ and boundedness of flag singular integrals on $BMO_F$ space: Proofs of Theorems 1.14, 1.16 and 1.18 }
\vskip .2cm

This section deals with the duality theory of flag Hardy spaces $\phpf$ for all $0<p\le 1$.
We first prove Theorem 1.14, the Plancherel-P\^olya inequalities for $CMO^p_F$ space.

\noindent {\bf Proof of Theorem 1.14:} The idea of the proof of this
theorem is, as in the proof of Theorem 1.9, again to use the
discrete Calder\'on reproducing formula and the almost orthogonality
estimate. For the reader's convenience we  choose to present the
proof of Theorem 1.14 in the case when $n=m=1.$ However, it will be
clear from the proof that its extension to general $n$ and $m$ is
straightforward. Moreover, to simplify notation, we denote
$f_{j,k}=f_R$ when $R=I\times J\subset R^2$ and $\vert
I\vert=2^{-j-N}, \vert J\vert=2^{-k-N}+2^{-j-N}$ are dyadic intervals respectively. Here
$N$ is the same as in Theorem 1.8. We also denote by $dist(I, I')$
the distance between intervals $I$ and $I'$,
$$S_{R}=\sup\limits_{u\in
I,v\in J}\vert \psi _{R}*f(u,v)\vert^2,\,\,\, T_{R}=\inf\limits_{u\in I,v\in J}\vert \phi _{R}*f(u,v)\vert^2.$$

With these notations, we can rewrite the discrete Calder\'on reproducing formula in (1.9) by
$$ f(x,y) =\sum \limits_{R=I\times J}
\vert I\vert \vert J\vert {\widetilde \phi }_{R}(x,y)
\phi _{R}*f(x_I,y_J),$$
where the sum runs over all rectangles $R=I\times J.$

Let $R'=I'\times J', \vert I'\vert =2^{-j^\prime-N}, \vert J'\vert=2^{-j^\prime-N}+2^{-k^\prime-N}, j^\prime>k^\prime$.
Applying the above discrete Calder\'on reproducing formula and the estimates in Corollary 3.5
yields for all $(u,v)\in R,$
$$\vert \psi _{R}*f(u,v)\vert^2\leq C\sum\limits_{R'=I'\times J',  j'> k'}
({{\vert I\vert}\over{\vert
 I'\vert}}
\wedge {{\vert I'\vert}\over{\vert I\vert}})
^{L}
({{\vert J\vert}\over{\vert J'\vert}}\wedge {{\vert J'\vert}\over{\vert J\vert}})
^{L}\times $$
$${{\vert I'\vert}^{K}\over{({\vert
 I'\vert}+\vert
u-x_{I'}\vert )^{(1+K)}}} {{\vert
 J'\vert}^{K}\over{({
\vert J'\vert}+\vert v-y_{J'}\vert
 )^{(1+K)}}}\vert I'\vert \vert J'\vert
\vert\phi _{R'}*f(x_{I'},y_{J'})\vert^2$$
$$+C\sum\limits_{R'=I'\times J', j'\leq k'}
({{\vert I\vert}\over{\vert
 I'\vert}}
\wedge {{\vert I'\vert}\over{\vert I\vert}})
^{L}
({{\vert J\vert}\over{\vert J'\vert}}\wedge {{\vert J'\vert}\over{\vert J\vert}})
^{L}\times $$
$${{\vert I'\vert}^{K}\over{({\vert
 I'\vert}+\vert
u-x_{I'}\vert )^{(1+K)}}} {{\vert
 I'\vert}^{K}\over{({\vert
 I'\vert}+\vert
v-y_{J'}\vert )^{(1+K)}}}\vert I'\vert \vert J'\vert
\vert\phi _{R'}*f(x_{I'},y_{J'})\vert^2,$$
where $K, L$ are any positive integers which can be chosen by $L, K>
{{2}\over {p}} -1$( for general $n,m, K$ can be chosen by $K>(n\vee m)({{2}\over {p}} -1))$,
the constant $C$ depends only on $K, L$ and
functions $\psi$ and $\phi$, here  $x_{I'}$ and $y_{J'}$ are
any fixed points in $I', J',$ respectively.

Adding up all the terms with multiplying $\vert I\vert\vert J\vert$ over $R\subseteq \Omega,$
we obtain
$$\sum\limits_{R\subseteq \Omega}\vert I\vert\vert J\vert S_{R}
\leq C \sum\limits_{R\subseteq
 \Omega}\sum\limits_{R'}\vert I'\vert\vert J'\vert
r(R,R')P(R,R')T_{R'},\tag5.1$$
where $$r(R,R')=({{\vert I\vert}\over{\vert
 I'\vert}}
\wedge {{\vert I'\vert}\over{\vert I\vert}})
^{L-1}({{\vert J\vert}\over{\vert
 J'\vert}}
\wedge {{\vert J'\vert}\over{\vert J\vert}})
^{L-1}$$
 and
$$P(R,R')={{1}\over {(1+{{dist(I,I')}\over{ \vert
 I'\vert}})^{1+K}
(1+{{dist(J,J')}\over{ \vert J'\vert}})^{1+K}}}$$
if $j'>k'$, and
$$P(R,R')={{1}\over {(1+{{dist(I,I')}\over{ \vert
 I'\vert}})^{1+K}
(1+{{dist(J,J')}\over{ \vert I'\vert}})^{1+K}}}$$
if $j'\leq k'$.

To estimate the right-hand side in the above inequality (5.1),
where we first consider $R'=I'\times J', \vert I'\vert =2^{-j^\prime-N}, \vert J'\vert=2^{-j^\prime-N}+2^{-k^\prime-N}, j^\prime>k^\prime$.

 Define
$$\Omega^{i,\ell}=\bigcup_{I\times J\subset \Omega}3(2^i I \times
2^{\ell} J)\,\,\text{for}\,\, i,\ell\geq 0.$$
Let $B_{i,\ell}$ be a collection of
dyadic rectangles $R'$   so that for $i, \ell\ge 1$
$$B_{i,\ell}=\lbrace
R'=I'\times J', 3(2^i I' \times 2^{\ell} J')\bigcap
\Omega^{i,\ell}\ne \emptyset\,\, \text{and}\,\, 3(2^{i-1} I' \times
2^{\ell-1} J')\bigcap \Omega^{i,\ell}=\emptyset \rbrace,$$
 and $$B_{0,\ell}=\lbrace R'=I'\times J', 3(I' \times
2^{\ell} J')\bigcap \Omega^{0,\ell}\ne \emptyset\,\, \text{and}\,\, 3( I'
\times 2^{\ell-1} J')\bigcap \Omega^{0,\ell}=\emptyset \rbrace\,\,\text{for}\,\,
\ell \geq 1,$$
 and
 $$B_{i,0}=\lbrace R'=I'\times J', 3(2^i I' \times
J')\bigcap \Omega^{i,0}\ne \emptyset\,\,  \text{and}\,\, 3(2^{i-1} I' \times
J')\bigcap \Omega^{i,0}=\emptyset \rbrace\,\,\text{for}\,\, i\geq 1,$$
 and
$$B_{0,0}=\lbrace R': R'=I'\times J', 3(I' \times  J')\bigcap
\Omega^{0,0}\ne \emptyset\rbrace.$$
 We write
$$\sum\limits_{R\subseteq
 \Omega}\sum\limits_{R'}\vert I'\vert\vert
 J'\vert r(R,R')P(R,R')T_{R'}
=\sum\limits_{i\geq 0,\ell\geq 0}\sum\limits_{R'\in
{B_{i,\ell}}}\sum\limits_{R\subseteq
 \Omega}\vert I'\vert\vert
 J'\vert r(R,R')P(R,R')T_{R'}.$$

To estimate the right-hand side of the above equality, we first consider the case when $i=\ell=0.$
Note that when $R'\in
 B_{0,0}$, $3R'\bigcap \Omega^{0,0}\ne \emptyset.$
For each integer $h\geq 1,$ let ${\Cal  F}_{h}=\lbrace R'=I'\times
J'\in B_{0,0}, \vert (3I'\times 3J') \bigcap \Omega^{0,0}\vert\geq
({{1}\over{2^h}})\vert 3I'\times 3J'\vert \rbrace.$ Let ${\Cal
D}_{h}={\Cal  F}_h\backslash {\Cal  F}_{h-1},$ and
$\Omega_h=\bigcup_{R'\in {\Cal
 D}_h}
R'.$ Finally, assume that the right-hand side in
(1.12) is finite, that is, for any open set $\Omega\subset R^2,$
$$\sum\limits_{R=I\times J\subseteq \Omega}
\vert I\vert \vert J\vert T_{R}\leq C {\vert \Omega \vert}^{{{2}\over
{p}}-1}.\eqno(5.2)$$

Since $B_{0,0}=\bigcup_{h\geq 1} {\Cal  D}_h$ and for each $R'\in B_{0,0},
P(R,R')\leq 1,$ thus,
$$\sum\limits_{R'\in {B_{0,0}}}\sum\limits_{R\subseteq
\Omega}\vert I'\vert\vert
 J'\vert r(R,R')P(R,R')T_{R'}$$
$$\leq \sum\limits_{h\geq 1}\sum\limits_{R'\subseteq \Omega_h}\sum\limits_{R\subseteq
\Omega}\vert I'\vert\vert
 J'\vert r(R,R') T_{R'}$$

For each $h\geq 1$ and $R'\subseteq \Omega_h$, we decompose $\lbrace
R: R\subseteq\Omega\rbrace$ into
$$A_{0,0}(R')=\big\lbrace R\subseteq\Omega:\,\,\text{dist}(I,I')\leq \vert I\vert \vee \vert I'\vert,\,
\text{dist}(J,J')\leq \vert J\vert \vee \vert J'\vert\big\rbrace; $$
$$A_{i',0}(R')= \big\lbrace R\subseteq\Omega:\,  2^{i'-1}(\vert I\vert \vee \vert I'\vert) < \,\text{dist}(I,I')
\leq 2^{i'}(\vert I\vert \vee \vert I'\vert),\, \text{dist}(J,J')\leq
\vert J\vert \vee \vert J'\vert\big\rbrace;  $$
$$A_{0,\ell'}(R')=\big\lbrace R\subseteq\Omega:\, \text{dist}(I,I')\leq \vert I\vert \vee \vert I'\vert,\
2^{\ell'-1}(\vert J\vert \vee \vert J'\vert)  < \text{dist}(J,J')
\leq 2^{\ell'}(\vert J\vert \vee \vert J'\vert)\big\rbrace; $$
$$A_{i',\ell'}(R')=\big\lbrace R\subseteq\Omega:\ 2^{i'-1}(\vert I\vert \vee \vert I'\vert) < \text{dist}(I,I')
\leq 2^{i'}(\vert I\vert \vee \vert I'\vert),$$
$$ 2^{\ell'-1}(\vert
J\vert \vee \vert J'\vert)  <
\text{dist}(J,J') \leq 2^{\ell'}(\vert
J\vert \vee \vert J'\vert)\big\rbrace, $$ where $i',\ell'\geq 1$.

Now we split $\sum\limits_{h\geq 1}\sum\limits_{R'\subseteq
\Omega_h}\sum\limits_{R\subseteq \Omega}\vert I'\vert\vert
 J'\vert r(R,R')P(R,R')T_{R'}$ into
$$ \sum\limits_{h\geq 1}\sum\limits_{R'\in \Omega_h}\bigg(\sum\limits_{R\in A_{0,0}(R')}+ \sum\limits_{i'\geq
1}\sum\limits_{R\in A_{i',0}(R')}+\sum\limits_{\ell'\geq
1}\sum\limits_{R\in A_{0,\ell'}(R')}
+\sum\limits_{i',\ell'\geq
1}\sum\limits_{R\in A_{i',\ell'}(R')}\bigg)\vert I'\vert\vert
 J'\vert  $$
$$\times r(R,R')P(R,R')T_{R'}=:I_1+I_2+I_3+I_4.$$

To estimate the term $I_1$, we only need to estimate
$\sum\limits_{R\in A_{0,0}(R')}r(R,R')$ since $P(R,R')\leq 1$ in
this case.

Note that $R\in A_{0,0}(R')$ implies
$3R\bigcap3R'\ne \emptyset$. For such $R$, there are four
cases:

Case 1: $\vert I'\vert \geq \vert I\vert $, $\vert J'\vert \leq \vert
J\vert $; Case 2: $\vert I'\vert \leq \vert I\vert $, $\vert J'\vert \geq \vert
J\vert $; Case 3: $\vert I'\vert \geq \vert I\vert $, $\vert J'\vert \geq \vert J\vert $; Case 4: $\vert I'\vert \leq \vert I\vert $, $\vert J'\vert \leq \vert J\vert $.

In each case, we can show $\sum_{R\in A_{0,0}}r(R,R') \leq C2^{-hL}$ by using a simple geometric argument similar to that of Chang-R. Fefferman [CF3].
This, together with (5.1), implies that $I_1$ is bounded by
$$ \sum\limits_{h\geq 1}2^{-hL}\vert \Omega_h\vert^{{{2}\over
{p}}-1}
\leq C \sum\limits_{h\geq 1}h^{{{2}\over{p}}-1}2^{-h(L-{{2}\over{p}}+1)}\vert \Omega^{0,0}\vert^{{{2}\over
{p}}-1}\leq C \vert \Omega\vert^{{{2}\over
{p}}-1},$$
since $\vert \Omega_h\vert\leq Ch2^h\vert \Omega^{0,0}\vert$ and $\vert \Omega^{0,0}\vert
\leq C \vert \Omega\vert.$

Thus it remains to estimate the term $I_4$, since estimates of $I_2$ and $I_3$ can be derived  using the same techniques as in
$I_1$ and $I_4$. This term is more complicated to estimate than term $I_1$.

For each $i',\ell'\geq 1$, when $R\in A_{i',\ell'}(R')$, we have
$P(R,R')\leq 2^{-i'(1+K)}2^{-\ell'(1+K)}$. Similar to estimating term $I_1$, we only need to estimate
the sum $\sum_{R\in A_{i',\ell'}}r(R,R')$. Note that $R\in A_{i',\ell'}(R')$ implies that $3(2^{i'} I\times
2^{\ell'} J)\cap 3(2^{i'} I'\times 2^{\ell'} J')\ne \emptyset$. We also split into four cases in
estimating this sum.

Case 1: $\vert2^{i'} I'\vert \geq \vert2^{i'} I\vert$, $\vert2^{\ell'} J'\vert \leq \vert2^{\ell'} J\vert$. Then
$$ {{\vert 2^{i'}I\vert }\over{\vert 3\cdot 2^{i'}I'\vert }}\vert 3(2^{i'}I'\times2^{\ell'}J' )\vert
\leq \vert 3(2^{i'}I'\times2^{\ell'}J' )\cap 3(2^{i'}I\times2^{\ell'}J )\vert$$
$$\leq C 2^{i'}2^{\ell'}\vert 3R'\cap \Omega^{0,0}\vert \leq C 2^{i'}2^{\ell'}{{1}\over{2^{h-1}}}\vert
3R'\vert \leq C{{1}\over{2^{h-1}}}|3(2^{i'}I'\times2^{\ell'}J' )|.$$
Thus $\vert 2^{i'}I'\vert =2^{h-1+n}\vert 2^{i'}I\vert $ for some $n\geq 0$. For each fixed
$n$, the number of such $2^{i'}I$'s must be $\leq 2^n\cdot 5$. As for
$\vert 2^{\ell'}J\vert=2^m\vert2^{\ell'} J'\vert$ for some $m\geq 0$, for each fixed $m$,
$3\cdot2^{\ell'}J\cap 3\cdot2^{\ell'}J' \ne \emptyset$ implies that the number of such $2^{\ell'}J'$ is less than 5.
Thus
$$ \sum_{R\in case1}r(R,R') \leq \sum_{m,n\geq 0} \bigg({1\over{2^{n+m+h-1}}} \bigg)^L 2^n \cdot 5^2 \leq C2^{-hL}.$$

Similarly, we can handle

Case 2: $\vert2^{i'} I'\vert \leq \vert2^{i'} I\vert$, $\vert2^{\ell'} J'\vert \geq \vert2^{\ell'} J\vert$,

Case 3: $\vert2^{i'} I'\vert \geq \vert2^{i'} I\vert$, $\vert2^{\ell'} J'\vert \geq \vert2^{\ell'} J\vert$ and

Case 4: $\vert2^{i'} I'\vert \leq \vert2^{i'} I\vert$, $\vert2^{\ell'} J'\vert \leq \vert2^{\ell'} J\vert$.

Combining the four cases, we have $ \sum_{R\in A_{i',\ell'}(R')}r(R,R') \leq C2^{-hL},$ which, together with the
estimate of $P(R,R'),$ implies that
$$I_4\leq C\sum\limits_{h\geq 1}\sum\limits_{i',\ell'\geq
          1}\sum\limits_{R'\subseteq \Omega_h}2^{-hL}2^{-i'(1+K)}2^{-\ell'(1+K)}\vert I'\vert\vert
           J'\vert T_{R'}. $$
Hence $I_4$ is bounded by
$$\sum\limits_{h\geq 1}2^{-hL}\vert \Omega_h\vert^{{{2}\over{p}}-1}
\leq C \sum\limits_{h\geq 1}h^{{{2}\over{p}}-1}2^{-h(L-{{2}\over{p}}+1)}\vert \Omega^{0,0}\vert^{{{2}\over
{p}}-1}\leq C \vert \Omega\vert^{{{2}\over {p}}-1},$$
since $\vert \Omega_h\vert\leq Ch2^h\vert \Omega^{0,0}\vert$ and $\vert \Omega^{0,0}\vert
\leq C \vert \Omega\vert.$

Combining $I_1$, $I_2$, $I_3$ and $I_4$, we have
$${1\over{\vert \Omega\vert^{{{2}\over {p}}-1}}}\sum_{R'\in {B_{0,0}}}\sum\limits_{R\subseteq
\Omega}\vert I'\vert\vert
 J'\vert r(R,R')P(R,R')T_{R'}\leq C \sup_{\bar{\Omega}}{1\over{\vert \bar{\Omega}\vert^{{{2}\over {p}}-1}}}
 \sum_{R'\subseteq\bar{\Omega}}\vert I'\vert \vert J'\vert T_{R'} . $$

Now we consider
 $$\sum\limits_{i,\ell\geq 1}\sum\limits_{R'\in
{B_{i,\ell}}}\sum\limits_{R\subseteq \Omega} \vert I'\vert\vert
J'\vert r(R,R')P(R,R')T_{R'}.$$
 Note that for $R'\in {B_{i,\ell}}$,
$3(2^iI'\times2^{\ell}J') \cap\Omega^{i,l}\ne \emptyset$. Let
$${\Cal {F}}_h^{i,\ell}=\big\lbrace R'\in B_{i,\ell}: \vert
3(2^iI'\times 2^{\ell}J') \cap\Omega^{i,l}\vert \geq {1\over
{2^h}}\vert3(2^iI'\times2^{\ell}J')\vert \big\rbrace,$$
$${\Cal {D}}_h^{i,\ell}={\Cal {F}}_h^{i,\ell}\setminus
{\Cal {F}}_{h-1}^{i,\ell}$$ and $$\Omega_h^{i,\ell}=\bigcup_{R'\in
{\Cal {D}}_h^{i,\ell}}R'.$$. Since $B_{i,\ell}=\bigcup_{h\geq 1}{\Cal {D}}_h^{i,\ell},$ we first estimate
$$\sum_{R'\in {\Cal {D}}_h^{i,\ell}}\sum_{R\subseteq \Omega}\vert I'\vert \vert J'\vert r(R,R')P(R,R')T_{R'}$$ for some
$i,\ell,h\geq 1$.

Note that for each $R'\in {\Cal {D}}_h^{i,\ell}$, $3(2^iI'\times2^{\ell}J') \cap\Omega^{i-1,l-1}=\emptyset$. So for any
$R\subseteq\Omega$, we have $2^i(\vert I\vert \vee \vert I'\vert)\leq \,\text{dist}(I,I')$ and $2^{\ell}
(\vert J\vert \vee \vert J'\vert)\leq \,\,\text{dist}(J,J')$.
We decompose $\lbrace R:\ R\subseteq\Omega \rbrace$ by
$$ A_{i',\ell'}(R')=\lbrace R\subset\Omega:\ 2^{i'-1} \cdot2^i(\vert I\vert \vee \vert I'\vert)
\leq  \,\text{dist}(I,I')\leq  2^{i'} \cdot2^i(\vert I\vert \vee \vert I'\vert),$$
$$2^{\ell'-1} \cdot 2^{\ell}(\vert J\vert \vee \vert J'\vert)
\leq  \,\text{dist}(J,J')\leq  2^{\ell'} \cdot 2^{\ell}(\vert J\vert \vee \vert J'\vert)\rbrace, $$
where $i',\ell'\geq 1$. Then we write
$$\sum_{R'\in {\Cal {D}}_h^{i,\ell}}\sum_{R\subseteq \Omega}\vert I'\vert \vert J'\vert r(R,R')P(R,R')T_{R'}=
\sum_{i',\ell'\geq 1}\sum_{R'\in {\Cal {D}}_h^{i,\ell}}\sum_{R\in A_{i',\ell'}(R')}\vert I'\vert \vert J'\vert
r(R,R')P(R,R')T_{R'}$$

Since $P(R,R')\leq 2^{-i(1+K)}2^{-\ell (1+K)}2^{-i'(1+K)}2^{-\ell' (1+K)}$ for $R'\in B_{i,\ell}$
and $R\in A_{i',\ell'}(R'),$ repeating the same proof
with $B_{0,0}$ replaced by $B_{i,\ell}$ and a necessary modification yields
$$\sum_{R'\in {\Cal {D}}_h^{i,\ell}}\sum_{R\in A_{i',\ell'}(R')}\vert I'\vert \vert J'\vert
r(R,R')P(R,R')T_{R'}\leq C 2^{-i(1+K)}2^{-\ell (1+K)}2^{-i'(1+K)}2^{-\ell' (1+K)}\times $$
$$ i^{ {2\over p}-1  } 2^{i({2\over p}-1)  }{\ell}^{ {2\over p}-1  } 2^{\ell({2\over p}-1)  }
h^{ {2\over p}-1  } 2^{-h(L-{2\over p}+1)  }
\sup_{\bar{\Omega}}{1\over{\vert \bar{\Omega}\vert^{{{2}\over {p}}-1}}}
\sum_{R'\subseteq\bar{\Omega}}\vert I'\vert \vert J'\vert T_{R'}
.$$
Adding over all $i,\ell,i',\ell',h \geq 1 ,$ we have
$${1\over{\vert \Omega \vert^{{{2}\over {p}}-1}}} \sum\limits_{i,\ell\geq 1}\sum\limits_{R'\in
{B_{i,\ell}}}\sum\limits_{R\subseteq \Omega} \vert I'\vert\vert
J'\vert r(R,R')P(R,R')T_{R'}\leq C\sup_{\bar{\Omega}}{1\over{\vert \bar{\Omega}\vert^{{{2}\over {p}}-1}}}
\sum_{R'\subseteq\bar{\Omega}}\vert I'\vert \vert J'\vert T_{R'}.$$

Similar estimates hold for
$$\sum\limits_{i\geq 1}\sum\limits_{R'\in
{B_{i,0}}}\sum\limits_{R\subseteq \Omega} \vert I'\vert\vert J'\vert r(R,R')P(R,R')T_{R'}$$
and
$$\sum\limits_{\ell\geq 1}\sum\limits_{R'\in
{B_{0,\ell}}}\sum\limits_{R\subseteq \Omega} \vert I'\vert\vert J'\vert r(R,R')P(R,R')T_{R'},$$
which, after adding over all $i,\ell\geq 0$, shows Theorem 1.21. We leave the details to the reader.
\hfill {\bf Q.E.D.}

As a consequence of Theorem 1.14, it is easy to see that the space $CMO^p_F$ is well defined.
Particularly, we have

\proclaim{Corollary 5.1} We have
$$ \Vert f \Vert_{CMO^p_F}
\approx
\sup_{\Omega}\left\lbrace{{1}\over{\vert \Omega \vert}^{{{2}\over{p}}-1}}\sum \limits_j\sum
 \limits_k\sum
\limits_{I\times J\subseteq \Omega}
\vert \psi _{j,k}*f(x_I,y_J)
\vert ^2\vert I\vert \vert J\vert \right\rbrace^{{1}\over {2}},$$
where $I, J$ are dyadic cubes in $R^n, R^m$ with length $2^{-j-N}, 2^{-j-N}+2^{-k-N},$
and $x_I, y_J$ are any fixed points in $I, J,$ respectively.
\endproclaim

Before we prove Theorem 1.16, we remark again that Theorem 1.16 with $p=1$ in the one-parameter setting was proved in [FJ] on
$R^n$ by use of
the distribution inequalities. This method is difficult to apply to multi-parameter case.
We will give a simpler and more constructive proof which also gives a new proof of the result in [FJ].
Moreover, this constructive
proof works also for other multi-parameter cases.

\noindent {\bf Proof of Theorem 1.16:} We first prove $c^p\subseteq (s^p)^*.$ Applying the
proof in Theorem 4.4, set
$$ s(x,y)=\lbrace \sum \limits_
{I\times J}\vert  s_{I\times J}\vert^2\vert I\vert^{-1}\vert J\vert^{-1}
\chi _I(x)\chi _J(y)\rbrace^{{1}\over{2}}$$
and
$$\Omega_i=\lbrace (x,y)\in R^n\times R^m: s(x,y)
>2^i\rbrace.$$
Denote $${\Cal  B}_i=\lbrace (I\times J): \vert (I\times
J)\cap\Omega_i\vert>{{1}\over
 {2}}\vert I\times J\vert, \vert (I\times J)\cap\Omega_{i+1}\vert\leq{{1}\over
 {2}}\vert I\times J\vert\rbrace,$$ where $I, J$ are dyadic cubes in $R^n, R^m,$ with side length
$2^{-j-N},$ and $2^{-j-N}+2^{-k-N},$ respectively.
Suppose $t=\lbrace t_{I\times J}\rbrace\in c^p$ and write
$$
\align &
 \vert \sum\limits_{I\times J} s_{I\times J} \overline {t}_{I\times J}\vert \vert \\ = &
\vert \sum\limits_i\sum\limits_{(I\times J)\in {\Cal  B}_i}s_{I\times J} \overline {t}_{I\times
 J}\vert\\ \leq  & \lbrace \sum\limits_i\lbrace \sum\limits_{(I\times J)\in {\Cal  B}_i}
\vert s_{I\times J}\vert^2\rbrace^{{p}\over{2}} \lbrace \sum\limits_{(I\times J)\in {\Cal  B}_i}
\vert t_{I\times J}\vert^2\rbrace^{{p}\over{2}}\rbrace ^{{1}\over{p}}\\
\leq &  C\Vert t \Vert_{c^p} \lbrace \sum\limits_i\vert \Omega_i\vert^{1-{{p}\over{2}}}\lbrace
 \sum\limits_{(I\times J)\in {\Cal   B}_i}\vert s_{I\times J}\vert^2\rbrace^{{p}\over{2}}\rbrace^{{1}\over {p}} \tag5.2
 \endalign
 $$
since if $I\times J\in{{\Cal  B}_i}$, then $$I\times J\subseteq {\widetilde{\Omega_i}}=
\lbrace (x,y): M_s(\chi_{\Omega_i})(x,y)>{{1}\over {2}}\rbrace, \vert {\widetilde{\Omega_i}}
\vert \leq C\vert {\Omega_i}\vert,$$
and $\lbrace t_{I\times J}\rbrace\in c^p$ yields $$\lbrace \sum\limits_{(I\times J)\in {\Cal  B}_i}
\vert t_{I\times J}\vert^2\rbrace^{{1}\over{2}}\leq C \Vert t
\Vert_{c^p}
 \vert {\Omega_i}\vert^{{{1}\over{p}}-{{1}\over{2}}}.$$

The same proof as in the claim of Theorem 4.4 implies
$$\sum\limits_{(I\times J)\in {\Cal   B}_i}
\vert s_{I\times J}\vert^2 \leq C2^{2i}\vert{\Omega_i}\vert.$$
Substituting the above term back to the last term in (5.2) gives $c^p\subseteq (s^p)^*.$

The proof of the converse is simple and is similar to one given in [FJ] for $p=1$ in the one-parameter setting on $R^n$. If
 $\ell\in
(s^p)^*,$ then it is clear that $\ell(s)=
\sum\limits_{I\times J}s_{I\times J}{\overline t}_{I\times J}$
 for some $t=\lbrace t_{I\times J}\rbrace$. Now fix an open set $\Omega\subset R^n\times R^m$
and let $S$ be the sequence space of all $s=\lbrace s_{I\times J}\rbrace$ such that $I\times
 J\subseteq
\Omega.$ Finally, let $\mu$ be a measure on $S$ so that the $\mu-$measure of the ``point'' $I\times J$ is
${{1}\over{{\vert \Omega\vert}^{{{2}\over
{p}}-1}}}.$ Then,
$$\lbrace {{1}\over {\vert \Omega\vert^{{{2}\over
{p}}-1}}}
\sum_{I\times J\subseteq \Omega}
\vert t_{I\times J}\vert ^2\rbrace^{{1}\over {2}}=\Vert t_{I\times J}
\Vert_{\ell^2(S,d\mu)}$$
$$=\sup_{\Vert s\Vert_{\ell^2(S,d\mu)}\leq 1}\vert {{1}\over {\vert \Omega\vert^{{{2}\over
{p}}-1}}}\sum\limits_{I\times J\subseteq\Omega}s_{I\times J}{\overline t}_{I\times J}\vert$$
$$\leq \Vert t \Vert_{(s^p)^*}\sup_{\Vert s\Vert_{\ell^2(S,d\mu)}\leq 1}
\Vert s_{I\times J}{{1}\over {\vert
 \Omega\vert^{{{2}\over
{p}}-1}}}\Vert_{s^p}.$$
By Holder's inequality,
$$
\align
& \Vert s_{I\times J}{{1}\over {\vert
 \Omega\vert^{{{2}\over
{p}}-1}}}\Vert_{s^p}\\ = &
{{1}\over {\vert \Omega\vert^{{{2}\over
{p}}-1}}}\lbrace \int\limits_{\Omega}(\sum\limits_{I\times J\subseteq\Omega}
\vert s_{I\times J}\vert^2\vert I\times J
\vert^{-1}\chi_I(x)\chi_J(y))^{{p}\over{2}}dxdy
\rbrace^{{1}\over{p}}\\
\leq  & \lbrace {{1}\over {\vert \Omega\vert^{{{2}\over
{p}}-1}}}\int\limits_{\Omega}\sum\limits_{I\times J\subseteq\Omega}
\vert s_{I\times J}\vert^2{\vert I\times J\vert}^{-1}
\chi_I(x)\chi_J(y)dxdy
\rbrace^{{1}\over{2}}\\ = & \Vert s\Vert_{\ell^2(S,d\mu)}\leq 1,
\endalign
$$
which shows $\Vert t\Vert_{c^p}\leq \Vert t\Vert_{(s^p)^*}.$ \hfill{\bf Q.E.D.}

In order to use Theorem 1.16 to show Theorem 1.17, we define a map $S$ which takes $f\in
({\Cal  S}_F)^\prime$ to the sequence of coefficients $\lbrace s_{I\times J}\rbrace=\lbrace\vert
 I\vert^{{1}\over{2}}\vert J\vert^{{1}\over{2}}\psi_{j,k}*f(x_I,y_J)\rbrace,$
where $I, J$ are cubes in $R^n, R^m,$ with side length $2^{-j-N}, 2^{-j-N}+2^{-k-N},$ and
$x_I, y_J$ are any fixed points in $I, J,$ respectively.
For any sequence $s=\lbrace s_{I\times J}\rbrace,$ we define the map $T$ which takes $s$ to
$T(s)=\sum \limits_j\sum \limits_k\sum \limits_J\sum
\limits_I \vert I\vert^{{1}\over{2}} \vert J\vert^{{1}\over{2}} {\widetilde \psi }_{j,k}(x,y)
s_{I\times J},$ where ${\widetilde \psi }_{j,k}$ are the same as in (1.9).

The following result together with Theorem 1.16 will show theorem 1.17.

\proclaim{Theorem 5.2:} The maps $S: H_F^p\rightarrow s^p$ and $ CMO^p_F\rightarrow
 c^p,$ and $T: s^p\rightarrow H_F^p$ and $c^p\rightarrow CMO^p_F$ are bounded, and $T\circ S$ is the identity on $H_F^p$ and $CMO^p_F.$
\endproclaim
\noindent {\bf Proof of Theorem 5.2:} The boundedness of $S$ on $H_F^p$ and $CMO^p_F$ follows
directly from the Plancherel-P\^olya inequalities, Theorem 1.9 and Theorem 1.14. The boundedness of
$T$ follows from the same proofs in Theorem 1.9 and 1.14. To be precise, to see $T$ is bounded
from $s^p$ to
$H_F^p,$ let $s=\lbrace s_{I\times J}\rbrace.$ Then, by Proposition 4.1,
$$ \Vert T(s)\Vert_{H_F^p}\leq C\Vert \lbrace \sum \limits_j\sum
 \limits_k\sum
\limits_J\sum \limits_I \vert \psi _{j,k}*T(s)(x,y)
\vert ^2\chi _I(x)\chi _J(y)\rbrace ^{{{1}\over{2}}}\Vert _p.$$
By adapting a similar the proof in Theorem 1.9, we have for some $0<r<p$
$$
\align &
\vert \psi_{j,k}*T(s)(x,y)\chi_I(x)\chi_J(y)\vert^2\\
= &
\vert \sum \limits_{j^\prime }\sum \limits_{k^\prime
}\sum \limits_{J^\prime }\sum \limits_{I^\prime }\vert I^\prime\vert\vert I^\prime\vert
\psi _{j,k}*{\widetilde \psi }_{{j^\prime },{k^\prime
}}(\cdot,\cdot)(x,y) s_{I^\prime\times J^\prime}\vert I^\prime\vert^{-{{1}\over{}2}}
\vert J^\prime\vert^{-{{1}\over{}2}}\chi_I(x)\chi_J(y)\vert^2\\
\leq  & C \sum \limits_{ k^\prime\leq j^\prime
}2^{-\vert j-{j^\prime }\vert K }2^{-\vert k-{k^\prime }\vert K}\lbrace M_s(
\sum \limits_{J^\prime }\sum \limits_{I^\prime }\vert
s_{I^\prime\times J^\prime}\vert\vert I^\prime\vert^{-1}
\vert J^\prime\vert^{-1} \chi _{J^{
\prime {\chi _{I^\prime }}}})^r\rbrace^{{{2}\over{r}}}(x,y)
\chi _I(x)\chi _J(y)\\ + &
\sum \limits_{ k^\prime > j^\prime
}2^{-\vert j-{j^\prime }\vert K }2^{-\vert k-{k^\prime }\vert K}
\lbrace M(
\sum \limits_{J^\prime }\sum \limits_{I^\prime }\vert
s_{I^\prime\times J^\prime}\vert\vert I^\prime\vert^{-1}
\vert J^\prime\vert^{-1} \chi _{J^{
\prime {\chi _{I^\prime }}}})^r\rbrace^{{{2}\over{r}}}(x,y)
\chi _I(x)\chi _J(y).
\endalign
$$
Repeating the proof in Theorem 1.9 gives the boundedness of $T$ from $s^p$ to $H_F^p.$ The
similar proof given in the proof of Theorem 1.14 applies to the boundedness of $T$ from $c^p$
to $CMO_F^p.$ We leave the
details to the reader. The discrete Calder\'on reproducing formula, Theorem 1.8, and Theorem 1.14
 show that $T\circ S$ is the identity on $H_F^p$ and $CMO^p_F.$

We are now ready to give the

\noindent {\bf Proof of Theorem 1.17:} If $f\in {\Cal  S}_F$ and $g\in CMO_F^p$, then the discrete
Calder\'on reproducing formula, Theorem 1.16 and Theorem 5.2 imply
$$\vert \ell_g\vert =\vert <f, g>\vert =$$
$$\vert \sum\limits_{R=I\times J}\vert I\vert \vert J\vert \psi_R*f(x_I,y_J)
{\widetilde \psi}_R(g)(x_I,y_J)\vert\leq C \Vert f\Vert_{H_F^p}\Vert g\Vert_{CMO^p_F}.$$

Because ${\Cal  S}_F$ is dense in $H_F^p,$ this shows that the map $\ell_g=<f, g>,$ defined initially for $f\in{\Cal  S}_F$ can be
extended to a continuous linear functional on $H_F^p$ and $\Vert \ell_g\Vert\leq C\Vert
 g\Vert_{CMO^p_F}$.

Conversely, let $\ell\in (H_F^p)^*$ and set $\ell_1=\ell \circ T,$ where $T$ is defined as in
Theorem 5.2. Then, by theorem 5.2, $\ell_1\in (s^p)^*,$ so by Theorem 1.16, there exists $t=\lbrace
t_{I\times J}\rbrace$ such that $\ell_1(s)= \sum\limits_{I\times J}s_{I\times J}{\overline t}_{
I\times J}$ for all $s=\lbrace s_{I\times J}\rbrace$ and $\Vert t\Vert_{c^p}\approx \Vert
 \ell_1\Vert\leq C\Vert \ell\Vert$ because $T$ is bounded. Again, by Theorem 1.16, $\ell =\ell\circ
 T\circ S=\ell_1\circ S.$ Hence, with $f\in {\Cal  S}_F$ and $g=\sum\limits_{I\times J}
t_{I\times J}\psi_R(x-x_I, y-y_J),$ where, without loss the generality, we may assume $\psi$ is a radial function,
$$\ell(f)=\ell_1(S(f))=<S(f), t>=<f, g>,$$
This proves $\ell =\ell_g$ and, by Theorem 1.16, $\Vert g\Vert _{CMO_F^p}\leq C \Vert t \Vert
_{c^p}\leq C \Vert \ell_g\Vert.$
 \hfill {\bf Q.E.D.}

{\bf Proof of Theorem 1.18}
As mentioned in section 4, since $H_F^1$ is a subspace of $L^1,$ by the duality of $H_F^1$ and
$BMO_F$, and the boundedness of flag singular integrals on $H_F^1$, one concludes that $L^{\infty}$ is a subspace of $BMO_F,$ and flag singular integrals are bounded
on $BMO_F$ and from $L^{\infty}$ to
$BMO_F.$ This shows Theorem 1.18. \hfill {\bf Q.E.D.}

\vskip 0.2cm

\noindent {\bf 6. Calder\'on-Zygmund decomposition and interpolation on flag Hardy spaces $\phpf$: Proofs of Theorems 1.19 and 1.20 }
\vskip .2cm

The main purpose of this section is to derive a Calder\'on-Zygmund decomposition using functions in flag Hardy spaces. Furthermore, we will prove an interpolation theorem on $\phpf$.

We first recall that Chang and R. Fefferman established the following Calder\'on-Zygmund decomposition on the pure  product domains $R^2_{+}\times R^2_{+}$ ([CF2]).

{\bf Calder\'on-Zygmund Lemma:} Let $\alpha>0$ be given and $f\in L^p(R^2)$, $1<p<2$. Then we may write $f=g+b$ where $g\in L^2(R^2)$ and $b\in H^1(R^2_{+}\times R^2_{+})$ with $||g||_2^2\le \alpha^{2-p}||f||_p^p$ and $||b||_{H^1(R^2_{+}\times R^2_{+})}\le C\alpha^{1-p}||f||_p^p$, where $c$ is an absolute constant.

We now prove the Calder\'on-Zygmund decomposition in the setting of flag Hardy spaces, namely we give
the

{\bf Proof of Theorem 1.19}
 We first assume $f\in L^2(R^{n+m})\cap \phpf.$ Let $\alpha>0$ and $\Omega_\ell=\{(x,y)\in \prods:
S(f)(x,y)>\alpha 2^\ell\}$, where, as in Corollary 4.6,
$$S(f)(x,y)=\left\{\sum\limits_{j,k}\sum\limits_{I,J}|\phi_{jk}*\left (T_N^{-1}(f)\right )(x_I,y_J)|^2\chi_I(x)\chi_J(y)\right\}^{\frac{1}{2}}.$$
It has been shown in Corollary 4.6 that $f\in L^2(R^{n+m})\cap \phpf$ then
$||f||_{H^p_F}\approx ||S(f)||_{p}$.

In the following we take $R=I\times J$ as all dyadic rectangles in $\prods$ with $|I|=2^{-j-N}$, $|J|=2^{-j-N}+2^{-k-N}$, where $j, k$ are integers and $N$ is large enough.

Let $${\Cal R}_0=\left\{R=I\times J,  \,\,\text{such that}\,\, |R\cap \Omega_0|<\frac{1}{2}|R|\right\}$$
and for $\ell\ge 1$
$${\Cal R}_\ell=\left\{R=I\times J, \,\,\text{such that}\,\, |R\cap \Omega_{\ell-1}|\ge \frac{1}{2}|R|\,\,\text{but}\,\, |R\cap \Omega_\ell|<\frac{1}{2}|R|\right\}.$$
By the discrete Calder\'on reproducing formula in Theorem 4.4,
$$
\align
f(x,y)= & \sum\limits_{j,k}\sum\limits_{I,J}|I||J|{\widetilde \phi}_{jk}(x-x_I,y-y_J)\phi_{jk}*\left (T_N^{-1}(f)\right )(x_I,y_J)\\
= & \sum\limits_{\ell\ge 1}\sum\limits_{I\times J\in {\Cal R}_\ell}|I||J|{\widetilde \phi}_{jk}(x-x_I,y-y_J)\phi_{jk}*\left (T_N^{-1}(f)\right )
(x_I,y_J)\\+ & \sum\limits_{I\times J\in {\Cal R}_0}|I||J|{\widetilde \phi}_{jk}(x-x_I,y-y_J)\phi_{jk}*\left (T_N^{-1}(f)\right )(x_I,y_J)\\= & b(x,y)+g(x,y)
\endalign
$$
When $p_1>1$, using duality argument, it is easy to show
$$||g||_{p_1}\le C ||\left\{ \sum\limits_{R=I\times J\in {\Cal R}_0}|\phi_{jk}*\left (T_N^{-1}(f)\right )(x_I,y_J)|^2\chi_I\chi_J\right\}^{\frac{1}{2}}||_{p_1}.
$$
Next, we estimate $||g||_{H^{p_1}_F}$ when $0<p_1\le 1$. Clearly, the duality argument will not work here. Nevertheless,
we can estimate the $H^{p_1}_F$  norm directly by using the discrete Calder\'on reproducing formula in Theorem 1.8. To this end, we note that
$$
||g||_{H^{p_1}_F}\le ||\left\{\sum\limits_{j^\prime, k^\prime}\sum\limits_{I^\prime, J^\prime}|\left(\psi_{j^\prime k^\prime}*g\right)(x_{I^\prime}, y_{J^\prime})|^2\chi_{I^\prime}(x)\chi_{J^\prime}(y)\right\}^{\frac{1}{2}}||_{L^{p_1}}.$$
\vskip -0.4cm
Since
$$\left(\psi_{j^\prime,k^\prime}*g\right)(x_{I^\prime}, y_{J^\prime})=\sum\limits_{I\times J\in {\Cal R}_0}|I||J|
\left(\psi_{j^\prime k^\prime}*{\widetilde \phi}_{jk}\right)(x_{I^\prime}-x_I,y_{J^\prime}-y_J)\phi_{jk}*\left (T_N^{-1}(f)\right )(x_I,y_J)$$
Repeating the same proof of Theorem 1.9, we have
$$
\align
&
||\left\{\sum\limits_{j^\prime, k^\prime}\sum\limits_{I^\prime, J^\prime}|\left(\psi_{j^\prime k^\prime}*g\right)
(x_{I^\prime}, y_{J^\prime})|^2\chi_{I^\prime}(x)\chi_{J^\prime}(y)\right\}^{\frac{1}{2}}||_{L^{p_1}}\\
 \le & C ||\left\{ \sum\limits_{R=I\times J\in {\Cal R}_0}|\phi_{jk}*\left (T_N^{-1}(f)\right )(x_I,y_J)|^2\chi_I\chi_J\right\}^{\frac{1}{2}}||_{p_1}.
\endalign
$$
This shows that for all $0<p_1<\infty$
$$||g||_{H^{p_1}_F}\le  C ||\left\{ \sum\limits_{R=I\times J\in {\Cal R}_0}|\phi_{jk}*\left (T_N^{-1}(f)\right )(x_I,y_J)|^2\chi_I\chi_J\right\}^{\frac{1}{2}}||_{p_1}.
$$
\vskip -0.4cm
{\bf Claim 1:}
$$\int_{S(f)(x,y)\le \alpha}S^{p_1}(f)(x,y)dxdy\ge C
||\left\{ \sum\limits_{R=I\times J\in {\Cal R}_0}|\phi_{jk}*\left (T_N^{-1}(f)\right )(x_I,y_J)|^2\chi_I\chi_J\right\}^{\frac{1}{2}}||_{p_1}.
$$
This claim implies
$$
\align
& ||g||_{p_1}\le C \int_{S(f)(x,y)\le \alpha}S^{p_1}(f)(x,y)dxdy\\ \le & C\alpha^{p_1-p}  \int_{S(f)(x,y)\le \alpha}S^{p}(f)(x,y)dxdy\\ \le & C\alpha^{p_1-p}||f||^p_{\phpf}.
\endalign
$$
To show Claim 1, we denote $R=I\times J\in {\Cal R}_0$. We choose $0<q<p_1$ and note that
$$
\align
& \int_{S(f)(x,y)\le \alpha}S^{p_1}(f)(x,y)dxdy\\
= & \int_{S(f)(x,y)\le \alpha}\left\{\sum\limits_{j,k}\sum\limits_{I,J}|\phi_{jk}*\left (T_N^{-1}(f)\right )(x_I,y_J)|^2\chi_I(x)\chi_J(y)\right\}^{\frac{p_1}{2}}dxdy\\
\ge & C\int_{\Omega_0^c} \left\{ \sum\limits_{R \in {\Cal R}_0}|\phi_{jk}*\left (T_N^{-1}(f)\right )(x_I,y_J)|^2\chi_I\chi_J\right\}^{\frac{p_1}{2}}dxdy\\
= & C\int_{\prods} \left\{ \sum\limits_{R \in {\Cal  R}_0}|\phi_{jk}*\left (T_N^{-1}(f)\right )(x_I,y_J)|^2\chi_{R\cap \Omega_0^c}(x,y)\right\}^{\frac{p_1}{2}}dxdy
\\ \ge & C \int_{\prods} \left\{\left\{ \sum\limits_{R \in {\Cal  R}_0}\left(M_s\left(|\phi_{jk}*\left (T_N^{-1}(f)\right )(x_I,y_J)|^q\chi_{R\cap \Omega_0^c}\right)(x,y)\right)^{\frac{2}{q}}\right\}^{\frac{q}{2}}\right\}^{\frac{p_1}{q}}dxdy\\
\ge & C \int_{\prods} \left\{ \sum\limits_{R \in {\Cal  R}_0}|\phi_{jk}*\left (T_N^{-1}(f)\right )(x_I,y_J)|^2\chi_R (x,y)\right\}^{\frac{p_1}{2}}dxdy
\endalign
$$
where in the last inequality we have used the fact that
$|\Omega_0^c\cap (I\times J)|\ge \frac{1}{2}|I\times J|$ for $I\times J\in {\Cal R}_0$, and thus
$$\chi_R(x,y)\le 2^{\frac{1}{q}} M_s(\chi_{R\cap\Omega_0^c})^{\frac{1}{q}}(x,y)$$
and in the second to the last inequality we have used the vector-valued Fefferman-Stein inequality for strong maximal functions
$$||\left(\sum\limits_{k=1}^\infty(M_s(f_k))^r\right)^{\frac{1}{r}}||_p\le C ||\left(\sum\limits_{k=1}^\infty|f_k|^r\right)^{\frac{1}{r}}||_p$$
with the exponents $r=2/q>1$ and $p=p_1/q>1$.
Thus the claim follows.

We now recall $\widetilde{\Omega_\ell}=\{(x,y)\in \prods: M_s(\chi_{\Omega_\ell})>\frac{1}{2}\}$.

{\bf Claim 2:} For $p_2\le 1$,
$$
||\sum\limits_{I\times J\in {\Cal R}_\ell}|I||J|{\widetilde \phi}_{jk}(x-x_I,y-y_J)\phi_{jk}*\left (T_N^{-1}(f)\right )(x_I,y_J)||^{p_2}_{H^{p_2}_F}\le C(2^{\ell}\alpha)^{p_2}|\widetilde{\Omega_{\ell-1}}|.
$$
Claim 2 implies
$$
\align
||b||^{p_2}_{H^{p_2}_F}\le & \sum\limits_{\ell\ge 1}(2^{\ell}\alpha)^{p_2}|\widetilde{\Omega_{\ell-1}}|\\ \le  & C \sum\limits_{\ell\ge 1}(2^{\ell}\alpha)^{p_2}|\Omega_{\ell-1}|
\le  C\int_{S(f)(x,y)>\alpha}S^{p_2}(f)(x,y)dxdy \\ \le & C\alpha^{p_2-p}\int_{S(f)(x,y)>\alpha}S^{p}(f)(x,y)dxdy\le C\alpha^{p_2-p}||f||^p_{H^p_F}.
\endalign
$$
\vskip -0.4cm
To show Claim 2, again we
have
$$
\align
&
||\sum\limits_{I\times J\in {\Cal R}_\ell}|I||J|{\widetilde \phi}_{jk}(x-x_I,y-y_J)\phi_{jk}*\left (T_N^{-1}(f)\right )(x_I,y_J)||^{p_2}_{H^{p_2}_F}\\
 \le & C || \left\{\sum\limits_{j^\prime k^\prime}\sum\limits_{I^\prime,J^\prime}\left\vert\sum\limits_{I\times J\in
{\Cal R}_\ell}|I||J|\left(\psi_{j^\prime k^\prime}*{\widetilde \phi}_{jk}\right)(x_{I^\prime}-x_I,y_{J^\prime}-y_J)\phi_{jk}*\left (T_N^{-1}(f)\right )(x_I,y_J)\right\vert^2\right\}^{\frac{1}{2}}||_{L^{p_2}}\\
  \le &  C ||\left\{ \sum\limits_{R=I\times J\in {\Cal R}_\ell}|\phi_{jk}*\left (T_N^{-1}(f)\right )(x_I,y_J)|^2\chi_I\chi_J\right\}^{\frac{1}{2}}||_{p_2}
\endalign
$$
where we can use a similar argument in the proof of Theorem 1.19 to prove the last inequality.

However,
$$
\align &
\sum\limits_{\ell=1}^\infty (2^\ell \alpha)^{p_2}|\widetilde{\Omega}_{\ell-1}|\\ \ge &
\int_{\widetilde{\Omega}_{\ell-1}\backslash \Omega_{\ell}}S(f)^{p_2}(x,y)dxdy\\ = & \
\int_{\widetilde{\Omega}_{\ell-1}\backslash \Omega_{\ell}}\left\{ \sum\limits_{j,k}\sum\limits_{I,J}|\phi_{jk}*\left (T_N^{-1}(f)\right )(x_I,y_J)|^2\chi_I(x)\chi_J(y)\right\}^{\frac{p_2}{2}}dxdy\\
= & \int_{\prods} \left\{ \sum\limits_{j,k}\sum\limits_{I,J}|\phi_{jk}*\left (T_N^{-1}(f)\right )(x_I,y_J)|^2\chi_{(I\times J)\cap \widetilde{\Omega}_{\ell-1}\backslash \Omega_{\ell})}(x,y) \right\}^{\frac{p_2}{2}}dxdy\\
\ge & \int_{\prods} \left\{ \sum\limits_{I\times J\in {\Cal R}_\ell}|\phi_{jk}*\left (T_N^{-1}(f)\right )(x_I,y_J)|^2\chi_{(I\times J)\cap \widetilde{\Omega}_{\ell-1}\backslash \Omega_{\ell})}(x,y) \right\}^{\frac{p_2}{2}}dxdy\\
\ge &\int_{\prods} \left\{ \sum\limits_{I\times J\in {\Cal R}_\ell}|\phi_{jk}*\left (T_N^{-1}(f)\right )(x_I,y_J)|^2\chi_{I}(x)\chi_J(y) \right\}^{\frac{p_2}{2}}dxdy
\endalign
$$

In the above string of inequalities, we have used the fact that for $R\in {\Cal R}_\ell$ we have
$$|R\cap \Omega_{\ell-1}|> \frac{1}{2}|R|\,\,\,\text{and}\,\,|R\cap \Omega_\ell|\le \frac{1}{2}|R|$$
and consequently $R\subset \widetilde{\Omega}_{\ell-1}$. Therefore $|R\cap (\widetilde{\Omega}_{\ell-1}\backslash \Omega_\ell)|>\frac{1}{2}|R|$.
Thus the same argument applies here to conclude the last inequality above. Finally, since $L^2(R^{n+m})$ is dense in $\phpf,$ Theorem 6.1 is proved.
\hfill {\bf Q.E.D.}

 We are now ready to prove the interpolation theorem on Hardy spaces $H^p_F$ for all $0<p<\infty$.

{\bf Proof of Theorem 1.20:} Suppose that $T$  is bounded from $H^{p_2}_F$ to $L^{p_2}$ and  from $H^{p_1}_F$ to  $L^{p_1}$. For any given $\lambda>0$ and $f\in H^p_F$, by the Calder\'on-Zygmund decomposition, $$f(x,y)=g(x,y)+b(x,y)$$ with
$$||g||^{p_1}_{H^{p_1}_F}\le C\lambda^{p_1-p}||f||_{H^p_F}^p\,\,\,\text{and}\,\, ||b||_{H^{p_2}_F}^{p_2}\le C\lambda^{p_2-p}||f||_{H^p_F}^p.$$

Moreover, we have proved the estimates
$$||g||^{p_1}_{H^{p_1}_F}\le C\int_{S(f)(x,y)\le \alpha}S(f)^{p_1}(x,y)dxdy$$ and
 $$||b||^{p_2}_{H^{p_2}_F}\le C\int_{S(f)(x,y)> \alpha}S(f)^{p_2}(x,y)dxdy$$
which implies that
$$
\align
||Tf||^p_p= &  p\int_0^\infty \alpha^{p-1} |\left\{(x,y): |Tf(x,y)|>\lambda\right\}|d \alpha\\
\le & p\int_0^\infty \alpha^{p-1}|\left\{(x,y): |Tg(x,y)|>\frac{\lambda}{2}\right\}|d\alpha+p\int_0^\infty \alpha^{p-1}|\left\{(x,y): |Tb(x,y)|>\frac{\lambda}{2}\right\}|d\alpha\\
\le &
p\int_0^\infty \alpha^{p-1}\int_{S(f)(x,y)\le \alpha}S(f)^{p_1}(x,y)dxdy d\alpha+p\int_0^\infty \alpha^{p-1}\int_{S(f)(x,y)>\alpha}S(f)^{p_2}(x,y)dxdy d\alpha\\ \le &
C||f||^p_{H^p_F}
\endalign
$$
Thus,
$$||Tf||_p\le C||f||_{H^p_F}$$
 for any $p_2<p<p_1$. Hence,  $T$ is bounded from $H^p_F$ to $L^p$.

To prove the second assertion that $T$ is bounded on $H^p_F$ for $p_2<p<p_1$,
 for any given $\lambda>0$ and $f\in H^p_F$, by the Calder\'on-Zygmund decomposition again
$$
\align
& |\left\{(x,y): |g(Tf)(x,y)|>\alpha\right\}|\\
\le & |\left\{(x,y): |g(Tg)(x,y)|>\frac{\alpha}{2}\right\}|+|\left\{(x,y): |g(Tb)(x,y)|>\frac{\alpha}{2}\right\}|\\
\le & C\alpha^{-p_1}||Tg||_{H^{p_1}_F}^{p_1}+ C\alpha^{-p_2}||Tb||_{H^{p_2}_F}^{p_2}\\
\le &  C\alpha^{-p_1}||g||^{p_1}_{H^{p_1}_F}+ C\alpha^{-p_2}||b||^{p_2}_{H^{p_2}_F}\\
\le  &  C \alpha^{-p_1}\int_{S(f)(x,y)\le \alpha}S(f)^{p_1}(x,y)dxdy +C\alpha^{-p_2}\int_{S(f)(x,y)> \alpha}S(f)^{p_2}(x,y)dxdy
\endalign
$$
which, as above, shows that
$||Tf||_{H^p_F}\le C || g(TF)||_p\le C||f||_{H^p_F}$
 for any $p_2<p<p_1$. \hfill{\bf Q.E.D.}

\vskip 0.1cm
\noindent {\bf References}

\noindent [Car] L. Carleson, A counterexample for measures bounded on $H^p$ for the bidisc, Mittag-Leffler Report No. 7, 1974.
\vskip -0.2cm
 \noindent [CW] R. R. Coifman and G. Weiss, Transference methods in analysis,
The Conference Board of the Mathematical Sciences by the AMS, 1977.

\noindent [Ch] S-Y. A. Chang, Carleson measures on the bi-disc, Ann. of Math. 109(1979), 613-620.

\noindent [CF1] S-Y. A. Chang and R. Fefferman, Some recent developments in Fourier analysis
and $H^p$ theory on product domains, Bull. A. M. S. 12(1985), 1-43.

\noindent [CF2] S-Y. A. Chang and R. Fefferman, The Calder\'on-Zygmund decomposition
on product domains, Amer. J. Math. 104(1982), 455-468.

\noindent [CF3] S-Y. A. Chang and R. Fefferman, A continuous version of duality of $H^1$ with
 $BMO$ on the bidisc, Annals of math. 112(1980), 179-201.

\noindent [F] R. Fefferman, Harmonic Analysis on product spaces, Ann. of math. 126(1987), 109-130.

\noindent [FS] C. Fefferman and E.M. Stein,  Some maximal inequalities, Amer. J. Math. 93 1971, 107--115.

\noindent [FS1] R. Fefferman and E. M.
Stein, singular integrals on product spaces, Adv. Math. 45(1982),
 117-143.

 \noindent [FJ] M. Frazier and B. Jawerth, A discrete transform and decompositions of distribution
 spaces, J. of Func. Anal. 93(1990), 34-170.

\noindent [FP] R. Fefferman and J. Pipher, Multiparameter operators and sharp weighted
 inequalities, Amer. J. Math., vol 11(1997), 337-369.

\noindent [GS] R. Gundy and E. M. Stein, $H^p$ theory for the polydisk, Proc. Nat.
Acad. Sci., 76(1979),

\noindent [HL] Y. Han and G. Lu, Endpoint estimates for singular integral operators and Hardy spaces associated with Zygmund dilations, to appear.

\noindent [JMZ] B. Jessen, J. Marcinkiewicz, and A. Zygmund, Note on the differentiability of
multiple integrals, Funda. Math. 25(1935), 217-234.

\noindent [J1] J. L. Journe, Calder\'on-Zygmund operators on product spaces, Rev. Mat.
Iberoamericana 1(1985), 55-92.

\noindent [J2] J. L. Journe, Two problems of Calder\'on-Zygmund theory on product spaces, Ann. Inst. Fourier (Grenoble) 38(1988), 111-132.

\noindent [MRS] D. Muller, F. Ricci, and E. M. Stein, Marcinkiewicz multipliers and multi-parameter
 structure on Heisenberg(-type) groups, I, Invent. math. 119(1995), 119-233.

\noindent [NRS] A. Nagel, F. Ricci, and E. M. Stein, Singular integrals with flag kernels and analysis
 on quadratic CR manifolds, J. Func. Anal. 181(2001), 29-118.

\noindent [P] J. Pipher, Journe's covering lemma and its extension to higher dimensions, Duke
 Mathematical Journal, 53(1986), 683-690.

\noindent [RS] F. Ricci and E. M. Stein, Multiparameter singular integrals and maximal functions,
 Ann. Inst. Fourier(Grenoble) 42(1992), 637-670.

\bye

\enddocument